\documentclass[12pt]{amsart}%
\usepackage{amsfonts}
\usepackage{amsmath}
\usepackage{amssymb}
\usepackage{graphicx}%

\usepackage{setspace}

\def\cJ{{\mathcal{J}}}
\def\cD{{\mathcal{D}}}
\def\cE{{\mathcal{E}}}
\def\cF{{\mathcal{F}}}
\def\cA{{\mathcal{A}}}
\def\cM{{\mathcal{M}}}
\def\cL{{\mathcal{L}}}
\def\cQ{{\mathcal{Q}}}
\def\PP{{\mathcal{P}}}
\def\cS{{\mathcal{S}}}
\def\cR{{\mathcal{R}}}
\def\cU{{\mathcal{U}}}
\def\bC{{\mathbf{C}}}

\def\bN{{\mathbb{N}}}
\def\bR{{\mathbb{R}}}
\def\mbY{{\mathbb{Y}}}
\def\mbS{{\mathbb{S}}}

\def\mbP{{\mathbb{P}}}
\def\mbF{{\mathbb{F}}}
\def\mbE{{\mathbb{E}}}
\def\mbD{{\mathbb{D}}}
\def\mbW{{\mathbb{W}}}

\def\f{\varphi}
\def\a{\alpha}
\def\sg{\sigma}
\def\pd{\partial}
\def\ds{\displaystyle}

\def\be{\begin{equation}}
\def\ee{\end{equation}}
\def\ds{\displaystyle}
\def\ba{\begin{array}}
\def\ea{\end{array}}

\newtheorem{theorem}{Theorem}
\theoremstyle{plain}
\newtheorem{corollary}[theorem]{Corollary}
\newtheorem{definition}[theorem]{Definition}
\newtheorem{example}[theorem]{Example}
\newtheorem{proposition}[theorem]{Proposition}
\newtheorem{lemma}[theorem]{Lemma}
\newtheorem{remark}[theorem]{Remark}
\numberwithin{equation}{section} \numberwithin{theorem}{section}

\begin{document}

\small{To appear in the  volume in honor of A. N. Shiryev}

\title[Wiener Chaos for Stochastic Equations]
{Stochastic Differential Equations: A Wiener Chaos Approach}
\author{S. V. Lototsky}
\curraddr[S. V. Lototsky]{Department of Mathematics, USC\\
Los Angeles, CA 90089} \email[S. V.
Lototsky]{lototsky@math.usc.edu}
\urladdr{http://math.usc.edu/$\sim$lototsky}
\author{B. L. Rozovskii}
\curraddr[B. L. Rozovskii]{Department of Mathematics, USC\\
Los Angeles, CA 90089} \email[B. L.
Rozovskii]{rozovski@math.usc.edu}
\urladdr{http://www.usc.edu/dept/LAS/CAMS/usr/facmemb/boris/main.htm}
\thanks{The work of S. V. Lototsky was partially supported by the Sloan
 Research Fellowship, by the NSF
CAREER award DMS-0237724, and by the ARO Grant DAAD19-02-1-0374.}
\thanks{The work of B. L. Rozovskii was partially supported by the ARO
Grant DAAD19-02-1-0374 and ONR Grant N0014-03-1-0027.}
\subjclass[2000]{Primary 60H15; Secondary 35R60, 60H40}
\keywords{Anticipating Equations, Generalized Random Elements,
Degenerate Parabolic Equations, Malliavin Calculus, Passive Scalar
Equation, Skorokhod Integral, S-transform, Weighted Spaces}

\begin{abstract}
A new method is described for constructing a generalized solution
for stochastic differential equations. The method is based on the
Cameron-Martin version of the Wiener Chaos expansion and provides
a unified framework for the study of ordinary and partial
differential equations driven by finite- or infinite-dimensional
noise with either adapted or anticipating input. Existence,
uniqueness, regularity, and probabilistic representation of this
Wiener Chaos solution is established for a large class of
equations. A number of examples are presented to illustrate the
general constructions. A detailed analysis is presented for the
various forms of the passive scalar equation and for the
first-order It\^{o} stochastic partial differential equation.
Applications to nonlinear filtering if diffusion processes and to
the stochastic Navier-Stokes equation are also discussed.
\end{abstract}
\maketitle

\begin{spacing}{0.15}

\noindent{\large\bf Contents}

  {\bf\ref{sec-:I}.
Introduction\dotfill\pageref{sec-:I}}

  {   \bf\ref{sec:.:clas}.
Traditional Solutions of Linear Parabolic
Equations\dotfill\pageref{sec:.:clas}}

   {   \bf\ref{sec:.:WN}. White
Noise Solutions of Stochastic Parabolic
Equations\dotfill\pageref{sec:.:WN}}

  {
  \bf\ref{secGF}. Generalized Functions on the Wiener
Chaos Space\dotfill\pageref{secGF}}

  {
  \bf\ref{secMD}. The
Malliavin Derivative and its Adjoint\dotfill\pageref{secMD}}

  {
  \bf\ref{secGS}. The Wiener Chaos Solution and the
Propagator\dotfill\pageref{secGS}}

  {
  \bf\ref{secWS}. Weighted Wiener Chaos Spaces and S-Transform
\dotfill\pageref{secWS}}

 {   \bf\ref{secPWC}. General Properties of the
Wiener Chaos Solutions \dotfill\pageref{secPWC}}

  {
  \bf\ref{secLH}. Regularity of the Wiener Chaos  Solution
\dotfill\pageref{secLH}}

  {
  \bf\ref{secPR}. Probabilistic Representation of Wiener Chaos
Solutions
\dotfill\pageref{secPR}}

  {   \bf\ref{sec:.:FLT}. Wiener
Chaos and Nonlinear Filtering \dotfill\pageref{sec:.:FLT}}

  {
  \bf\ref{secPS}. Passive Scalar in a Gaussian Field
\dotfill\pageref{secPS}}

  {   \bf\ref{sec:.:SNS}.
Stochastic Navier-Stokes Equation \dotfill\pageref{sec:.:SNS}}

  {   \bf\ref{sec:.:frd}.
First-Order It\^{o} Equations \dotfill\pageref{sec:.:frd}}

\end{spacing}

\section{Introduction}
\label{sec-:I}

Consider a stochastic evolution equation
\begin{equation}
\label{eq:intr}du(t)=({\mathcal{A}} u(t)+f(t))dt+ ({\mathcal{M}}
u(t) + g(t))dW(t),
\end{equation}
where ${\mathcal{A}}$ and ${\mathcal{M}}$ are differential
operators, and $W$ is a noise process on a stochastic basis
${\mathbb{F}}=(\Omega, {\mathcal{F}},
\{{\mathcal{F}}_{t}\}_{t\geq0}, {\mathbb{P}})$. Traditionally,
this equation is studied under the following assumptions:

\begin{enumerate}
\item[(i)] The operator ${\mathcal{A}}$ is elliptic, the order of
the operator ${\mathcal{M}}$ is at most half the order of
${\mathcal{A}}$, and a special parabolicity condition holds.

\item[(ii)] The functions $f$ and $g$ are predictable with respect
to the filtration $\{{\mathcal{F}}_{t}\}_{t\geq0}$, and the
initial condition is ${\mathcal{F}}_{0}$-measurable.

\item[(iii)] The noise process $W$ is sufficiently regular.
\end{enumerate}

Under these assumptions, there exists a unique predictable
solution $u$ of (\ref{eq:intr}) so that $u \in
L_{2}(\Omega\times(0,T); H)$ for $T>0$ and a suitable function
space $H$ (see, for example, Chapter 3 of \cite{Roz}). Moreover,
there are examples showing that the parabolicity condition and the
regularity of noise are necessary to have a square integrable
solution of (\ref{eq:intr}).

The objective of the current paper is to study stochastic
differential equations of the type (\ref{eq:intr}) without making
the above assumptions (i)--(iii). We show that, with a suitable
definition of the solution, solvability of the stochastic equation
is essentially equivalent to solvability of a deterministic
evolution equation $dv=({\mathcal{A}} v + \varphi)dt$ for certain
functions $\varphi$; the operator ${\mathcal{A}}$ does not even
have to be elliptic.

Generalized solutions have been introduced and studied for
stochastic differential equations, both ordinary and with partial
derivatives, and definitions of such solutions relied on various
forms of the Wiener Chaos decomposition. For stochastic ordinary
differential equations, Krylov and Veretennikov \cite{KV} used
multiple Wiener integral expansion to study Ito diffusions with
non-smooth coefficients, and more recently, LeJan and Raimond
\cite{LeJan} used a similar approach in the construction of
stochastic flows. Various versions of the Wiener chaos appear in a
number of papers on nonlinear filtering and related topics
\cite[etc.]{BdKl,LMR,MR,Ocn,Wong} The book by Holden et al.
\cite{HOUZ} presents a systematic approach to the stochastic
differential equations based on the white noise theory. See also
\cite{HKPS}, \cite{Pot} and the references therein.

For stochastic partial differential equations, most existing
constructions of the generalized solution rely on various
modifications of the Fourier
transform in the infinite-dimensional Wiener Chaos space $L_{2}({\mathbb{W}%
})=L_{2}(\Omega,{\mathcal{F}}^{W}_{T},{\mathbb{P}})$. The two main
modifications are known as the S-transform \cite{HKPS} and the
Hermite transform \cite{HOUZ}. The key elements in the development
of the theory are the spaces of the test functions and the
corresponding distributions. Several constructions of these spaces
were suggested by Hida \cite{HKPS}, Kondratiev \cite{Kond}, and
Nualart and Rozovskii \cite{NR}. Both S- and Hermite transforms
establish a bijection between the space of generalized random
elements and a suitable space of analytic functions. Using the
S-transform, Mikulevicius and Rozovskii \cite{MR} studied
stochastic parabolic equations with non-smooth coefficients, while
Nualart and Rozovskii \cite{NR} and Potthoff et. al \cite{Pot}
constructed generalized solutions for the equations driven by
space-time white noise in more than one spacial dimension. Many
other types of equations have been studied, and the book
\cite{HOUZ} provides a good overview of literature the
corresponding results.

In this paper, generalized solutions of (\ref{eq:intr}) are
defined in the spaces that are even larger than Hida or Kondratiev
distribution. The Wiener Chaos space is a separable Hilbert space
with a Cameron-Martin basis \cite{CM}. The elements of the space
with a finite Fourier series expansion
provide the natural collection of test functions ${\mathcal{D}}(L_{2}%
({\mathbb{W}}))$, an analog of the space
${\mathcal{D}}({\mathbb{R}}^{d})$ of smooth compactly supported
functions on ${\mathbb{R}}^{d}$. The corresponding space of
distributions ${\mathcal{D}}^{\prime}(L_{2}({\mathbb{W}}))$ is the
collection of generalized random elements represented by formal
Fourier series. A generalized solution $u=u(t,x)$ of
(\ref{eq:intr}) is constructed as an element of
${\mathcal{D}}^{\prime}(L_{2}({\mathbb{W}}))$ so that the
generalized Fourier coefficients satisfy a system of deterministic
evolution equations, known as the {\ propagator}. If the equation
is linear the propagator is a lower-triangular system. We call
this solution a Wiener Chaos solution.

The propagator was first introduced by Mikulevicius and Rozovskii
in \cite{MR0}, and further studied in \cite{LMR}, as a numerical
tool for solving the nonlinear filtering problem. The  propagator
can also be derived for certain nonlinear equations; in
particular, it was used in \cite{MR_AP,MR1,MR2} to study the
stochastic Navier-Stokes equation.

The propagator approach to defining the solution of
(\ref{eq:intr}) has two advantages over the S-transform approach.
First, the resulting construction is more general: there are
equations for which the Wiener Chaos solution is not in the domain
of the S-transform. Indeed, it is shown in Section \ref{sec:.:frd}
that, for certain initial conditions, equation $du=u_{x}dW_{t}$
has a Wiener Chaos solution for which the S-transform is not
defined. On the other hand, by Theorem \ref{th:SW1} below, if the
generalized solution of (\ref{eq:intr}) can be defined using the
S-transform, then this solution is also a Wiener Chaos solution.
Second, there is no problem of inversion: the propagator provides
a direct approach to studying the properties of Wiener Chaos
solution and computing both the sample trajectories and
statistical moments.

Let us emphasize also the following important features of the
Wiener Chaos approach:

\begin{itemize}
\item The Wiener Chaos solution  is a strong solution in the
probabilistic sense, that is, it is uniquely determined by the
coefficients, free terms, initial condition, and the Wiener
process.

\item The solution exists under minimal regularity conditions on
the coefficients in the stochastic part of the equation and no
special measurability restriction on the input.

\item The Wiener Chaos solution  often serves as a convenient
first step in the investigation of the traditional solutions or
solutions in weighted stochastic Sobolev spaces that are much
smaller then the spaces of Hida or Kondratiev distributions.
\end{itemize}

To better understand the connection between the Wiener Chaos
solution and other notions of the solution, recall that,
traditionally, by a  solution of a stochastic equation we
understand
 a random process or field satisfying  the equation
 for almost all elementary outcomes. This  solution can be  either strong or weak
 in the probabilistic sense.

Probabilistically strong  solution is constructed on a prescribed
probability space  with a specific  noise process. Existence of
strong  solutions requires certain regularity of the coefficients
and the noise in the equation. The tools for constructing strong
solutions often come from the
 theory of the corresponding deterministic equations.

 Probabilistically  weak solution includes not only
the solution process but also the  stochastic basis and the noise
process. This freedom  to choose the probability space and the
noise process makes the conditions for existence of weak solutions
less restrictive than the similar conditions  for strong
solutions.  Weak solutions can be obtained  either by considering
the corresponding martingale problem or by constructing a suitable
Hunt process using the theory of the  Dirichlet forms.

There exist  equations  that have neither weak nor strong
solutions in the traditional sense. An example is the bi-linear
stochastic heat equation driven by a
 multiplicative space-time white noise in two or more spatial dimensions:
 the irregular nature of the noise  prevents the existence of a random
 field  that would satisfy the equation for individual elementary outcomes.
 For such equations, the solution must be  defined as a generalized random element
 satisfying the equation after the randomness has been averaged out.

 White noise theory provides one
 approach for constructing these generalized solutions. The approach is
 similar to the Fourier integral method for deterministic equations.
  The white noise solution is  constructed on a special white noise
 probability space by inverting an integral transform;
  the special structure of the probability space is
 essential  to
 carry out the inversion. We can therefore say that the  white noise solution
  extends the notion  of the probabilistically weak solution.
  Still, this extension is not a true generalization:
   when the equation satisfies the necessary regularity conditions,
   the connection between the white noise and the traditional weak solution
   is often not clear.

The Wiener chaos approach provides the means for constructing a
  generalized solution on  a prescribed probability space.
  The Wiener Chaos solution is a formal Fourier
  series in the corresponding Cameron-Martin basis.
   The coefficients in the series are uniquely determined by the
  equation via the propagator system. This representation
   provides a convenient way for computing numerically the solution and
  its statistical moments. As a result,  the Wiener Chaos solution
  extends the notion  of the probabilistically strong  solution.
  Unlike the white noise approach, this is a bona fide  extension: when the
  equation satisfies the necessary regularity conditions, the Wiener Chaos
  solution coincides with the traditional strong solution.

After the general discussion of the Wiener Chaos space in Sections
\ref{secGF} and \ref{secMD}, the Wiener Chaos solution for
equation (\ref{eq:intr}) and the main properties of the solution
are studied in Section \ref{secGS}. Several examples illustrate
how the Wiener Chaos solution provides a uniform treatment of
various types of equations: traditional parabolic, non-parabolic,
and anticipating. In particular, for equations with
non-predictable input, the Wiener Chaos solution corresponds to
the Skorohod integral interpretation of
the equation. The initial solution space ${\mathcal{D}}^{\prime}({\mathbb{W}%
})$ is too large to provide much\textit{ of} interesting
information about the solution. Accordingly, Section \ref{secWS}
discusses various weighted Wiener Chaos spaces. These weighted
spaces provide the necessary connection between the Wiener Chaos,
white noise, and traditional solutions. This connection is studied
in Section \ref{secPWC}. In Section \ref{secLH}, the Wiener Chaos
solution is constructed for degenerate linear parabolic equations
and new regularity results are obtained for the solution.
Probabilistic representation of the Wiener Chaos solution is
studied in Section \ref{secPR}, where a Feynmann-Kac type formula
is derived. Sections \ref{sec:.:FLT}, \ref{secPS},
\ref{sec:.:SNS}, and \ref{sec:.:frd} discuss the applications of
the general results to particular equations: the Zakai filtering
equation, the stochastic transport equation, the stochastic
Navier-Stokes equation, and a first-order It\^{o} SPDE.

The following notation will be in force throughout the paper:
$\Delta$ is the Laplace operator, $D_{i}=\partial/\partial x_{i}$,
$i=1,\ldots,d$, and summation over the repeated indices is
assumed. The space of continuous functions is denoted by
${\mathbf{C}}$, and $H_{2}^{\gamma}$, $\gamma \in{\mathbb{R}}$, is
the Sobolev space
\[
\left\{  f:\int_{{\mathbb{R}}}|\hat{f}(y)|^{2}(1+|y|^{2})^{\gamma}%
dy<\infty\right\}  ,\ \mathrm{where\ }\hat{f}%
\mathrm{\ is\ the\ Fourier\ transform\ of\ }f.
\]

\section{Traditional Solutions of Linear Parabolic Equations}
\label{sec:.:clas} \setcounter{equation}{0}

Below is a summary of the Hilbert space theory of linear
stochastic parabolic equations. The  details can be found in the
books \cite{DaPr1} and \cite{Roz}; see also \cite{Kr_Lp1}.  For a
Hilbert space $X$, $(\cdot, \cdot)_X$ and $\|\cdot\|_X$ denote the
inner product and the norm in $X$.

\begin{definition}
\label{def:Ntrip}
 The triple $(V, H, V')$ of Hilbert spaces is called normal if and only
if

\begin{enumerate}
\item $V\hookrightarrow H \hookrightarrow V'$ and both  embeddings
$V\hookrightarrow H $ and $H \hookrightarrow V'$ are dense and
continuous; \item The space $V'$ is the dual of $V$ relative to
the inner product in $H$; \item There exists a constant $C>0$ so
that $|(h, v)_H| \leq C\|v\|_V\|h\|_{V'}$ for all $v\in V$ and
$h\in H$.
\end{enumerate}
\end{definition}
For example,  the Sobolev spaces $(H^{\ell + \gamma}_2(\bR^d),
H^{\ell}_2(\bR^d), H^{\ell - \gamma}_2(\bR^d))$, $\gamma >0$,
$\ell \in \bR$,
 form a normal triple.

 Denote by $\langle v', v\rangle$, $v'\in V'$, $v\in V$,
  the duality between $V$ and $V'$ relative to
 the inner product in $H$. The properties of the normal triple imply that
 $|\langle v', v\rangle|\leq C\|v\|_V\|v'\|_{V'}$, and,
 if $v'\in H$ and $v\in V$, then $\langle v', v\rangle = (v',v)_H;$

Let $\mbF=(\Omega, \cF, \{\cF_t\}_{t\geq 0}, \mbP)$ be a
stochastic basis with the usual assumptions. In particular, the
sigma-algebras $\cF$  and $\cF_0$ are $\mbP$-complete, and the
filtration  $\{\cF_t\}_{t\geq 0}$ is right-continuous; for
details, see \cite[Definition I.1.1]{LSh.m}.
 We assume that $\mbF$ is rich enough to carry  a collection
 $w_k=w_k(t),\; k\geq 1,\; t\geq 0$ of
independent standard Wiener processes.

Given a normal triple $(V,H,V')$ and a family of  linear bounded
operators $\cA(t): V \to V'$, $\cM_k(t): V\to H$, $t\in [0,T]$,
consider the following equation:
\begin{equation}
\label{eq:ParabCl} u(t)=u_0+\int_0^t (\cA u(s)+f(s))ds + \int_0^t
(\cM_k u (s)+g_k(s))
 dw_k(s),\ 0\leq t\leq T,
\end{equation}
where $T<\infty$ is fixed and non-random and the  summation
convention  is in force.

Assume that,  for all $v \in V$,
\begin{equation}\label{eq:MCl}
\sum_{k\geq 1}\|\cM_k (t)v\|_H^2 < \infty, \ t\in [0,T].
\end{equation}
The
input data $u_0, f$, and $g_k$ are chosen so that
\begin{equation}
\label{eq:InputCl} \mbE\left(\|u_0\|_H^2 + \int_0^T\|f(t)\|_{V'}^2
dt + \sum_{k\geq 1} \int_0^T\|g_k(t)\|_H^2dt\right) < \infty,
\end{equation}
$u_0$ is $\cF_0$-measurable, and the processes $f,g_k$ are
$\cF_t$-adapted, that is, $f(t)$ and each $g_k(t)$ are
$\cF_t$-measurable for each $t\geq 0$.

\begin{definition}
\label{def:ParabCl} An $\cF_t$-adapted
 process $u\in L_2(\mbF;L_2((0,T); V))$ is called a traditional, or
 square-integrable,  solution of equation
(\ref{eq:ParabCl}) if, for every $v\in V$, there exists  a
measurable sub-set $\Omega'$ of $\Omega$ with $\mbP(\Omega')=1$,
so that, the equality
\begin{equation}\label{eq:trad-def}
(u(t),v)_H=(u_0,v)_H + \int_0^t \langle \cA u(s)+ f(s), v \rangle
ds + \sum_{k\geq 1} (\cM_k u(s)+g_k(s),v)_Hdw_k(s)
\end{equation}
holds on $\Omega'$ for all $0\leq t\leq T$.
\end{definition}

Existence and uniqueness of the traditional solution for
(\ref{eq:ParabCl})
 can be  established when the equation is parabolic.

 \begin{definition}
\label{def:parab} Equation (\ref{eq:ParabCl}) is called
\textbf{strongly parabolic} if there exists a positive number
$\varepsilon$ and a real number $C_{0}$ so that, for all $v\in V$
and $t\in [0,T]$,
\begin{equation}
2\langle{\mathcal{A}}(t)v,v\rangle+\sum_{k\geq1}\Vert{\mathcal{M}}(t)_{k}v\Vert
_{H}^{2}+\varepsilon\Vert v\Vert_{V}^{2}\leq C_{0}\Vert v\Vert_{H}%
^{2}.\label{eq:ParabCll}%
\end{equation}
Equation (\ref{eq:ParabCl}) is called \textbf{weakly} \textbf{\
parabolic} (or degenerate parabolic) if condition
(\ref{eq:ParabCll}) holds with $\varepsilon=0$.
\end{definition}

\begin{theorem}
\label{th:ParabCl} If (\ref{eq:InputCl}) and (\ref{eq:ParabCll})
hold, then there exists a unique traditional solution of
(\ref{eq:ParabCl}). The  solution process $u$ is an element of the
space
$$
L_2(\mbF; L_2((0,T); V))\bigcap L_2(\mbF; \bC((0,T), H))
$$ and satisfies
\begin{equation}\label{eq:enestCL}
\begin{split}
&\mbE\left( \sup_{0<t<T} \|u(t)\|_H^2 + \int_0^T\|u(t)\|_V^2dt \right)\\
&\leq C(C_0, \delta,T) \mbE\left(\|u_0\|_H^2 +
\int_0^T\|f(t)\|_{V'}^2 dt + \sum_{k\geq 1}
\int_0^T\|g_k(t)\|_H^2dt\right).
\end{split}
\end{equation}
\end{theorem}

\begin{proof} This follows, for example, from Theorem 3.1.4 in
\cite{Roz}.
\end{proof}

A somewhat different solvability result holds for weakly parabolic
equations \cite[Section 3.2]{Roz}.

As an application of Theorem \ref{th:ParabCl}, consider equation
\begin{equation}
\label{eq:ParabEqCl}
\begin{split}
du(t,x)&= (a_{ij}(t,x)D_iD_ju(t,x) + b_i(t,x)D_iu(t,x) +c(t,x)u(t,x)+f(t,x))dt\\
&+ (\sigma_{ik}(t,x)D_iu(t,x)+\nu_k(t,x)u(t,x)+g_k(t,x))dw_k(t)
\end{split}
\end{equation}
with $0<t\leq T,\ x \in \bR^d,$ and
 initial condition $u(0,x)=u_0(x)$. Assume that
\begin{enumerate}
\item[(CL1)] The functions $a_{ij}$ are bounded and Lipschitz
continuous, the functions $b_i$, $c$, $\sigma_{ik}$, and $\nu$ are
bounded measurable. \item[(CL2)] There exists a positive number
$\varepsilon>0$ so that
$$
(2a_{ij}(x)-\sigma_{ik}(x)\sigma_{jk}(x))y_iy_j \geq \varepsilon
|y|^2,
 \ x,y \in  \bR^d,\ t\in [0,T].
$$
\item[(CL3)] There exists a positive number $K$ so that, for all
$x \in \bR^d$, $\sum_{k\geq 1}|\nu_k(x)|^2 \leq K.$ \item[(CL4)]
The initial condition $u_0\in L_2(\Omega; L_2(\bR^d))$ is
$\cF_0$-measurable,  the processes $f\in L_2(\Omega\times[0,T];
H^{-1}_2(\bR^d))$ and $g_k\in L_2(\Omega\times[0,T]; L_2(\bR^d))$
are $\cF_t$-adapted, and $\sum_{k\geq 1}\int_0^T\mbE
\|g_k\|_{L_2(\bR^d)}^2(t)dt < \infty$.
\end{enumerate}
\begin{theorem}
\label{th:ParabCl1} Under assumptions (CL1)--(CL4), equation
(\ref{eq:ParabEqCl}) has a unique traditional solution
$$
u \in L_2(\mbF;L_2((0,T); H^1_2(\bR^d)))\bigcap L_2(\mbF;
\bC((0,T), L_2(\bR^d))),
$$ and the solution satisfies
\begin{equation}\begin{split}
&\mbE\left( \sup_{0<t<T} \|u\|_{L_2(\bR^d)}^2(t) +
\int_0^T\|u\|_{H^1_2(\bR^d)}^2(t)dt \right)\\
&\leq C(K, \varepsilon,T) \mbE\left(\|u_0\|_{L_2(\bR^d)}^2 +
\int_0^T\!\! \|f\|_{H^{-1}_2(\bR^d)}^2(t) dt \!+\! \sum_{k\geq 1}
\int_0^T\!\! \|g_k\|_{L_2(\bR^d)}^2(t)dt\right).
\end{split}
\end{equation}
\end{theorem}
\begin{proof} Apply Theorem \ref{th:ParabCl}
in the normal triple \\
$(H^1_2(\bR^d),L_2(\bR^d),H^{-1}_2(\bR^d))$; condition
(\ref{eq:ParabCll}) in this case is equivalent to
 assumption (CL2). The
details of the proof are  in \cite[Section 4.1]{Roz}.
\end{proof}

Condition (\ref{eq:ParabCll}) essentially means that the
deterministic part of the equation dominates the stochastic part.
Accordingly, there are two main ways to violate
(\ref{eq:ParabCll}):
\begin{enumerate}
\item The order of the operator $\cM$ is more than half the order
of the operator $\cA$. Equation $du=u_xdw(t)$ is an example. \item
The value of  $\sum_k\|\cM_k(t)v\|_H^2$ is too large. This value
can be either finite, as in equation
$du(t,x)=u_{xx}(t,x)dt + 5u_x(t,x)dw(t)$
or infinite, as in equation
\begin{equation}
\label{eq:stwn} du(t,x)=\Delta u(t,x)dt+\sigma_k(x)udw_k, \
\sigma_k - {\rm \ CONS\ in\  } L_2(\bR^d),\ d\geq 2.
\end{equation}
 Indeed,  it is shown in \cite{NR} that,
for equation (\ref{eq:stwn}), we have
$$
\sum_{k\geq
1}\|\cM_k(t)v\|_H^2=\infty
$$
 in every Sobolev space $H^{\gamma}$.
\end{enumerate}

Without condition (\ref{eq:ParabCll}), analysis of equation
(\ref{eq:ParabCl}) requires new technical tools  and a different
notion of solution. The white noise theory provides one possible
collection  of such tools.

\section{White Noise Solutions of Stochastic Parabolic Equations}
\label{sec:.:WN} \setcounter{equation}{0} \setcounter{theorem}{0}

The central part of the white noise theory is  the mathematical
model for the
 derivative of the Brownian motion. In particular,
 the It\^{o} integral $\int_0^t f(s) dw(s)$
is replaced with the integral  $\int_0^tf(s)\diamond\dot{W}(s)ds$,
where $\dot{W}$ is the white noise process and $\diamond$ is the
Wick product.  The white noise  formulation is very different from
the Hilbert space approach of the previous section, and requires
several new   constructions. The book \cite{HKPS} is a general
reference about the white noise theory, while  \cite{HOUZ}
presents
 the white noise analysis of stochastic partial
differential equations. Below is the summary of the main
definitions and results.

Denote by $\cS=\cS(\bR^{\ell})$ the Schwartz space of rapidly
decreasing functions and by $\cS'=\cS'(\bR^{\ell})$, the Schwartz
space of tempered distributions.
 For the properties of the spaces $\cS$ and $\cS'$ see \cite{Rudin}.

\begin{definition}
\label{def:WNSpace} The white noise probability space is the
triple
$$
\mbS=(\cS', {\mathcal{B}}(\cS'), \mu),
$$
 where
${\mathcal{B}}(\cS')$ is the Borel sigma-algebra of subsets of
$\cS'$, and $\mu$ is the normalized Gaussian measure on
${\mathcal{B}}(\cS')$.
\end{definition}
The measure $\mu$ is characterized by the property
$$
\int_{\cS'}e^{\sqrt{-1}\langle \omega,
\varphi\rangle}d\mu(\omega)=
e^{-\frac{1}{2}\|\varphi\|_{L_2(\bR^d)}^2},
$$
where $\langle \omega, \varphi \rangle$, $\omega \in \cS'$,
$\varphi \in \cS$, is the duality between $\cS$ and $\cS'$.
Existence of this measure   follows from the Bochner-Minlos
theorem \cite[Appendix A]{HOUZ}.

Let $\{\eta_k, k\geq 1\}$ be the Hermite basis in
$L_2(\bR^{\ell})$, consisting of  the normalized eigenfunctions
of the operator
\begin{equation}
\label{eq:HermOp} \Lambda= -\Delta +|x|^2, \ x\in \bR^{\ell}.
\end{equation}
Each $\eta_k$ is an element of $\cS$ \cite[Section 2.2]{HOUZ}.

Consider the collection of multi-indices
 $$
 \cJ_1=\Big\{ \alpha =(\alpha_i,\ i\geq 1),\ \alpha_{i}\in
\{0,1,2,\ldots\},\ \sum_{i} \alpha_i<\infty \Big\}.
$$
The set $\cJ_1$ is countable, and, for every $\alpha \in \cJ$,
only finitely many of $\alpha_i$ are not equal to zero. For
$\alpha \in {\mathcal J}_1$, write $\alpha!=\prod_{i} \alpha_i!$
and define
\begin{equation}
\label{eq:CMBWN} \xi_{\alpha}(\omega)=\frac{1}{\sqrt{\alpha!}}
\prod_{i}H_{\alpha_{i}}(\langle \omega, \eta_{i}\rangle),\ \omega
\in \cS',
\end{equation}
where $\langle \cdot, \cdot \rangle$ is the duality between $\cS$
and $\cS'$, and
\begin{equation}
\label{eq:HerPol}
 H_{n}(t)=(-1)^ne^{t^{2}/2}
\frac{d^{n}}{dt^{n}}e^{-t^{2}/2}
\end{equation}
is $n^{\rm th}$
 Hermite polynomial. In particular,  $H_1(t)=1$, $H_1(t)=t$, $H_2(t)=t^2-1$.
If, for example,   $\alpha=(0,2,0,1,3,0,0,\ldots)$ has   three
non-zero entries, then
$$
\xi_{\alpha}(\omega)=\frac{H_2(\langle \omega,
\eta_2\rangle)}{2!}\cdot \langle \omega, \eta_4\rangle \cdot
\frac{H_3(\langle \omega, \eta_5\rangle)}{3!}.
$$

\begin{theorem}
\label{th:CMWN} The collection $\{\xi_{\alpha}, \; \alpha \in
\cJ_1\}$ is an orthonormal basis in $L_2(\mbS)$.
\end{theorem}

\begin{proof} This is a version of the classical result of Cameron
and Martin \cite{CM}. In this particular form, the result is
stated and proved in \cite[Theorem 2.2.3]{HOUZ}.
\end{proof}

By Theorem \ref{th:CMWN}, every element $\varphi$ of $L_2(\mbS)$
is represented as a Fourier series  $\varphi=\sum_{\alpha}
\varphi_{\alpha}\xi_{\alpha}$, where
$\varphi_{\alpha}=\int_{\cS'}\varphi(\omega)\xi_{\alpha}(\omega)d\mu$,
and $\|\varphi\|_{L_2(\mbS)}^2=\sum_{\alpha \in \cJ_1}
|\varphi_{\alpha}|^2$.

For $\alpha \in \cJ_1$ and $q\in \bR$, we write
$$
(2\bN)^{q\alpha}=\prod_j (2j)^{q\alpha_j}.
$$

\begin{definition}
\label{def:KS1} For $\rho\in [0,1]$ and $q\geq 0$,
\begin{enumerate}
\item the space $(\cS)_{\rho, q}$ is the collection of elements
$\varphi$ from $L_2(\mbS)$ so that
$$
\|\varphi\|_{\rho, q}^2=\sum_{\alpha \in \cJ_1}
(\alpha!)^{\rho}(2\bN)^{q\alpha} |\varphi_{\alpha}|^2 < \infty;
$$
\item the space $(\cS)_{-\rho, -q}$ is the closure of $L_2(\mbS)$
relative to the norm
\begin{equation}
\label{eq:KondWN3} \|\varphi\|_{-\rho, -q}^2=\sum_{\alpha \in
\cJ_1} (\alpha!)^{-\rho}(2\bN)^{-q\alpha} |\varphi_{\alpha}|^2;
\end{equation}
\item the space $(\cS)_{\rho}$ is the projective  limit of
$(\cS)_{\rho, q}$ as $q$ changes over all non-negative integers;
\item the space $(\cS)_{-\rho}$ is the inductive   limit of
$(\cS)_{-\rho,-q}$ as $q$ changes over all non-negative integers.
\end{enumerate}
\end{definition}

It follows that
\begin{itemize}
\item For each $\rho\in [0,1]$ and $q\geq 0$, $((\cS)_{\rho,
q},L_2(\mbS),  (\cS)_{-\rho, -q})$ is a normal triple of Hilbert
spaces. \item The space $(\cS)_{\rho}$ is a Frechet space with
topology generated by the countable family of norms
$\|\cdot\|_{\rho, n}$, $n=0,1,2,\ldots$, and $\varphi \in
(\cS)_{\rho}$ if and only if $\varphi \in (\cS)_{\rho, q}$ for
every $q \geq 0$. \item The space $(\cS)_{-\rho}$ is the dual of
$(\cS)_{\rho}$ and $\varphi \in (\cS)_{-\rho}$ if and only if
$\varphi \in (\cS)_{-\rho, -q}$ for some  $q \geq 0$. Every
element $\varphi $ from $(\cS)_{\rho}$ is identified with a formal
sum $\sum_{\alpha \in \cJ_1} \varphi_{\alpha} \xi_{\alpha}$ so
that (\ref{eq:KondWN3}) holds for some $q\geq 0$. \item For
$0<\rho<1$,
$$
(\cS)_1 \subset (\cS)_{\rho} \subset (\cS)_0 \subset L_2(\mbS)
 \subset (\cS)_{-0} \subset
(\cS)_{-\rho} \subset (\cS)_{-1},
$$
with all inclusions strict.
\end{itemize}

The spaces $(\cS)_{0}$ and $(\cS)_{1}$ are known   as the spaces
of Hida and Kondratiev test functions. The spaces $(\cS)_{-0}$ and
$(\cS)_{-1}$ are known  as the spaces of Hida and Kondratiev
distributions. Sometimes, the spaces  $(\cS)_{\rho}$ and
$(\cS)_{-\rho}$, $0<\rho\leq 1$, go under the name of Kondratiev
test functions and Kondratiev distributions, respectively.

Let $h \in \cS$ and $h_k=\int_{\bR^{\ell}} h(x)\eta_k(x)dx.$ Since
the asymptotics of $n^{\rm th}$ eigenvalue of the operator
$\Lambda$ in (\ref{eq:HermOp}) is $n^{1/d}$ \cite[Chapter
21]{HlPh} and $\Lambda^k h \in \cS$ for every positive integer
$k$, it follows that
\begin{equation}
\label{eq:StrWN1} \sum_{k\geq 1}|h_k|^2k^q <\infty
\end{equation}
for every $q\in \bR$.

For $\alpha \in \cJ_1$ and $h_k$ as above, write
$h^{\alpha}=\prod_j (h_j)^{\alpha_j}$, and define the stochastic
exponential
\begin{equation}
\label{eq:StrWN2} \cE(h)=\sum_{\alpha \in \cJ_1}
\frac{h^{\alpha}}{\sqrt{\alpha!}}\xi_{\alpha}
\end{equation}

\begin{lemma}
\label{lm:StochExp} The stochastic exponential $\cE=\cE(h)$, $h
\in \cS$,
 has the following properties:
\begin{itemize}
\item $\cE(h)\in (\cS)_{\rho}$, $0<\rho<1$; \item For every $q>0$,
there exists a $\delta>0$ so that $\cE(h)\in (\cS)_{1,q}$ as long
as $\sum_{k\geq 1}|h_k|^2< \delta$.
\end{itemize}
\end{lemma}

\begin{proof} Both properties are verified by direct calculation
\cite[Chapter 2]{HOUZ}. \end{proof}

\begin{definition} The S-transform $S\varphi(h)$ of an element
 $\varphi=\sum_{\alpha\in \cJ}\varphi_{\alpha}\xi_{\alpha}$
from $ (\cS)_{-\rho}$ is the number
\begin{equation}
\label{eq:StrWN3} S\varphi(h)=\sum_{\alpha \in \cJ_1}
\frac{h^{\alpha}}{\sqrt{\alpha!}}\varphi_{\alpha},
\end{equation}
where  $h=\sum_{k\geq 1}h_k\eta_k\in \cS$ and $h^{\alpha}=\prod_j
(h_j)^{\alpha_j}$.
\end{definition}

The definition implies that if $\varphi \in (\cS)_{-\rho, -q}$ for
some $q\geq 0$, then  $S\varphi(h)=\langle \varphi, \cE(h)
\rangle$, where $\langle \cdot, \cdot \rangle$ is the duality
between $(\cS)_{\rho,q}$ and $(\cS)_{-\rho,-q}$ for suitable $q$.
Therefore, if $\rho<1$, then
 $S\varphi(h)$ is
well-defined for all $h\in \cS$,  and, if $\rho=1$, the
$S\varphi(h)$ is well-defined for $h$ with sufficiently small
$L_2(\bR^{\ell})$ norm. To give a complete characterization of the
S-transform, one additional construction is necessary.

Let  $\cU^{\rho}$, $0\leq \rho<1$, be the collection of mappings
$F$ from $\cS$ to the complex numbers so that
\begin{enumerate}
\item[1.] For every $h_1, h_2 \in \cS$, the function $F(h_1+zh_2)$
is an analytic function of the complex variable $z$. \item[2.]
There exist positive numbers $K_1, K_2$ and an integer number $n$
so that, for all $h\in \cS$ and all complex number $z$,
$$
|F(zh)|\leq K_1
\exp\left(K_2\|\Lambda^nh\|_{L_2(\bR^d)}^{\frac{2}{1-\rho}}
|z|^{\frac{2}{1-\rho}}\right).
$$
\end{enumerate}
For $\rho=1$, let $\cU^1$ be the collection of mappings $F$ from
$\cS$ to the complex numbers so that
\begin{enumerate}
\item[1$'$.] There exist $\varepsilon>0$ and a positive  integer
$n$ so that, for all $h_1, h_2\in \cS$ with $\|\Lambda^n
h_1\|_{L_2(\bR^{\ell})}<\varepsilon$, the function of a complex
variable  $z\mapsto F(h_1+h_2z)$ is analytic at zero, and
\item[2$'$.] There exists a positive number $K$ so that, for all
$h\in \cS$ with $\|\Lambda^n h\|_{L_2(\bR^{\ell})}<\varepsilon$,
$|F(h)|\leq K$.
\end{enumerate}
Two mappings $F,G$ with properties $1'$ and $2'$ are identified
with the same element of $\cU^1$ if $F=G$ on an open neighborhood
of zero in $\cS$.

The following result holds.
\begin{theorem}
\label{th:StrWN} For every $\rho\in [0,1]$, the S-transform is a
bijection from $(\cS)_{-\rho}$ to $\cU^{\rho}$.
\end{theorem}
In other words, for every $\varphi\in (\cS)_{-\rho}$, the
S-transform $S\varphi$ is an element of $\cU^{\rho}$, and, for
every $F\in \cU^{\rho}$, there exists a unique $\varphi \in
(\cS)_{-\rho}$ so that $S\varphi =F$. This result  is proved in
\cite{HKPS} when $\rho=0$, and in \cite{Kond} when $\rho=1$.

\begin{definition}
\label{def:WickProd} For $\varphi$ and $\psi$ from
$(\cS)_{-\rho}$, $\rho\in [0,1]$, the Wick product
$\varphi\diamond \psi$ is the unique element of $(\cS)_{-\rho}$
whose S-transform is $S\varphi \cdot S\psi$.
\end{definition}

 If $S^{-1}$ is the inverse S-transform, then
$$
\varphi\diamond \psi=S^{-1}(S\varphi \cdot S\psi),
$$
Note that, by Theorem \ref{th:StrWN}, the Wick product is well
defined, because the space $\cU^{\rho}$, $\rho\in [0,1]$ is closed
under the point-wise multiplication. Theorem \ref{th:StrWN} also
ensures the correctness of the following definition of the white
noise.
\begin{definition}
\label{def:WNWN} The white noise $\dot{W}$ on $\bR^{\ell}$ is the
unique element of $(\cS)_{0}$ whose $S$ transform satisfies
$S\dot{W}(h)=h$.
\end{definition}

\begin{remark}
If $g \in L_p(\mbS)$, $p>1$, then $g \in (\cS)_{-0}$
\cite[Corollary 2.3.8]{HOUZ}, and
 the Fourier transform
$$
\hat{g}(h)=\int_{\cS'}\exp\left(\sqrt{-1}\langle \omega,
h\rangle\right) g(\omega)d\mu(\omega)
$$
is defined. Direct calculations \cite[Section 2.9]{HOUZ} show
that, for those  $g$,
$$
Sg(\sqrt{-1}\,h)=\hat{g}(h)\,e^{\frac{1}{2}\|h\|_{L_2(\bR^{\ell})}^2}.
$$
As a result,  the Wick product can be interpreted  as a
convolution on the infinite-dimensional space $(\cS)_{-\rho}$.
\end{remark}

In the  study of stochastic parabolic equations,  $\ell=d+1$ so
that the generic point from $\bR^{d+1}$ is written as $(t,x), \ t
\in \bR, \ x \in \bR^d$. As was mentioned earlier, the terms of
the type $fdW(t)$
  become  $f\diamond \dot{W}dt$.
The precise connection between the It\^{o} integral and Wick
product
 is discussed, for example, in \cite[Section 2.5]{HOUZ}.

As an example, consider the following equation:
\begin{equation}
\label{eq:spdeWN}
u_t(t,x)=a(x)u_{xx}(t,x)+b(x)u_x(t,x)+u_x(t,x)\diamond
\dot{W}(t,x),\ 0<t< T,\ x \in \bR,
\end{equation}
with initial condition $u(0,x)=u_0(x)$. In (\ref{eq:spdeWN}),
\begin{enumerate}
\item[(WN1)] $\dot{W}$ is the   white noise process on $\bR^{2}$.
\item[(WN2)] The initial condition $u_0$ and the coefficients
 $a$,  $b$ are bounded and
have continuous bounded derivatives up to second order.
\item[(WN3)] There exists a positive number $\varepsilon$ so that
 $a(x)\geq \varepsilon $, $x\in \bR$.
\item[(WN4)] The second-order  derivative of  $a$ is uniformly
H\"{o}lder continuous.
\end{enumerate}

The equivalent It\^{o} formulation of (\ref{eq:spdeWN}) is
\begin{equation}
\label{eq:ClWn} du(t,x)=(a(x)u_{xx}(t,x)+b(x)u_{x}(t,x))dt +
e_k(x)u_x(t,x)dw_k(x),
\end{equation}
where $\{e_k,\; k\geq 1\}$ is the Hermite basis in $L_2(\bR)$.

With $\cM_kv=e_kv_x$, we see that condition (\ref{eq:MCl}) does
not hold in any Sobolev space $H^{\gamma}_2(\bR)$. In fact,
 no traditional
solution exists in any normal triple of Sobolev space.
 On the other hand, with a suitable
definition of solution, equation (\ref{eq:spdeWN}) is solvable in
the space $(\cS)_{-0}$ of Hida distributions.

\begin{definition}
\label{def:defWN} A mapping $u:\bR^d\to (S)_{-\rho}$ is called
weakly differentiable with respect to  $x_i$ at a point $x^*\in
\bR^{\ell}$ if and only if there exists a $U_i(x^*)\in
(S)_{-\rho}$ so that, for all $\varphi \in (S)_{\rho}$,
$D_i\langle u(x), \varphi\rangle|_{x=x^*}=\langle U_i(x^*),
\varphi \rangle$. In that case, we write $U_i(x^*)=D_iu(x^*)$.
\end{definition}

\begin{definition}
\label{def:WNSol} A mapping $u$ from $[0,T]\times \bR$ to
$(\cS)_{-0}$ is called  a white noise solution of
(\ref{eq:spdeWN}) if and only if
\begin{enumerate}
\item The weak derivatives $u_t$, $u_x$, and $u_{xx}$ exist, in
the sense of Definition \ref{def:defWN}, for all $(t,x)\in
(0,T)\times \bR$. \item Equality (\ref{eq:spdeWN}) holds for all
$(t,x)\in (0,T)\times \bR^d$. \item $\lim_{t\downarrow 0}
u(t,x)=u_0(x)$ in the topology of $(\cS)_{-0}$.
\end{enumerate}
\end{definition}

\begin{theorem}
\label{th:spdeWN} Under assumptions (WN1)--(WN4), there exists a
white noise solution of (\ref{eq:spdeWN}). This solution is unique
in the class of weakly measurable mappings $v$ from $(0,T)\times
\bR$ to $(\cS)_{-0}$, for which there exists a non-negative
integer $q$ and a positive number $K$ so that
$$
\int_0^T \int_{\bR} \|v(t,x)\|_{-0,-q}e^{-Kx^2}dxdt < \infty.
$$
\end{theorem}
\begin{proof} Consider the  S-transformed equation
\begin{equation}
\label{eq:StrPdeWN}
F_{t}(t,x;h)=a(x)F_{xx}(t,x;h)+b(x)F_x(t,x;h)+F_x(t,x;h)h,
\end{equation}
$0<t< T,\ x \in \bR,\ h \in \cS(\bR),$ with initial condition
$F(0,x;h)=u_0(x)$. This a deterministic parabolic equation, and
one can  show, using the probabilistic representation of $F$, that
$F, F_t, F_x, $ and $F_{xx}$ belong to $\cU^{0}$. Then the inverse
S-transform of $F$ is a solution of (\ref{eq:spdeWN}), and the
uniqueness follows from the uniqueness for equation
(\ref{eq:StrPdeWN}). The details of the proof are in \cite{Pot},
where a similar equation is considered for $x \in\bR^d$.
\end{proof}

Even though the initial condition in (\ref{eq:spdeWN}) is
deterministic,
 there are no measurability restrictions on $u_0$ for the white noise solution
  to exist; see \cite{HOUZ} for more details.

With appropriate modifications, the white noise solution can be
defined for equations more general than (\ref{eq:spdeWN}). The
solution $F=F(t,x;h)$ of the  corresponding S-transformed equation
determines the regularity of the white noise solution
 \cite[Section 4.1]{HOUZ}.

Two main advantages of the white noise approach over the Hilbert
space approach are
\begin{enumerate}
\item no need for parabolicity condition; \item no measurability
restrictions on the input data.
\end{enumerate}
Still, there are substantial limitations:
\begin{enumerate}
\item There seems to be little or no connection between the white
noise solution and the traditional solution. While white noise
solution can, in principle, be constructed for equation
(\ref{eq:ParabEqCl}), this solution will be very different from
the traditional solution. \item There are no clear ways of
computing the solution numerically, even with available
representations of the Feynmann-Kac type \cite[Chapter 4]{HOUZ}.
\item The white noise  solution, being constructed on a special
white noise probability space,  is weak in the probabilistic
sense. Path-wise uniqueness does not apply to such solutions
because of the "averaging" nature of the solution spaces.
\end{enumerate}

\section{Generalized  Functions on the  Wiener Chaos Space}
\label{secGF}
\setcounter{equation}{0}
\setcounter{theorem}{0}

The objective of this section is to introduce the space of
generalized random elements on an arbitrary stochastic basis.

Let $\mbF=(\Omega, \cF, \{\cF_t\}_{t\geq 0}, \mbP)$ be a
stochastic basis with the usual assumptions and $Y$,  a separable
Hilbert space with inner product $( \cdot, \cdot )_Y$ and an
orthonormal basis $\{ y_k, \ k \geq 1\}$.
  On $\mbF$ and $Y$, consider a cylindrical Brownian motion $W$, that is, a family
  of continuous $\cF_t$-adapted Gaussian martingales  $W_y(t)$, $y\in Y$,
   so that $W_y(0)=0$ and
  $\mbE(W_{y_1}(t)W_{y_2}(s))=\min(t,s)(y_1,y_2)_Y$. In particular,
\begin{equation}
\label{eq:w_k} w_k(t)=W_{y_k}(t), \ k\geq 1, \ t\geq 0,
\end{equation}
 are independent standard Wiener processes on $\mbF$.

 Equivalently, instead of the process $W$, the starting point can be a system of
 independent standard Wiener processes $\{w_k,\; k\geq 1\}$ on $\mbF$.
 Then, given a separable Hilbert space $Y$ with an orthonormal basis
  $\{ y_k, \ k \geq 1\}$, the corresponding cylindrical Brownian motion
  $W$ is defined by
  \begin{equation}
  \label{eq:WCyl}
  W_y(t)=\sum_{k\geq 1} (y,y_k)_Y\,w_k(t).
  \end{equation}

 Fix a non-random  $T\in (0,\infty)$ and denote by $\cF^W_T$ the sigma-algebra
 generated by $w_k(t),\ k\geq 1, \ 0<t<T$.
Denote by $L_2(\mbW)$ the collection of $\cF^W_T$-measurable
square integrable random variables.

We now review construction of the Cameron-Martin basis in the
Hilbert space $L_2(\mbW)$.

Let  ${\mathfrak{m}}=\{m_k,\ k\geq 1\}$ be an orthonormal basis in
$L_2((0,T))$ so that each $m_k$ belongs to $L_{\infty}((0,T))$.
 Define the independent standard Gaussian random variables
$$
\xi_{ik}=\int_0^T m_i(s)dw_k(s).
$$
Consider the collection of multi-indices
 $$
 \cJ=\Big\{ \alpha =(\alpha_i^k,\ i,k\geq 1),\ \alpha_{i}^k\in
\{0,1,2,\ldots\},\ \sum_{i,k} \alpha_i^k<\infty \Big\}.
$$
The set $\cJ$ is countable, and, for every $\alpha \in \cJ$, only
finitely many of $\alpha_i^k$ are not equal to zero. The upper and
lower  indices in  $\alpha^k_i$ represent, respectively, the
 space and time components of the noise process $W$.
For $\alpha \in {\mathcal J}$, define
$$
|\alpha|=\sum_{i,k} \alpha_i^k, \ \alpha!=\prod_{i,k}\alpha_i^k!,
$$
and
\begin{equation}
\label{eq:CMB}
\xi_{\alpha}=\frac{1}{\sqrt{\alpha!}}\prod_{i,k}H_{\alpha_{i}^{k}}(\xi_{ik}),
\end{equation}
where $H_n$ is $n^{\rm th}$ Hermite polynomial. For example, if
$$
\alpha=\left(
\begin{array}{cccccccc}
 0 & 1 & 0 & 3 & 0 & 0 &\cdots\\
2  &  0 & 0 & 0 & 4 & 0 &\cdots\\
0 &0&0&0&0&0 & \cdots\\
\vdots &\vdots &\vdots &\vdots  &\vdots & \vdots&\cdots
\end{array}
\right)
$$
with four non-zero entries $\alpha^1_2=1; \ \alpha^1_4=3;\
\alpha^2_1=2; \ \alpha^2_5=4$, then
$$
 \xi_{\alpha}=\xi_{2,1}\cdot\frac{H_3(\xi_{4,1})}{\sqrt{3!}}\cdot
\frac{H_2(\xi_{1,2})}{\sqrt{2!}}\cdot
\frac{H_4(\xi_{5,2})}{\sqrt{4!}}.
$$

There are two main differences between (\ref{eq:CMBWN}) and
(\ref{eq:CMB}):
\begin{enumerate}
\item The basis (\ref{eq:CMB}) is constructed on an arbitrary
probability space. \item In (\ref{eq:CMB}), there is a clear
separation of the time and space components of the noise, and
explicit presence of the time-dependent functions $m_i$
facilitates the analysis of evolution equations.
\end{enumerate}

\begin{definition} The  space $L_2(\mbW)$ is called the Wiener Chaos space.
The $N$-th Wiener Chaos is the linear subspace of $L_2(\mbW)$,
generated by $\xi_{\alpha},\ |\alpha|=N$.
\end{definition}

The following is  another version of the classical results of
Cameron and Martin \cite{CM}.
\begin{theorem}
\label{th:CM} The collection $\Xi=\{ \xi_{\alpha},\ \alpha \in
\cJ\}$ is  an orthonormal basis in  $L_2(\mbW)$.
\end{theorem}
We refer to $\Xi$ as the Cameron-Martin basis in $L_2(\mbW)$. By
Theorem \ref{th:CM}, every element $v$ of $L_2(\mbW)$ can be
written as
$$
v=\sum_{\alpha \in \cJ} v_{\alpha}\xi_{\alpha},
$$
where $v_{\alpha}=\mbE (v\xi_{\alpha})$.

We now define the space $\cD(L_2(\mbW))$ of test functions and the
space $\cD'(L_2(\mbW);X)$ of $X$-valued generalized random
elements.

\begin{definition}
${ }$

(1) The space $\cD(L_2(\mbW))$ is the collection of elements from
$L_2(\mbW)$ that can be written in the form
$$
v=\sum_{\alpha \in \cJ_v} v_{\alpha}\xi_{\alpha}
$$
for some $v_{\alpha}\in \bR$ and a
 finite  subset $\cJ_v$ of $\cJ$. \\
(2) A sequence $v_n$ converges to $v$ in $\cD(L_2(\mbW))$ if and
only if $\cJ_{v_n}\subseteq \cJ_v$ for all $n$ and
$\lim\limits_{n\to \infty}|v_{n,\alpha} -v_{\alpha}| = 0$ for all
$\alpha$.
\end{definition}

\begin{definition}
\label{def:gen-el} For a linear topological space  $X$ define the
space\\
 $\cD'(L_2(\mbW);X)$ of  $X$-valued generalized random
elements as the collection of continuous linear maps from the
linear topological space $\cD(L_2(\mbW))$ to $X$. Similarly, the
elements of $\cD'(L_2(\mbW); L_1((0,T);X))$ are called $X$-valued
generalized random processes.
\end{definition}
 The element $u$ of $ \cD'(L_2(\mbW);X)$ can be identified with a formal Fourier
 series
$$
u=\sum_{\alpha \in \cJ} u_{\alpha} \xi_{\alpha},
$$
where $u_{\alpha} \in X$ are the  {\em generalized Fourier
coefficients} of $u$.
 For such a  series and for $v\in \cD(L_2(\mbW))$, we have
$$
u(v)=\sum_{\alpha \in \cJ_v} v_{\alpha}u_{\alpha}.
$$
Conversely, for $u\in\cD'(L_2(\mbW);X)$, we  define the formal
Fourier series of $u$ by setting  $u_{\alpha}=u(\xi_{\alpha})$.
 If $u \in L_2(\mbW)$, then $u \in \cD'(L_2(\mbW); \bR)$ and $u(v)=\mbE(uv)$.

 By Definition \ref{def:gen-el}, a  sequence $\{u_n, \; n\geq 1\}$
 converges to $u$ in
 $\cD'(L_2(\mbW);X)$ if and only if $u_{n}(v)$
 converges to $u(v)$ in the topology  of $X$ for every $v\in \cD(\mbW)$. In terms
 of generalized Fourier coefficients, this is equivalent to
  $\lim\limits_{n\to \infty}u_{n,\alpha}=u_{\alpha}$ in the topology of $X$ for every
 $\alpha\in \cJ$.

The construction of the space $\cD'(L_2(\mbW);X)$ can be extended
to Hilbert spaces other than $L_2(\mbW)$. Let $H$ be a real
separable Hilbert space with an orthonormal basis $\{e_k, \ k\geq
1\}$.  Define the space
$$
\cD(H)=\Big\{v\in H: v=\sum_{k\in \cJ_v} v_k e_k,\ v_k\in \bR,\
\cJ_v \ - \ {\rm a\  finite \ subset\ of \ } \{1,2,\ldots\}
\Big\}.
$$
By definition, $v_n$ converges to $v$ in $\cD(H)$ as $n\to \infty$
if and only if $\cJ_{v_n} \subseteq \cJ_v$ for all $n$ and
$\lim\limits_{n\to \infty}|v_{n,k}-v_k|=0$ for all $k$.

For  a linear topological  space $X$, $\cD'(H;X)$ is the space of
continuous linear maps from $\cD(H)$ to $X$. An element $g$ of
$\cD'(H;X)$ can be identified with a formal series $\sum_{k\geq 1}
g_k\otimes e_k$ so that $g_k=g(e_k)\in X$ and,  for $v\in \cD(H)$,
$g(v)=\sum_{k\in \cJ_v} g_kv_k$. If $X=\bR$ and $\sum_{k\geq 1}
g_k^2 < \infty$, then $g=\sum_{k\geq 1}g_ke_k \in H$ and $g(v)=(
g, v)_H$, the inner product in $H$. The space $X$ is naturally
imbedded into $\cD'(H; X)$: if  $u\in X$, then $\sum_{k\geq 1}
u\otimes e_k \in \cD'(H; X)$.

A sequence $g_n=\sum_{k\geq 1}g_{n,k}\otimes e_k,\ n\geq 1,$
converges to $g=\sum_{k\geq 1}g_{k}\otimes e_k$ in $\cD'(H;X)$ if
and only if, for every $k\geq 1,$ $\lim\limits_{n\to \infty}
g_{n,k}=g_k$ in the topology of $X$.

A collection   $\{\cL_k, \ k\geq 1\}$ of linear operators  from
$X_1$ to  $X_2$ naturally defines  a linear operator $\cL$ from
$\cD'(H; X_1)$ to $\cD'(H; X_2)$:
$$
\cL\left( \sum_{k\geq 1}g_k\otimes e_k\right) = \sum_{k\geq
1}\cL_k(g_k)\otimes e_k.
$$
Similarly, a linear operator $\cL: \cD'(H; X_1)\to\cD'(H; X_2)$
can be identified with a collection $\{\cL_k, \ k\geq 1\}$ of
linear operators  from  $X_1$ to  $X_2$ by setting
$\cL_k(u)=\cL(u\otimes e_k)$. Introduction of spaces $\cD'(H;X)$
and the corresponding operators makes it possible to avoid
conditions of the type (\ref{eq:MCl}).

\section{The Malliavin Derivative and its Adjoint}
\label{secMD}
\setcounter{equation}{0}
\setcounter{theorem}{0}

In this section, we define an analog of the It\^{o} stochastic
integral for generalized random processes.

All notations from the previous section will remain in force. In
particular,
 $Y$ is  a separable Hilbert space with a fixed
orthonormal basis $\{y_k, \ k\geq 1\}$, and  $\Xi=\{\xi_{\alpha},
\ \alpha \in \cJ\}$, the  Cameron-Martin  basis in $L_2(\mbW)$
defined in (\ref{eq:CMB}).

We start with a brief  review of the Malliavin calculus
\cite{Nualart}.

The {\em Malliavin derivative}  $\mbD$ is a continuous linear
operator from
\begin{equation}
\label{eq:MC-space} L_2^1(\mbW)= \Big\{ u \in
L_2(\mbW):\sum_{\alpha \in \cJ} |\alpha|u_{\alpha}^2< \infty\Big\}
\end{equation}
to $L_2\left(\mbW; (L_2((0,T))\times Y)\right)$. In particular,
\begin{equation}
\label{eq:MD} (\mbD\xi_{\alpha})(t)=\sum_{i,k}\sqrt{\alpha_i^k}
\xi_{\alpha^-(i,k)} m_i(t)y_k,
\end{equation}
where  $\alpha^-(i,k)$ is the multi-index with the  components
$$
\Big(\alpha^-(i,k)\Big)_j^l=\left\{
\begin{array}{ll}
\max(\alpha_i^k-1,0),& {\rm \ if \ } i=j \ {\rm and \ } k=l,\\
\alpha_j^l, & {\rm \ otherwise}.
\end{array}
\right.
$$
Note that, for each $t \in [0,T]$, $\mbD \xi_{\alpha}(t)\in
\cD(L_2(\mbW)\times Y).$ Using (\ref{eq:MD}), we extend the
operator $\mbD$  by linearity to  the space $\cD'(L_2(\mbW))$:
$$
\mbD\left(\sum_{\alpha\in \cJ} u_{\alpha} \xi_{\a}\right)=
\sum_{\alpha \in \cJ}\left(u_{\alpha} \sum_{i,k}\sqrt{\alpha_i^k}
\xi_{\alpha^-(i,k)} m_i(t)y_k\right).
$$

 For the sake of completeness and to justify  further definitions, let us establish
 connection between the Malliavin derivative and the stochastic It\^{o} integral.

  If  $u$ is an $\cF_t^W$-adapted process from $L_2\left(\mbW;  L_2((0,T);Y)\right)$,
   then $u(t)=\sum_{k\geq 1}u_k(t)y_k$, where  the random
 variable $u_k(t)$ is  $\cF_t^W$-measurable for each $t$ and $k$, and
 $$
 \sum_{k\geq 1} \int_0^T \mbE|u_k(t)|^2dt < \infty.
 $$
 We  define the stochastic It\^{o}  integral
 \begin{equation}
 \label{eq:int0}
U(t)= \int_0^t ( u(s), dW(s))_Y = \sum_{k\geq 1} \int_0^t
u_k(s)dw_k(s).
\end{equation}
Note that    $U(t)$ is $\cF_t^W$-measurable and $\mbE |U(t)|^2 =
\sum_{k\geq 1}\int_0^t\mbE|u_k(s)|^2ds$.

The next  result establishes a connection between the Malliavin
derivative
 and the stochastic It\^{o} integral.
\begin{lemma}
\label{lm:adj} Suppose that  $u$ is an $\cF^W_t$-adapted process
from \\
$L_2\left(\mbW; L_2((0,T);Y)\right)$, and define the process
$U$ according to (\ref{eq:int0}). Then, for every $0<t\leq T$ and
$\alpha \in \cJ$,
 \begin{equation}
 \label{eq:adj}
 \mbE (U(t)\xi_{\alpha})=\mbE\int_0^t ( u(s), (\mbD \xi_{\alpha})(s)
 )_Yds.
 \end{equation}
\end{lemma}

\begin{proof} Define $\xi_{\alpha}(t) =
\mbE(\xi_{\alpha}|\cF_t^W)$. It is known (see \cite{MR} or Remark
\ref{rm:xi(t)} below) that
\begin{equation}
\label{eq:xi(t)} d\xi_{\alpha}(t) = \sum_{i,k} \sqrt{\alpha_i^k}
 \xi_{\alpha^-(i,k)}(t)m_i(t) dw_k(t).
\end{equation}
Due to  $\cF_t^W$-measurability of $u_k(t)$, we have
 \begin{equation}
 \label{eq:u_pred}
 u_{k,\alpha}(t)=\mbE\Big( u_k(t) \mbE(\xi_{\alpha}|\cF^W_t)\Big)
 =\mbE(u_k(t)\xi_{\alpha}(t)).
 \end{equation}

 The definition of $U$ implies $dU(t)=\sum_{k\geq 1} u_k(t)dw_k(t)$, so that,
 by (\ref{eq:xi(t)}),  (\ref{eq:u_pred}), and the It\^{o} formula,
 \begin{equation}
 \label{eq:int1}
 U_{\alpha}(t) = \mbE (U(t)\xi_{\alpha})=
 \int_0^t\sum_{i,k} \sqrt{\alpha_i^k} u_{k,\alpha^-(i,k)}(s)m_i(s) ds.
\end{equation}
Together with (\ref{eq:MD}), the last equality implies
(\ref{eq:adj}). Lemma \ref{lm:adj} is proved.
\end{proof}

Note that the coefficients $u_{k,\alpha}$ of $u\in L_2(\mbW;
L_2((0,T);H))$ belong to $ L_2((0,T))$. We therefore define
$u_{k,\alpha,i}=\int_0^Tu_{k,\alpha}(t) m_i(t)dt$. Then, by
(\ref{eq:int1}),
\begin{equation}
\label{eq:int3} U_{\alpha}(T)=\sum_{i,k}\sqrt{\alpha_i^k}
u_{k,\alpha^-(i,k),i}.
\end{equation}
Since $U(T)=\sum_{\alpha \in \cJ} U_{\alpha}(T)\xi_{\alpha}$, we
shift the summation index in (\ref{eq:int3}) and    conclude that
\begin{equation}
\label{eq:plus} U(T)=\sum_{\alpha \in \cJ} \sum_{i,k}
\sqrt{\alpha_i^k+1} u_{k,\alpha,i} \xi_{\alpha^+(i,k)},
\end{equation}
where
\begin{equation}
\label{eq:xiplus} \Big(\alpha^+(i,k)\Big)_j^l=\left\{
\begin{array}{ll}
\alpha_i^k+1,& {\rm \ if \ } i=j \ {\rm and \ } k=l,\\
\alpha_j^l, & {\rm \ otherwise}.
\end{array}
\right.
\end{equation}
As a result, $U(T)=\delta(u)$, where $\delta$ is the adjoint of
the Malliavin derivative, also known as the Skorokhod integral;
 see \cite{Nualart} or \cite{NR} for details.

 Lemma \ref{lm:adj} suggests the following definition.
 For an $\cF^W_t$-adapted process $u$ from
 $  L_2\left(\mbW; L_2((0,T))\right)$, let  $\mbD_k^*u$ be
the $\cF^W_t$-adapted  process from  $L_2\left(\mbW;
L_2((0,T))\right)$ so that
\begin{equation}
\label{eq:opD1} (\mbD^*_k u)_{\alpha}(t) = \int_0^t\sum_{i}
\sqrt{\alpha_i^k} u_{\alpha^-(i,k)}(s)m_i(s) ds.
\end{equation}
If $u \in L_2\left(\mbW;  L_2((0,T);Y)\right)$ is
 $\cF^W_t$-adapted, then $u$ is in the domain of the operator $\delta$
 and $\delta(uI(s<t))=\sum_{k\geq  1}(\mbD_k^* u_k)(t)$.

We now extend the  operators $\mbD^*_k$
  to the  generalized random processes. Let $X$ be a Banach
space with norm $\|\cdot\|_X$.

\begin{definition}
If  $u$ is an $X$-valued generalized random process, then
$\mbD^*_ku$ is the $X$-valued generalized random process  so that
\begin{equation}
\label{eq:Dk}
(\mbD^*_ku)_{\alpha}(t)=\sum_i\int_0^tu_{\alpha^-(i,k)}(s)\sqrt{\alpha_i^k}
m_i(s)ds.
\end{equation}
If $ g \in \cD'\Big(Y;\cD'\left(L_2(\mbW);L_1((0,T);
X)\right)\Big) $, then $\mbD^*g$ is the $X$-valued generalized
random process so that, for $g=\sum_{ k\geq 1} g_{k}\otimes y_k,\;
g_k \in \cD'(L_2(\mbW);L_1((0,T); X))$,
\begin{equation}
\label{eq:D} (\mbD^*g)_{\alpha}(t)=
\sum_{k}(\mbD_k^*g_k)_{\alpha}(t)=
\sum_{i,k}\int_0^tg_{k,\alpha^-(i,k)}(s)\sqrt{\alpha_i^k}
m_i(s)ds.
\end{equation}
\end{definition}

Using (\ref{eq:MD}), we get a generalization of equality
(\ref{eq:adj}):
\begin{equation}
\label{eq:Da} (\mbD^*g)_{\alpha}(t) = \int_0^t g(\mbD
\xi_{\alpha}(s))(s)ds.
\end{equation}
Indeed, by linearity,
$$
g_k\left(\sqrt{\alpha_i^k}m_i(s)\xi_{\alpha^-(i,k)}\right)(s)=
\sqrt{\alpha_i^k}m_i(s)g_{k,\alpha^-(i,k)})(s).
$$

\begin{theorem}
\label{th:cont_op} If $T<\infty$, then $\mbD_k^*$ and $\mbD^*$ are
continuous linear operators.
\end{theorem}

\begin{proof} It is enough to show that, if
 $u,u_n \in \cD'\left(L_2(\cF_T^W); L_1((0,T); X)\right)$ and \\
$\lim_{n\to \infty}\|u_{\alpha}-u_{n,\alpha}\|_{L_1((0,T); X)}=0$
for every $\alpha\in \cJ$,
then, for every $k\geq 1$ and $\alpha \in \cJ$,\\
 $\lim_{n\to \infty}
  \|(\mbD^*_ku)_{\alpha}-(\mbD^*_ku_n)_{\alpha}\|_{L_1((0,T); X)}=0$.

 Using (\ref{eq:Dk}), we find
 $$
\|(\mbD^*_ku)_{\alpha}-(\mbD^*_ku_n)_{\alpha}\|_{X}(t) \leq
\sum_i\int_0^T\sqrt{\alpha_i^k}
\|u_{\alpha^-(i,k)}-u_{n,\alpha^-(i,k)}\|_X(s) |m_i(s)|ds.
$$
Note that the sum contains finitely many terms. By assumption,
$|m_i(t)|\leq C_i$, and so
\begin{equation*}
\label{eq:cont_op_pr}
\|(\mbD^*_ku)_{\alpha}-(\mbD^*_ku_n)_{\alpha}\|_{L_1((0,T);
X)}\!\leq\! C(\alpha)
\!\sum_i\!\!\sqrt{\alpha_i^k}\|u_{\alpha^-(i,k)}-u_{n,\alpha^-(i,k)}\|_{L_1((0,T);
X)}.
\end{equation*}
Theorem \ref{th:cont_op} is proved.
\end{proof}

\section{The Wiener Chaos Solution and the Propagator}
\label{secGS}
\setcounter{equation}{0}
\setcounter{theorem}{0}

In this section we build on the ideas from \cite{LMR} to
 introduce the Wiener Chaos solution and the corresponding
propagator for a general stochastic evolution equation. The
notations from Sections \ref{secGF} and \ref{secMD} will remain in
force. It will be convenient to interpret the cylindrical Brownian
motion  $W$ as a collection $\{w_k,\ k\geq 1\}$ of independent
standard Wiener processes. As before, $T\in (0,\infty)$ is fixed
and non-random.
 Introduce the following objects:
\begin{itemize}
\item The Banach spaces $A$, $X$, and  $U$   so that $U \subseteq
X$. \item Linear operators
\begin{eqnarray*}
\cA: L_1((0,T); A) \to L_1((0,T); X)\ {\rm and}\\
\cM_k: L_1((0,T); A) \to L_1((0,T); X).
\end{eqnarray*}
\item Generalized random processes $f \in
\cD'\left(L_2(\mbW);L_1((0,T); X)\right)$
 and\\
$g_k\in  \cD'\left(L_2(\mbW); L_1((0,T); X)\right).$ \item The
initial condition $u_0 \in \cD'\left(L_2(\mbW);U\right)$.
\end{itemize}

Consider the deterministic equation
\begin{equation}
\label{eq:det} v(t)=v_0+\int_0^t(\cA v)(s)ds + \int_0^t
\varphi(s)ds,
\end{equation}
where $v_0\in U$ and $\varphi \in L_1((0,T); X)$.

\begin{definition}
\label{def:det_main} A function $v$  is called a $w(A,X)$
solution of (\ref{eq:det}) if and only if $v\in L_1((0,T); A)$ and
equality (\ref{eq:det}) holds in the space $L_1((0,T); A)$.
\end{definition}

\begin{definition}
\label{def:main} An $A$-valued  generalized random process $u$ is
called a {\bf  $w(A,X)$ Wiener Chaos solution}  of the stochastic
differential equation
\begin{equation}
\label{eq:gen_eq} du(t)=(\cA u(t) + f(t))dt+(\cM_k u(t) +
g_k(t))dw_k(t),\ 0<t\leq T,\ u|_{t=0}=u_0,
\end{equation}
if and only if the equality
\begin{equation}
\label{eq:def1} u(t)=u_{0} + \int_0^t (\cA u +f)(s)ds+ \sum_{k\geq
1}(\mbD_k^*(\cM_k u + g_k))(t)
\end{equation}
 holds in $\cD'\left(L_2(\mbW); L_1((0,T); X)\right)$.
\end{definition}
Sometimes, to stress the dependence of the Wiener Chaos solution
on the terminal time $T$, the notation $w_T(A,X)$ will be used.

Equalities  (\ref{eq:def1}) (\ref{eq:D}) mean that, for every
$\alpha\in \cJ$, the generalized Fourier coefficient  $u_{\alpha}$
of $u$  satisfies
\begin{equation}
\label{eq:def2} u_{\alpha}(t)=u_{0,\alpha} + \int_0^t (\cA u
+f)_{\alpha}(s)ds+ \int_0^t\sum_{i,k}\sqrt{\alpha_i^k}(\cM_k u +
g_{k})_{\alpha^-(i,k)}(s)m_i(s)ds.
\end{equation}

  \begin{definition} System (\ref{eq:def2}) is called   the
propagator for equation (\ref{eq:gen_eq}).
\end{definition}

The propagator is a lower triangular system. Indeed, If
$\alpha=(0)$, that is, $|\alpha|=0$, then the corresponding
equation in (\ref{eq:def2}) becomes
\begin{equation}\label{eq:S(0)}
u_{(0)}(t)=u_{0,(0)} + \int_0^t
(\cA u_{(0)}(s) +f_{(0)}(s))ds.
\end{equation} If
$\alpha=(j\ell)$, that is, $\alpha^{\ell}_j=1$  for some fixed $j$
and $\ell$ and $\alpha^k_i=0$ for all other $i,k\geq 1$, then the
corresponding equation in (\ref{eq:def2}) becomes
\begin{equation}\label{eq:S(1)}
\begin{split}
u_{(j\ell)}(t)&=u_{0,(j\ell)} +
\int_0^t (\cA u_{(j\ell)}(s) +f_{(j\ell)}(s))ds
\\
&+ \int_0^t(\cM_k
u_{(0)}(s) + g_{\ell,(0)}(s))m_j(s)ds.
\end{split}
\end{equation}
 Continuing in this way,
we conclude that (\ref{eq:def2})
  can be solved by induction on $|\alpha|$
as long as the corresponding deterministic equation (\ref{eq:det})
is solvable. The precise result is as follows.

\begin{theorem}
\label{th:soft}
 If, for every $v_0\in U$ and $\varphi\in L_1((0,T);X)$, equation (\ref{eq:det})
has a unique $w(A,X)$ solution $v(t)=V(t,v_0, \varphi)$, then
equation (\ref{eq:gen_eq}) has a unique $w(A,X)$ Wiener Chaos
solution so that
\begin{equation}
\label{eq:S(a)}
\begin{split}
u_{\alpha}(t)&=V(t,u_{0,\alpha},f_{\alpha})+
\sum_{i,k}\sqrt{\alpha^k_i}V(t,0,m_i\cM_ku_{\alpha^-(i,k)})\\
&+
\sum_{i,k}\sqrt{\alpha^k_i}V(t,0,m_ig_{k,\alpha^-(i,k)}).
\end{split}
\end{equation}
\end{theorem}

\begin{proof} Using the assumptions of the theorem and linearity,
we conclude that (\ref{eq:S(a)}) is the unique solution of
(\ref{eq:def2}).
\end{proof}


To derive a more  explicit formula for $u_{\alpha}$, we need some
additional constructions. For every multi-index $\a$ with
$|\a|=n$, define  the {\bf characteristic set} $K_{\alpha}$ of
$\alpha$ so that
$$
K_{\a}= \{ (i_1^{\a},k_1^{\a}), \ldots, (i_n^{\a},k_n^{\a}) \},
$$
$i_1^{\a} \leq i_2^{\a}\leq \ldots \leq i_n^{\a}$, and if
$i_j^{\a}=i_{j+1}^{\a}$, then $k_j^{\a} \leq k_{j+1}^{\a}$. The
first pair $(i_1^{\a},k_1^{\a})$ in $K_{\a}$ is the position
numbers of the first nonzero element of $\a$. The second pair is
the same as the first if the first nonzero element of $\a$ is
greater than one; otherwise, the second pair is the position
numbers of the second nonzero element
 of $\a$ and so on. As a result, if $\a^k_i>0$, then  exactly $\a_i^k$
 pairs in $K_{\a}$ are equal to $(i,k)$.  For example, if
$$ \a=\left( \ba{llllllll}
0 &1 &0 &2 &3 &0 &0 &\cdots \\
1 &2 &0 &0 &0 &1 &0 &\cdots\\
0 & 0&0 &0 &0 &0 & 0& \cdots\\
\vdots&\vdots&\vdots&\vdots&\vdots&\vdots&\vdots&\cdots \ea
\right) $$ with  nonzero elements
  $$
  \a_1^2=\a_2^1=\a_1^6=1,\ \a_2^2=\a_4^1=2,\
 \a_5^1=3,
 $$
 then  the  characteristic set is
$$
K_{\a}\!=\!\{ (1,2),\,(2,1),\,(2,2),\,(2,2),\,(4,1),\,(4,1),\,
(5,1),\,(5,1),\,(5,1),\,(6,2) \}.
$$

 \begin{theorem}
\label{S.th} Assume that
\begin{enumerate}
\item for every $v_0\in U$ and $\varphi\in L_1((0,T);X)$, equation
(\ref{eq:det}) has a unique $w(A,X)$ solution $v(t)=V(t,v_0,
\varphi)$, \item the input data in (\ref{eq:def2}) satisfy
 $g_k=0$ and $f_{\alpha}=u_{0,\alpha}=0$ if $|\alpha|>0$.
\end{enumerate}
Let $u_{(0)}(t)=V(t,u_0,0)$ be the solution of (\ref{eq:def2}) for
$|\alpha|=0$. For $\alpha \in \cJ$ with $|\alpha|=n\geq 1$ and the
characteristic set $K_{\alpha}$, define functions
$F^n=F^n(t;\alpha)$ by induction as follows:
\begin{equation}
\begin{split}
F^1(t;\alpha)&=V(t,0,m_i\cM_{k}u_{(0)}) \ {\rm  if} \  K_{\a}=\{(i,k)\};\\
F^n(t;\alpha)&=\sum_{j=1}^nV(t,0,m_{i_j}\cM_{k_j}F^{n-1}(\cdot;\alpha^-(i_j,k_j)))
\\
 &{\rm if \ } K_{\alpha}=\{(i_1, k_1), \ldots, (i_n, k_n)\}.
\end{split}
\end{equation}

Then
\begin{equation}
\label{eq:ind}
u_{\alpha}(t)=\frac{1}{\sqrt{\alpha!}}F^n(t;\alpha).
\end{equation}
\end{theorem}

\begin{proof}
  If $|\alpha|=1$, then  representation (\ref{eq:ind}) follows from (\ref{eq:S(1)}).
 For $|\alpha|>1$, observe that
\begin{itemize}
\item If $\bar{u}_{\alpha}(t)=\sqrt{\alpha!}u_{\alpha}$ and
$|\alpha|\geq 1$, then (\ref{eq:def2}) implies
$$
\bar{u}(t)=\int_0^t\cA \bar{u}_{\alpha}(s)ds +
\sum_{i,k}\int_0^t\alpha^k_im_i(s)\cM_k\bar{u}_{\alpha^-(i,k)}(s)ds.
$$
\item If $K_{\alpha}=\{(i_1, k_1), \ldots, (i_n, k_n)\}$, then,
for every $j=1,\ldots, n$, the characteristic set
$K_{\alpha^-(i_j,k_j)}$ of $\alpha^-(i_j,k_j)$ is obtained from
$K_{\alpha}$ by removing the pair $(i_j,k_j)$. \item By  the
definition of the characteristic set,
$$
\sum_{i,k}\alpha^k_im_i(s)\cM_k\bar{u}_{\alpha^-(i,k)}(s)=
\sum_{j=1}^n m_{i_j}(s)\cM_{k_j}\bar{u}_{\alpha^-(i_j,k_j)}(s).
$$
\end{itemize}
As a result, representation  (\ref{eq:ind}) follows by
induction on $|\alpha|$ using (\ref{eq:S(a)}): \\
if  $|\alpha|=n>1$, then \begin{equation}\begin{split}
\bar{u}_{\alpha}(t)&=
\sum_{j=1}^nV(t,0,m_{i_j}\cM_{k_j}\bar{u}_{\alpha^-(i_j,k_j)})\\
&=\sum_{j=1}^n V(t,0,m_{i_j}\cM_{k_j}
 F^{(n-1)}(\cdot;\alpha^-(i_j,k_j))=
 F^n(t;\alpha).
 \end{split}
 \end{equation}

 Theorem \ref{S.th} is proved.
 \end{proof}

\begin{corollary}
\label{cor:ind} Assume that the operator $\cA$ is a generator of a
strongly continuous semi-group $\Phi=\Phi_{t,s}, \ t\geq s\geq
0,$ in some Hilbert space $H$ so that $A\subset H$, each $\cM_k$
is a bounded operator from $A$ to $H$, and
 the solution $V(t,0,\varphi)$
of equation (\ref{eq:det}) is written as
\begin{equation}
\label{eq:semigr} V(t,0,\varphi)=\int_0^T \Phi_{t,s} \f(s)ds, \ \
\varphi \in L_2((0,T); H)).
\end{equation}
Denote by $\PP^n$ the permutation group of $\{1, \ldots, n\}.$ If
$u_{(0)}\in L_2((0,T); H))$, then, for $|\alpha|=n>1$ with the
characteristic  set $K_{\alpha}=\{(i_1, k_1), \ldots, (i_n,
k_n)\},$ representation (\ref{eq:ind}) becomes
\begin{equation}
\begin{split}
\label{eq:ind_sg} u_{\a}(t)&=\frac{1}{\sqrt{\a!}}\sum_{\sigma \in
\PP^n}
\int_0^t\int_0^{s_{n}}\ldots \int_0^{s_2}\\
&\Phi_{t,s_n}\cM_{k_{\sigma(n)}}\cdots
\Phi_{s_2,s_1}\cM_{k_{\sigma(1)}}u_{(0)}(s_1)
m_{i_{\sigma(n)}}(s_n) \cdots m_{i_{\sigma(1)}}(s_1) ds_1\ldots
ds_n.
\end{split}
\end{equation}
Also,
\begin{equation}
\label{eq:iter_int}%
\begin{split}
\sum_{|\alpha|=n}u_{\alpha}(t)\xi_{\alpha}  &  = \sum_{k_{1},
\ldots, k_{n}\geq1}\int_{0}^{t}\int_{0}^{s_{n}}\ldots\int_{0}^{s_{2}}\\
& \!\!\!\!\!\!\!\!\!\!  \!\!\!\!\!\!\!
\Phi_{t,s_{n}}{{\mathcal{M}}}_{k_{n}}\cdots
\Phi_{s_{2},s_{1}}\left(
{{\mathcal{M}}}_{k_{1}}u_{(0)}+g_{k_{1}}(s_{1})\right)
dw_{k_{1}}(s_{1})\cdots dw_{k_{n}}(s_{n}), \ n\geq1,
\end{split}
\end{equation}
and, for every Hilbert space $X$,  the following energy equality
holds:
\begin{equation}
\label{eq:ind_en}
\begin{split}
\sum_{|\alpha|=n}\|u_{\alpha}(t)\|^2_{X} & =
\sum_{k_1, \ldots, k_n=1}^{\infty}\int_0^t\int_0^{s_{n}}\ldots \int_0^{s_2}\\
& \|\Phi_{t,s_n}\cM_{k_n}\cdots
\Phi_{s_2,s_1}\cM_{k_{1}}u_{(0)}(s_1)\|_X^2 ds_1\ldots ds_n;
\end{split}
\end{equation}
both sides in the last equality can be infinite. For $n=1$,
formulas (\ref{eq:ind_sg}) and (\ref{eq:ind_en}) become
\begin{equation}
\label{eq:ind_sg1}
u_{(ik)}(t)=\int_0^t\Phi_{t,s}\cM_{k}u_{(0)}(s)\ m_i(s)ds;
\end{equation}
\begin{equation}
\label{eq:ind_en1} \sum_{|\alpha|=1} \|u_{\alpha}(t)\|_X^2 =
\sum_{k=1}^{\infty} \int_0^t \|\Phi_{t,s}\cM_{k}u_{(0)}(s)\|^2_X
ds.
\end{equation}
\end{corollary}

\begin{proof} Using the semi-group representation
(\ref{eq:semigr}),
 we conclude that  (\ref{eq:ind_sg}) is just an
expanded version of  (\ref{eq:ind}).

Since $\{m_i, \; i\geq 1\}$ is an orthonormal basis in $L_2(0,T)$,
 equality (\ref{eq:ind_en1}) follows from (\ref{eq:ind_sg1})
and the Parcevall identity. Similarly, equality (\ref{eq:ind_en})
will follow from (\ref{eq:ind_sg}) after an application of an
appropriate Parcevall's identity.

To carry out the necessary arguments when $|\a|>1$, denote by
$\cJ_1$ the collection of one-dimensional multi-indices
$\beta=(\beta_1, \beta_2, \ldots)$ so that each $\beta_i$ is a
non-negative integer and $|\beta|=\sum_{i\geq 1} \beta_i <
\infty$. Given a $\beta \in \cJ_1$ with $|\beta|=n$, we define
$K_{\beta}=\{ i_1, \ldots, i_n\}$, the characteristic set of
$\beta$ and the function
\begin{equation}
E_{\beta}(s_1,\ldots, s_n)=\frac{1}{\sqrt\beta!n!}\sum_{\sigma \in
\PP^n} m_{i_1}(s_{\sigma(1)}) \cdots m_{i_n}(s_{\sigma(n)}).
\end{equation}
By construction, the collection $\{E_{\beta}, \beta\in \cJ_1,
|\beta|=n\}$ is an orthonormal basis in the sub-space of
symmetric functions in $L_2((0,T)^n; X)$.

Next, we re-write (\ref{eq:ind_sg}) in a symmetrized form. To make
the notations shorter, denote by $s^{(n)}$ the ordered set
$(s_1,\ldots, s_n)$ and write $ds^n=ds_1\ldots ds_n$. Fix $t\in
(0,T]$ and the set   $k^{(n)}=\{k_1, \ldots, k_n\}$ of the second
components of the characteristic set $K_{\alpha}$. Define the
symmetric function
\begin{equation}
\label{eq:GGG}
\begin{split}
&G(t,k^{(n)};s^{(n)})\\
&=\frac{1}{\sqrt{n!}}\sum_{\sigma\in \PP^n}
\Phi_{t,s_{\sg(n)}}\cM_{k_n}\cdots
\Phi_{s_{\sg(2)},s_{\sg(1)}}\cM_{k_1}
u_{(0)}(s_{\sg(1)})1_{s_{\sg(1)}<\cdots<s_{\sg(n)}<t}(s^{(n)}).
\end{split}
\end{equation}
Then (\ref{eq:ind_sg}) becomes
\begin{equation}
\label{eq:619}
u_{\alpha}(t)=\int_{[0,T]^n}G(t,k^{(n)};s^{(n)})E_{\beta(\alpha)}(s^{(n)})ds^n,
\end{equation}
where the multi-indices $\alpha$ and $\beta(\alpha)$ are related
via their characteristic sets: if
$$
K_{\alpha}=\{(i_1, k_1), \ldots, (i_n, k_n)\},
$$
then
$$
K_{\beta(\alpha)}=\{i_1, \ldots, i_n\}.
$$
Equality (\ref{eq:619}) means that,
  for fixed $k^{(n)}$, the function $u_{\alpha}$ is a Fourier
coefficient of the symmetric function $G(t,k^{(n)};s^{(n)})$ in
the space $L_2((0,T)^n; X)$.
 Parcevall's identity and  summation over all possible $k^{(n)}$
 yield
$$
\sum_{|\alpha|=n} \|u_{\a}(t)\|^2_X =\frac{1}{n!}\sum_{k_1,
\ldots, k_n=1}^{\infty} \int_{[0,T]^n}
\|G(t,k^{(n)};s^{(n)})\|^2_{X}ds^n,
$$
which, due  to (\ref{eq:GGG}),  is the same as (\ref{eq:ind_en}).

To prove equality (\ref{eq:iter_int}), relating the Cameron-Martin
and multiple It\^{o} integral expansions of the solution, we use
the following result \cite[Theorem 3.1]{Ito}:
$$
\xi_{\a}=\frac{1}{\sqrt{\a!}}\int_0^T\int_0^{s_n}\cdots
\int_0^{s_2} E_{\beta(\alpha)}(s^{(n)})dw_{k_1}(s_1)\cdots
dw_{k_n}(s_n);
$$
see also \cite[pp. 12--13]{Nualart}. Since the collection of all
$E_{\beta}$ is an orthonormal basis, equality
  (\ref{eq:iter_int}) follows from (\ref{eq:619}) after summation
over al $k_1,\ldots, k_n$.

Corollary \ref{cor:ind} is proved.
\end{proof}


We now present several examples   to illustrate the general
results.
\begin{example}
\label{ex1} {\rm Consider the following equation:
\begin{equation}
\label{eq:classical} du(t,x)=(au_{xx}(t,x)+f(t,x))dt+(\sigma
u_x(t,x)+g(t,x))dw(t),\ t>0,\ x \in \bR,
\end{equation}
where $a>0$, $\sigma \in \bR$, $f \in L_2((0,T); H^{-1}_2(\bR))$,
$g  \in L_2((0,T); L_2(\bR))$, and $u|_{t=0}=u_0 \in L_2(\bR)$.  By Theorem
\ref{th:ParabCl1},
 if $\sigma^2<2a$, then
equation (\ref{eq:classical}) has a unique traditional solution
$u\in L_2\left(\mbW; L_2((0,T); H^1_2(\bR))\right)$.

By $\cF^W_t$-measurability of $u(t)$, we have
$$
\mbE(u(t)\xi_{\a})=\mbE(u(t)\mbE(\xi_{\a}|\cF^W_t)).
$$
Using the relation (\ref{eq:xi(t)}) and the It\^{o} formula, we
find that $u_{\a}$ satisfy
$$
du_{\a}=a(u_{\a})_{xx}dt+\sum_i\sqrt{\a_i}\sigma
(u_{\a^{-}(i)})_xm_i(t)dt,
$$
which is precisely the propagator for equation
(\ref{eq:classical}). In other words, if $2a>\sigma^2$, then the
traditional
 solution of (\ref{eq:classical}) coincides with the
Wiener Chaos solution.

On the other hand,  the heat equation
$$
v(t,x)=v_0(x)+\int_0^tv_{xx}(s,x)ds+\int_0^t\varphi(s,x)ds,\ v_0
\in L_2(\bR)
$$
with $ \varphi \in L_2((0,T); H^{-1}_2(\bR))$ has a unique
$w(H^1_2(\bR),H^{-1}_2(\bR))$ solution. Therefore,  by Theorem
\ref{th:soft}, the unique $w(H^1_2(\bR), H^{-1}_2(\bR))$ Wiener
Chaos solution of (\ref{eq:classical}) exists for all $\sigma \in
\bR$.}
\end{example}

In the next example, the equation, although not parabolic,
 can be solved explicitly.
\begin{example}
\label{ex2} {\rm Consider the following equation:
\begin{equation}
\label{eq:nclassical} du(t,x)= u_x(t,x)dw(t),\ t>0,\ x \in \bR; \
\ u(0,x)=x.
\end{equation}
Clearly, $u(t,x)=x+w(t)$ satisfies  (\ref{eq:nclassical}).

To find the Wiener Chaos solution of (\ref{eq:nclassical}), note
that, with one-dimensional Wiener process, $\alpha_i^k=\alpha_i$,
and the propagator in this case becomes
$$
u_{\alpha}(t,x)=xI(|\alpha|=0)+\int_0^t \sum_{i}
\sqrt{\alpha_i}(u_{\alpha^-(i)}(s,x))_xm_i(s)ds.
$$
Then  $u_{\alpha}=0$ if $|\alpha|>1$, and
\begin{equation}
\label{eq:ex2} u(t,x)=x+\sum_{i\geq 1}
\xi_i\int_0^tm_i(s)ds=x+w(t).
\end{equation}

Even though  Theorem \ref{th:soft} does not apply, the above
arguments show that
 $u(t,x)=x+w(t)$ is
  the unique  $w(A,X)$ Wiener Chaos solution of (\ref{eq:nclassical})
for suitable spaces $A$ and $X$, for example,
$$
X=\left\{f: \int_{\bR} (1+x^2)^{-2}f^2(x)dx < \infty\right\}\ {\rm
and}\ A=\{f: f, f' \in X\}.
$$
Section \ref{sec:.:frd} provides a  more detailed analysis of
equation (\ref{eq:nclassical}). }
\end{example}

 If equation (\ref{eq:gen_eq}) is anticipating, that is, the initial
 condition is not deterministic  and/or the free terms $f,g$ are not
 $\cF^W_t$-adapted, then the Wiener Chaos solution generalizes the Skorohod
 integral interpretation of the equation.

 \begin{example}
 \label{ex:ant}
 {\rm Consider the equation
 \begin{equation}
 \label{eq:ant}
 du(t,x)=\frac{1}{2}u_{xx}(t,x)dt + u_x(t,x)dw(t), \ t\in (0,T],\ x\in \bR,
 \end{equation}
 with initial condition $u(0,x)=x^2w(T)$. Since $w(T)=\sqrt{T} \xi_1$, we
 find
 \begin{equation}
 \label{eq:Sant}
 (u_{\alpha})_t(t,x)=\frac{1}{2} (u_{\alpha})_{xx}(t,x)+ \sum_{i} \sqrt{\alpha_i}
 m_i(t)(u_{\alpha^-(i)})_x(t,x)
 \end{equation}
 with initial condition $u_{\alpha}(0,x)=\sqrt{T}x^2I(|\alpha|=1, \alpha_1=1)$.
 By Theorem \ref{th:soft}, there exists a
   unique  $w(A,X)$ Wiener Chaos solution of (\ref{eq:ant})
for suitable spaces $A$ and $X$. For example, we can take
$$
X=\left\{f: \int_{\bR} (1+x^2)^{-8}f^2(x)dx < \infty\right\}\ {\rm
and}\ A=\{f: f, f', f'' \in X\}.
$$
System (\ref{eq:Sant}) can be solved explicitly. Indeed,
$u_{\alpha}\equiv 0$ if $|\alpha|=0$ or $|\alpha|>3$ or if
$\alpha_1=0$. Otherwise, writing $M_i(t)=\int_0^tm_i(s)ds$, we
find:
\begin{equation*}
\begin{split}
u_{\alpha}(t,x)&=(t+x^2)\sqrt{T}, \ {\rm if } \ |\alpha|=1,\ \alpha_1=1;\\
u_{\alpha}(t,x)&=2\sqrt{2}\;xt, \ {\rm if } \ |\alpha|=2,\ \alpha_1=2;\\
u_{\alpha}(t,x)&=2\sqrt{T}\;xM_i(t), \ {\rm if } \ |\alpha|=2,\
\alpha_1=\alpha_i=1, \
1<i;\\
u_{\alpha}(t,x)&=\sqrt{\frac{6}{T}}\;t^2, \ {\rm if } \ |\alpha|=3,\ \alpha_1=3;\\
u_{\alpha}(t,x)&=2\sqrt{2T}\;M_1(t)M_i(t), \ {\rm if } \
|\alpha|=3,\ \alpha_1=2,\
\alpha_i=1, \ 1<i;\\
u_{\alpha}(t,x)&=\sqrt{2T}\;M^2_i(t), \ {\rm if } \ |\alpha|=3,\
\alpha_1=1,\
\alpha_i=2, \ 1<i;\\
u_{\alpha}(t,x)&=2\sqrt{T}\;M_i(t)M_j(t), \ {\rm if } \
|\alpha|=3,\ \alpha_1= \alpha_i=\alpha_j=1, \ 1<i<j.
\end{split}
\end{equation*}
Then
\begin{equation}
\label{eq:antsol} u(t,x)=\sum_{\alpha\in \cJ}
u_{\alpha}\xi_{\alpha}=w(T)w^2(t)-2tw(t)+2(W(T)w(t)-t)x+ x^2w(T)
\end{equation}
is the Wiener Chaos solution of (\ref{eq:ant}). It can be verified
using the properties of the Skorohod integral \cite{Nualart} that
the function $u$ defined by (\ref{eq:antsol}) satisfies
\begin{equation*}
u(t,x)=x^2w(T)+\frac{1}{2}\int_0^t u_{xx}(s,x)ds + \int_0^t
u_x(s,x)dw(s),\ t\in [0,T],\ x\in \bR,
\end{equation*}
where the stochastic integral is in the sense of Skorohod. }
\end{example}

\section{Weighted Wiener Chaos Spaces and S-Transform}
\label{secWS}
\setcounter{equation}{0}
\setcounter{theorem}{0}

The space $\cD'(L_2(\mbW); X)$ is too big to provide any
reasonable information about regularity of the Wiener Chaos
solution.
 Introduction of weighted Wiener
chaos spaces makes it possible to resolve this difficulty.

As before, let  $\Xi=\{\xi_{\a},\ \a \in \cJ\}$ be the
Cameron-Martin basis in $L_2(\mbW)$, and $\cD(L_2(\mbW);X)$, the
collection of finite linear combinations of $\xi_{\a}$ with
coefficients in a Banach space $X$.

\begin{definition}
\label{def:ws-gen} Given a collection $\{r_{\a},\ \a \in \cJ\}$ of
positive numbers, the space $\cR L_2(\mbW;X)$ is the closure of
$\cD(L_2(\mbW);X)$ with respect to the norm
$$
\|v\|_{\cR L_2(\mbW;X)}^2:=\sum_{\a \in \cJ} r_{\a}^2
\|v_{\a}\|^2_X.
$$
\end{definition}

The operator $\cR$ defined by $(\cR v)_{\a}:=r_{\a}v_{\a}$ is a
linear homeomorphism from $\cR L_2(\mbW;X)$ to $L_2(\mbW;X)$.

There are several  special choices of the weight sequence
$\cR=\{r_{\alpha}, \; \a \in \cJ\}$ and  special notations for the
corresponding weighted Wiener chaos spaces.
\begin{itemize}
\item If $Q=\{q_1, q_2, \ldots\}$ is  a sequence  of
  positive numbers,  define
 $$
  q^{\alpha} = \prod_{i,k} q_k^{\alpha_i^k}.
  $$
  The operator $\cR$, corresponding to $r_{\alpha}=q^{\alpha}$, is denotes by
  $\cQ$. The space $\cQ L_2(\mbW;X)$ is denoted by $L_{2,Q}(\mbW; X)$
  and is called a  {\em  Q-weighted Wiener chaos space.} The significance of this
  choice of weights will be explained shortly (see, in particular, Proposition
  \ref{prop:QT}).
  \item If
  $$
r_{\alpha}^2=(\alpha!)^{\rho}\prod_{i,k}(2ik)^{\gamma\alpha^k_i},\
\rho,\gamma \in \bR,
$$
then the corresponding space $\cR L_2(\mbW;X)$ is denoted by
$(\cS)_{\rho,\gamma}(X)$. As always, the argument $X$ will be
omitted if $X=\bR$. Note the analogy with Definition
\ref{def:KS1}.
\end{itemize}

The structure of weights in the spaces $L_{2,Q}$ and
$(\cS)_{\rho,\gamma}$ is different, and in general these two
classes of spaces are not related. There exist  generalized random
elements that belong to some $L_{2,Q}(\mbW;X)$, but do not belong
to any  $(\cS)_{\rho,\gamma}(X)$. For example, $u=\sum_{k\geq
1}e^{k^2}\xi_{1,k}$ belongs to $L_{2,Q}(\mbW)$ with
$q_k=e^{-2k^2}$, but to no $(\cS)_{\rho,\gamma}$, because the sum
$\sum_{k\geq 1}e^{2k^2}(k!)^{\rho}(2k)^{\gamma}$ diverges for
every $\rho,\gamma\in\bR$.  Similarly, there exist  generalized
random elements that belong to some  $(\cS)_{\rho,\gamma}(X)$,
but to no $L_{2,Q}(\mbW;X)$. For example, $u=\sum_{n\geq 1}
\sqrt{n!}\xi_{(n)}$, where $(n)$ is the multi-index with
$\alpha^1_1=n$ and $\alpha^k_i=0$ elsewhere, belongs to
$(S)_{-1,-1},$
 but does not belong to any
  $L_{2,Q}(\mbW)$, because the sum $\sum_{n\geq 1}q^nn!$ diverges for
 every $q>0$.

The next result  is the space-time analog of Proposition 2.3.3 in
\cite{HOUZ}.

\begin{proposition}
\label{prop:ConvSum} The sum
$$
\sum_{\alpha \in \cJ}\prod_{i,k\geq 1}(2ik)^{-\gamma\alpha^k_i}
$$
converges if and only if $\gamma>1$.
\end{proposition}
\begin{proof} Note that
\begin{equation}
\label{eq:CONV} \sum_{\alpha \in \cJ}\prod_{i,k\geq
1}(2ik)^{-\gamma\alpha^k_i}= \prod_{i,k\geq 1} \left( \sum_{n\geq
0} ((2ik)^{-\gamma})^n \right) =\prod_{i,k}
\frac{1}{(1-(2ik)^{-\gamma})},\ \gamma>0
\end{equation}
The infinite product on the right of (\ref{eq:CONV}) converges if
and only if each of the sums $\sum_{i\geq 1}i^{-\gamma}$,
$\sum_{k\geq 1}k^{-\gamma}$ converges, that is, if an only if
$\gamma>1$.
\end{proof}

\begin{corollary}
\label{cor:u-space}
For every $u \in \cD'(\mbW; X)$, there exists an operator $\cR$ so that\\
$\cR u \in L_2(\mbW; X)$.
\end{corollary}
\begin{proof} Define
$$
r_{\alpha}^2=\frac{1}{1+\|u_{\a}\|_X^2}\prod_{i,k\geq
1}(2ik)^{-2\alpha^k_i}.
$$
Then
\begin{equation*}
\|\cR u\|_{L_2(\mbW; X)}^2= \sum_{\a \in \cJ}
\frac{\|u_{\a}\|_X^2}{1+\|u_{\a}\|_X^2} \prod_{i,k\geq
1}(2ik)^{-2\alpha^k_i} \leq \sum_{\a \in \cJ}\prod_{i,k\geq
1}(2ik)^{-2\alpha^k_i} < \infty.
\end{equation*}
\end{proof}

The importance of the operator $\cQ$ in the study of stochastic
equations is due to the fact that the operator $\cR$ maps a Wiener
Chaos solution to a Wiener Chaos solution if and only $\cR=\cQ$
for some sequence $Q$. Indeed, direct calculations show that the
functions $u_{\alpha}, \alpha\in \cJ,$ satisfy the propagator
(\ref{eq:def2}) if and only if $v_{\alpha}=(\cR u)_{\alpha}$
satisfy
\begin{equation}
\label{eq:def2R}
\begin{split}
v_{\alpha}(t)&=(\cR u_{0})_{\alpha} + \int_0^t (\cA v +\cR f)_{\alpha}(s)ds\\
&+
\int_0^t\sum_{i,k}\sqrt{\alpha_i^k}\frac{\rho_{\alpha}}{\rho_{\a^-(i,k)}}
(\cM_k \cR u + \cR g_{k})_{\alpha^-(i,k)}(s)m_i(s)ds.
 \end{split}
\end{equation}
Therefore, the operator $\cR$ preserves the structure of the
propagator if and only if
$$
\frac{\rho_{\alpha}}{\rho_{\a^-(i,k)}}=q_k,
$$
that is,  $\rho_{\alpha}=q^{\alpha}$ for some sequence $Q$.

Below is the summary of the main properties of the operator $\cQ$.
\begin{proposition}
\label{prop:QT}
$           $\\
\begin{enumerate}
\item If $q_k\leq q<1$ for all $k\geq 1$, then
$L_{2,Q}(\mbW)\subset (\cS)_{0,-\gamma}$ for some $\gamma>0$.
\item If $q_k\geq q >1$ for all $k$, then $L_{2,Q}(\mbW)\subset
L_2^n(\mbW)$ for all $n\geq 1$, that is, the elements of
$L_{2,Q}(\mbW)$ are infinitely differentiable in the  Malliavin
sense. \item If $u \in L_{2,Q}(\mbW; X)$ with generalized  Fourier
coefficients
 $u_{\alpha}$ satisfying the propagator (\ref{eq:def2}),  and $v=\cQ u$, then
the corresponding system for the generalized Fourier coefficients
of $v$  is
\begin{equation}
\label{eq:def2Q}
\begin{split}
v_{\alpha}(t)&=(\cQ u_{0})_{\alpha} + \int_0^t (\cA v +\cQ f)_{\alpha}(s)ds\\
&+ \int_0^t\sum_{i,k}\sqrt{\alpha_i^k}(\cM_k v +
 \cQ g_{k})_{\alpha^-(i,k)}(s)q_km_i(s)ds.
 \end{split}
\end{equation}
\item The function  $u$ is a Wiener Chaos solution of
\begin{equation}
\label{eq:clas1} u(t)= u_0+ \int_0^t(\cA u(s) + f(s))dt + \int_0^t
( \cM u(s) + g(s), dW(s) )_Y
\end{equation}
if and only if $v=\cQ u$ is a Wiener Chaos solution of
\begin{equation}
\label{eq:clas} v(t)= (\cQ u)_0+ \int_0^t(\cA v(s) + \cQ f(s))dt +
\int_0^t ( \cM v(s) + \cQ g(s), dW^Q(s) )_Y,
\end{equation}
where, for $h\in Y$,  $W^Q_h(t)=\sum_{k\geq 1}(h,y_k)_Yq_kw_k(t)$.
\end{enumerate}
\end{proposition}

The following examples demonstrate how the operator $\cQ$ helps
with the analysis of various stochastic evolution equations.

 \begin{example}
 \label{ex4}
 {\rm
 Consider the   $w(H^1_2(\bR),H^{-1}_2(\bR))$
  Wiener Chaos solution $u$ of equation
  \begin{equation}
  \label{eq:classical55}
du(t,x)=(au_{xx}(t,x)+f(t,x))dt+\sigma u_x(t,x)dw(t),\ 0<t\leq T,\
x\in \bR,
\end{equation}
with $f \in L_2(\Omega\times(0,T); H^{-1}_2(\bR))$, $g\in
L_2(\Omega\times(0,T); L_2(\bR))$, and $u|_{t=0}=u_0\in L_2(\bR)$.
   Assume that $\sigma>0$ and define the sequence
  $Q$ so that $q_k=q$ for all $k\geq 1$ and $q< \sqrt{2a}/\sigma$. By Theorem
  \ref{th:ParabCl1}, equation
$$
dv=(av_{xx}+f)dt+(q\sigma u_x + g)dw
$$
 with  $v|_{t=0}=u_0$, has a unique traditional solution
$$
v \in L_{2}\left(\mbW; L_2((0,T);H^1_2(\bR) )\right)\bigcap
L_2\left(\mbW; \bC((0,T); L_2(\bR))\right).
$$
By Proposition \ref{prop:QT}, the $w(H^{1}_2(\bR), H^{-1}_2(\bR))$
Wiener Chaos
 solution $u$ of equation (\ref{eq:classical55}) satisfies
 $u=\cQ^{-1}v$ and
$$
u \in L_{2,Q}\left(\mbW; L_2((0,T);H^1_2(\bR) )\right)\bigcap
L_{2,Q}\left(\mbW; \bC((0,T); L_2(\bR))\right).
$$
Note that if equation (\ref{eq:classical55}) is strongly
parabolic, that is, $2a>\sigma^2$, then the weight $q$ can be
taken bigger than one, and, according  to the first statement of
Proposition \ref{prop:QT},
   regularity of the solution is better than the one guaranteed
by Theorem \ref{th:ParabCl1}. }
\end{example}

\begin{example}
\label{ex5} {\rm The  Wiener Chaos solutions can be constructed
for stochastic ordinary differential equations. Consider, for
example, \begin{equation}\label{eq:ODE55}
u(t)=1+\int_0^t\sum_{k\geq 1} u(s)dw_k(s),
\end{equation}
which clearly does not have a traditional  solution. On the other
hand, the unique  $w(\bR,\bR)$ Wiener Chaos solution of this
equation
 belongs to
$L_{2,Q}\left(\mbW; L_2((0,T)\right)$ for every $Q$ satisfying
$\sum_{k}q_k^2 < \infty$. Indeed, for (\ref{eq:ODE55}),  equation
(\ref{eq:clas}) becomes
$$
v(t)=1+\int_0^t\sum_{k} v(s)q_kdw_k(s).
$$
If $\sum_{k}q_k^2 < \infty$, then the traditional
  solution of this equation exists   and
 belongs to $L_2\left(\mbW; L_2((0,T))\right)$.
}
\end{example}

There exist equations for which the
 Wiener Chaos solution does not belong to any weighted
Wiener chaos space $L_{2,Q}$. An example is given below in Section
\ref{sec:.:frd}.

To define the S-transform, consider the following analog of the
stochastic exponential (\ref{eq:StrWN2}).
\begin{lemma}
\label{lm:StExp} If $h\in \cD\left(L_2((0,T); Y)\right)$ and
$$
\cE(h)=\exp\left(
 \int_0^T( h(t), dW(t))_Y-\frac{1}{2} \int_0^T\|h(t)\|_Y^2dt\right),
 $$
 then
 \begin{itemize}
 \item $\cE(h)\in L_{2,Q}(\mbW) $ for every sequence $Q$.
 \item $\cE(h)\in (\cS)_{\rho, \gamma}$ for $0\leq \rho<1$ and $\gamma\geq 0$.
 \item $\cE(h)\in (\cS)_{1,\gamma}$, $\gamma\geq 0$,  as long as
  $\|h\|_{L_2((0,T);Y)}^2$ is sufficiently small.
 \end{itemize}
 \end{lemma}
\begin{proof} Recall that, if $h\in \cD(L_2((0,T); Y))$, then
$h(t)=\sum_{i,k\in I_h} h_{k,i} m_i(t)y_k$, where $I_h$ is a
finite set. Direct computations show that
$$
\cE(h)=\prod_{i,k} \left( \sum_{n\geq 0}
\frac{H_n(\xi_{ik})}{n!}(h_{k,i})^n \right)= \sum_{\alpha \in \cJ}
\frac{h^{\alpha}}{\sqrt{\alpha!}} \xi_{\alpha}
$$
where $h^{\alpha}  = \prod_{i,k}h_{k,i}^{\alpha_i^k}$. In
particular,
\begin{equation}
\label{eq:cE} (\cE(h))_{\alpha}=\frac{h^{\alpha}}{\sqrt{\alpha!}}.
\end{equation}
Consequently, for every sequence $Q$ of positive numbers,
\begin{equation}
\label{eq:univ} \|\cE(h)\|_{L_{2,Q}(\mbW)}^2 = \exp\left(
\sum_{i,k \in I_h} {h_{k,i}^2}{q^2_k} \right) < \infty.
\end{equation}
Similarly, for $0\leq \rho<1$ and $\gamma\geq 0$,
\begin{equation}
\label{eq:univ55} \|\cE(h)\|_{(\cS)_{\rho, \gamma}}^2 =
\sum_{\alpha \in
\cJ}\prod_{i,k}\frac{((2ik)^{\gamma}h_{k,i})^{2\alpha^k_i}}
{(\a^k_i!)^{1-\rho}}= \prod_{i,k\in I_h}\left(\sum_{n\geq
0}\frac{((2ik)^{\gamma} h_{k,i})^{2n}}{(n!)^{1-\rho}} \right)
<\infty,
\end{equation}
and, for $\rho=1$,
\begin{equation}
\label{eq:univ66} \|\cE(h)\|_{(\cS)_{1, \gamma}}^2 = \sum_{\alpha
\in \cJ}\prod_{i,k}((2ik)^{\gamma}h_{k,i})^{2\alpha^k_i} =
\prod_{i,k\in I_h}\left(\sum_{n\geq
0}((2ik)^{\gamma}h_{k,i})^{2n}\right) <\infty,
\end{equation}
if $2\left(\max_{(m,n)\in
I_h)}(mn)^{\gamma}\right)\sum_{i,k}h_{k,i}^2 < 1$. Lemma
\ref{lm:StExp} is proved.
\end{proof}

\begin{remark}\mbox{}\label{rm:dense}
 It is well-known {\rm(}see, for example,
\cite[Proof of Theorem 5.5]{LSh_st}{\rm )} that the family
$\{\cE(h), h\in \cD\left(L_2((0,T); Y)\right)\}$ is dense in
$L_2(\mbW)$ and consequently in every $L_{2,Q}(\mbW)$ and every
$(\cS)_{\rho, \gamma}$,  $-1<\rho\leq 1$,  $\gamma\in \bR$.
\end{remark}

\begin{definition}
\label{def:str} If  $u \in L_{2,Q}(\mbW; X)$ for some $Q$, or if
$u\in \bigcup_{q\geq 0} (\cS)_{-\rho,-\gamma}(X)$, $0\leq \rho\leq
1$, then the deterministic function
\begin{equation}
\label{eq:str55} Su(h)=\sum_{\alpha \in
\cJ}\frac{u_{\a}h^{\a}}{\sqrt{\a!}}\in X
\end{equation}
  is called
 the {\bf S-transform} of $u$. Similarly,
 for $g\in \cD'\left(Y;L_{2,Q}(\mbW; X)\right)$ the S-transform
 $Sg(h) \in \cD'(Y; X)$ is defined by  setting $(Sg(h))_k=(Sg_k)(h)$.
\end{definition}

Note that if $u \in L_2(\mbW;X)$, then $Su(h)=\mbE(u\cE(h))$. If
$u$ belongs to  $L_{2,Q}(\mbW; X)$ or to  $ \bigcup_{q\geq 0}
(\cS)_{-\rho,-\gamma}(X)$, $0\leq \rho < 1$, then $Su(h)$ is
defined for all $ h \in \cD\left(L_2((0,T); Y)\right).$ If $u\in
\bigcup_{\gamma\geq 0}(\cS)_{-1,-\gamma}(X)$, then $Su(h)$ is
defined only for $h$ sufficiently close to zero.

By Remark \ref{rm:dense}, an  element $u$ from  $L_{2,Q}(\mbW; X)$
or $ \bigcup_{\gamma\geq 0} (\cS)_{-\rho,-\gamma}(X)$, $0\leq \rho
< 1$, is uniquely determined by the collection of deterministic
functions $Su(h),\; h \in \cD\left(L_2((0,T); Y)\right).$ Since
$\cE(h)>0$ for all $h \in \cD\left(L_2((0,T); Y)\right)$, Remark
\ref{rm:dense} also suggests the following definition.

\begin{definition}
\label{def:pos} An element $u$ from $ L_{2,Q}(\mbW)$ or $
\bigcup_{\gamma\geq 0} (\cS)_{-\rho,-\gamma}$, $0\leq \rho < 1$ is
called non-negative  ($u\geq 0$) if and only if $Su(h)\geq 0$ for
all $h\in \cD\left(L_2((0,T); Y)\right)$.
\end{definition}

The definition of the operator $\cQ$ and Definition \ref{def:pos}
imply the following result.

\begin{proposition}
\label{prop:pos} A generalized random element  $u$ from $
L_{2,Q}(\mbW)$ is non-negative
 if and only if $\cQ u\geq 0$.
\end{proposition}

For example, the solution of equation (\ref{eq:ODE55}) is
non-negative because
$$
\cQ u(t)=\exp\left(\sum_{k\geq 1}(q_kw_k(t)-(1/2)q_k^2)\right).
$$

We conclude this section with one technical remark.

Definition \ref{def:str} expresses the S-transform in terms of the
generalized Fourier coefficients. The following results makes it
possible to recover  generalized Fourier coefficients from the
corresponding S-transform.
\begin{proposition}
\label{prop:SFK55} If $u$  belongs to some $ L_{2,Q}(\mbW;X)$ or $
\bigcup_{\gamma\geq 0} (\cS)_{-\rho,-\gamma}(X)$, $0\leq \rho \leq
1$, then
\begin{equation}
\label{eq:diffFK} u_{\alpha} =\frac{1}{\sqrt{\a!}}
\left.\left(\prod_{i,k} \frac{\partial^{\alpha^k_i}Su(h)}{\partial
h_{k,i}^{\alpha_i^k}} \right)\right|_{h=0}.
\end{equation}
\end{proposition}
\begin{proof} For each  $\alpha\in \cJ$ with $K$ non-zero entries,
equality (\ref{eq:str55}) and Lemma \ref{lm:StExp} imply that the
function $Su(h)$, as a function of $K$ variables $h_{k,i}$, is
analytic in some neighborhood of zero. Then (\ref{eq:diffFK})
follows after differentiation
 of the series (\ref{eq:str55}).
\end{proof}

\section{General Properties of the Wiener Chaos Solutions}
\label{secPWC}
\setcounter{equation}{0}
\setcounter{theorem}{0}

Using notations and assumptions from Section \ref{secGS}, consider
the linear evolution equation
\begin{equation}
\label{eq:gen_eq55} du(t)=(\cA u(t) + f(t))dt+(\cM u(t) + g(t) ,
dW(t))_Y,\ 0<t\leq T,\ u|_{t=0}=u_0.
\end{equation}
The objective of this section is to study how  the Wiener Chaos
compares with the traditional and white noise solutions.

To make the presentation shorter, call an X-valued  generalized
random element {\em S-admissible} if and only if it belongs to
$L_{2,Q}(\cF^W;X)$ for some $Q$ or to $(\cS)_{\rho,q}(X)$ for some
$\rho\in [-1,1]$ and $q \in \bR$. It was shown in Section
\ref{secWS} that, for every S-admissible $u$,
 the S-transform $Su(h)$
  is defined when $h=\sum_{i,k}h_{k,i}m_iy_k\in \cD(L_2((0,T);Y))$ and is
   an analytic function of $h_{k,i}$  in some neighborhood of $h=0$.

 The next result describes the S-transform of the  Wiener Chaos solution.

\begin{theorem}
\label{th:SW1} Assume that
\begin{enumerate}
\item there exists a unique  $w(A,X)$ Wiener Chaos solution $u$ of
(\ref{eq:gen_eq55}) and $u$ is $S$-admissible; \item For each
$t\in [0,T]$, the linear operators $\cA(t), \cM_k(t)$ are bounded
from $A$ to $X$; \item  the generalized random elements
 $u_0,  f, g_k$ are S-admissible.
 \end{enumerate}
  Then, for
 every $h\in \cD(L_2((0,T);Y))$ with $\|h\|_{L_2((0,T);Y)}^2$  sufficiently small,
 the function $v=Su(h)$ is a   $w(A,X)$ solution of the deterministic equation
\begin{equation}
\label{eq:S} v(t) = Su_0(h) + \int_0^t \Big( \cA v + Sf(h)+ (\cM_k
v + Sg_k(h))h_k\Big)(s)ds.
\end{equation}
\end{theorem}
\begin{proof} By assumption, $Su(h)$ exists for suitable functions
$h$. Then the S-transformed equation (\ref{eq:S})  follows from
the definition of the S-transform (\ref{eq:str55}) and the
propagator equation (\ref{eq:def2}) satisfied by the generalized
Fourier coefficients of $u$. Indeed, continuity of operator $\cA$
implies
$$
S(\cA u)(h)=\sum_{\alpha}\frac{h^{\a}}{\sqrt{\a!}}\cA u_{\alpha}=
\cA \sum_{\alpha}\frac{h^{\a}}{\sqrt{\a!}} u_{\alpha}=\cA(Su(h)).
$$
Similarly,
\begin{equation*}
\begin{split}
\sum_{\a}\frac{h^{\a}}{\sqrt{\a!}}\sum_{i,k}\sqrt{\a^k_i}\cM_k
u_{\alpha^-(i,k)}m_i=
\sum_{\a}\sum_{i,k}\frac{h^{\a^-(i,k)}}{\sqrt{\a^-(i,k)!}}
\cM_ku_{\alpha^-(i,k)}m_ih_{k,i}\\
= \sum_{i,k}\left(\sum_{\alpha}\frac{h^{\a}}{\sqrt{\a}}
\cM_ku_{\alpha}\right)m_ih_{k,i}=\cM_k(Su(h))h_k.
\end{split}
\end{equation*}
Computations for the other terms are similar. Theorem \ref{th:SW1}
is proved.
\end{proof}

\begin{remark}
\label{rm:ExpEvol} If $h\in \cD(L_2((0,T); Y))$ and
\begin{equation}
\label{eq:ExpEvol1} \cE_t(h)=\exp\left(\int_0^t(h(s), dW(s))_Y -
\frac{1}{2}\int_0^t \|h(t)\|_Y^2dt \right),
\end{equation}
then, by the It\^{o} formula,
\begin{equation}
\label{eq:ExpEvol} d\cE_t(h)=\cE_t(h)(h(t), dW(t))_Y.
\end{equation}
If $u_0$ is deterministic, $f$ and $g_k$ are $\cF^W_t$-adapted,
and $u$ is a square-integrable  solution of (\ref{eq:gen_eq55}),
then equality (\ref{eq:S}) is  obtained by multiplying  equations
(\ref{eq:ExpEvol}) and (\ref{eq:gen_eq55}) according to the
It\^{o} formula and taking the expectation.
\end{remark}

\begin{remark}
\label{rm:xi(t)} Rewriting (\ref{eq:ExpEvol}) as
$$
d\cE_t(h)=\cE_t(h)h_{k,i}m_i(t)dw_k(t)
$$
and using the relations
$$
\cE_t(h)=\mbE(\cE_T(h)|\cF^W_t),\ \xi_{\a}=\frac{1}{\sqrt{\a!}}
\left.\left(\prod_{i,k}
\frac{\partial^{\alpha^k_i}\cE_T(h)}{\partial
h_{k,i}^{\alpha_i^k}} \right)\right|_{h=0},
$$
we arrive at representation  (\ref{eq:xi(t)}) for
$\mbE(\xi_{\a}|\cF^W_t)$.
\end{remark}
A partial converse of Theorem \ref{th:SW1} is that,  under some
regularity conditions, the Wiener Chaos solution can be recovered
from the solution of the S-transformed equation (\ref{eq:S}).

\begin{theorem}
\label{th:SW111} Assume that the linear operators $\cA(t),\
\cM_k(t)$, $t\in [0,T]$,  are bounded from $A$ to $X$, the input
data
 $u_0$, $f$, $g_k$ are S-admissible, and,   for
 every $h\in \cD(L_2((0,T);Y))$ with $\|h\|_{L_2((0,T);Y)}^2$  sufficiently small,
 there exists a  $w(A,X)$ solution $v=v(t;h)$ of  equation (\ref{eq:S}).
 We write $h=h_{k,i}m_iy_k$ and consider $v$ as a function of the variables
 $h_{k,i}$. Assume that all  the derivatives of $v$   at the point
 $h=0$ exists,  and, for $\a \in \cJ$,  define
\begin{equation}
\label{eq:diffFK111} u_{\alpha}(t) =\frac{1}{\sqrt{\a!}}
\left.\left(\prod_{i,k}
\frac{\partial^{\alpha^k_i}v(t;h)}{\partial h_{k,i}^{\alpha_i^k}}
\right)\right|_{h=0}.
\end{equation}
Then the generalized random process $u(t)=\sum_{\alpha \in
\cJ}u_{\a}(t)\xi_{\a}$ is a $w(A,X)$ Wiener Chaos solution of
(\ref{eq:gen_eq55}).
\end{theorem}

\begin{proof} Differentiation of   (\ref{eq:S}) and application of
 Proposition \ref{prop:SFK55} show that the functions  $u_{\alpha}$
 satisfy the propagator (\ref{eq:def2}).
 \end{proof}

\begin{remark}
The central part in the construction of the white noise solution
of (\ref{eq:gen_eq55}) is proving that the solution of
(\ref{eq:S}) is an S-transform of a suitable generalized random
process. For many particular cases of equation
(\ref{eq:gen_eq55}), the corresponding analysis is carried out in
\cite{HKPS,HOUZ,MR,Pot}. The consequence  of Theorems \ref{th:SW1}
and \ref{th:SW111} is  that a white noise solution of
(\ref{eq:gen_eq55}), if exists, must coincide with the Wiener
Chaos solution.
\end{remark}

The next theorem establishes the connection between the Wiener
Chaos solution and the traditional solution. Recall that the
traditional, or square-integrable,
 solution of (\ref{eq:gen_eq55}) was introduced in
  Definition \ref{def:ParabCl}. Accordingly, the notations from Section
  \ref{sec:.:clas} will be used.

\begin{theorem}
\label{th:WC-trad} Let $(V,H,V')$ be a normal triple of Hilbert
spaces. Take a deterministic function  $u_0$ and $\cF^W_t$-adapted
random processes function, $f$ and $g_k$  so that
(\ref{eq:InputCl}) holds. Under these assumptions we have the
following two statements.
\begin{enumerate}
\item An $\cF^W_t$-adapted  traditional solution of
(\ref{eq:gen_eq55}) is also a Wiener Chaos solution. \item If  $u$
is a $w(V,V')$ Wiener Chaos solution of (\ref{eq:gen_eq55}) so
that
\begin{equation}
\label{eq:WC-trad1} \sum_{\alpha\in
\cJ}\left(\int_0^T\|u_{\a}(t)\|_V^2dt +\sup_{0\leq t\leq
T}\|u_{\a}(t)\|_H^2\right)< \infty,
\end{equation}
then $u$ is an $\cF^W_t$-adapted  traditional solution of
(\ref{eq:gen_eq55}).
\end{enumerate}
\end{theorem}

\begin{proof} (1) If $u=u(t)$ is an  $\cF^W_t$-adapted traditional
solution, then
$$
u_{\a}(t)=\mbE(u(t)\xi_{\alpha})=\mbE\left(u(t)\mbE(\xi_{\a}|\cF^W_t)\right)=
\mbE (u(t)\xi_{\a}(t)).
$$
Then  the propagator  (\ref{eq:def2}) for $u_{\a}$ follows after
applying the It\^{o} formula to the product $u(t)\xi_{\a}(t)$ and
using (\ref{eq:xi(t)}).

(2) Assumption (\ref{eq:WC-trad1}) implies
$$
u \in L_2(\Omega\times (0,T); V)\bigcap L_2(\Omega; \bC((0,T);
H)).
$$
Then, by Theorem \ref{th:SW1},  for every $\varphi\in V$ and $h
\in \cD((0,T); Y)$, the S-transform $u_h$ of $u$ satisfies
\begin{equation*}
\begin{split}
(u_h(t),\varphi)_H&=(u_0,\varphi)_H+\int_0^t\langle \cA u_h(s),
\varphi \rangle ds + \int_0^t\langle f(s),\varphi\rangle ds \\
&+\sum_{\alpha\in \mathcal{J}}\frac{h^{\alpha}}{\alpha!} \sum_{i,k}
\int_0^t \sqrt{\alpha^k_i}m_i(s)\big((\cM_k
u_{\alpha^-(i,k)}(s),\varphi)_H \\
&+(g_k(s),
\varphi)_HI(|\alpha|=1)\big)ds.
\end{split}%
\end{equation*}
If  $I(t)=\int_0^t(\cM_k u(s), \varphi)_Hdw_k(s)$, then
\begin{equation}
 \label{eq:int1H}
\mathbb{E} (I(t)\xi_{\alpha}(t))= \int_0^t\sum_{i,k} \sqrt{\alpha_i^k}%
m_i(s)( \cM_k u_{\alpha^-(i,k)}(s),\varphi)_H ds.
\end{equation}
Similarly,
$$
\mathbb{E}
\left(\xi_{\alpha}(t)\int_0^t(g_k(s),\varphi)_Hdw_k(s)\right)=
\sum_{i,k} \int_0^t \sqrt{\alpha^k_i}m_i(s)(g_k(s),
\varphi)_HI(|\alpha|=1)ds.
$$
Therefore,
\begin{equation*}
\begin{split}
\sum_{\alpha\in \mathcal{J}}\frac{h^{\alpha}}{\alpha!} \sum_{i,k}
\int_0^t
\sqrt{\alpha^k_i}m_i(s)(\cM_k u_{\alpha^-(i,k)}(s),\varphi)_Hds\\
= \mathbb{E}\left(\mathcal{E}(h)\int_0^t\left( (\cM_k
u(s),\varphi)_H+ (g_k(s),\varphi)_H\right) dw_k(s)\right).
\end{split}
\end{equation*}
As a result,
\begin{equation}
 \label{eq:str5H}
\begin{split}
\mathbb{E}\left(\mathcal{E}(h)(u(t), \varphi)_H\right)&=
\mbE\left(\cE(h)(u_0,\varphi)_H\right)\\
&+\mathbb{E}\left(\mathcal{E}(h)\int_0^t\langle \cA u(s), \varphi\rangle ds \right)
+\mathbb{E}\left(\mathcal{E}(h)\int_0^t\langle f(s), \varphi
\rangle ds\right) \\
&+  \mathbb{E}\left(\mathcal{E}(h)\int_0^t\left(
(\cM_k u(s),\varphi)_H + (g_k(s),\varphi)_H\right) dw_k(s)\right).
\end{split}%
\end{equation}
Equality (\ref{eq:str5H}) and Remark \ref{rm:dense}
  imply that, for each $t$ and each $\varphi$, (\ref{eq:trad-def})
  holds with probability one. Continuity of $u$ implies that, for each $\varphi$,
   a single  probability-one set can be chosen for all $t\in [0,T]$.
   Theorem \ref{th:linH2} is proved. \end{proof}

\section{Regularity of the Wiener Chaos Solution}
\label{secLH}
\setcounter{equation}{0}
\setcounter{theorem}{0}

Let $\mbF=(\Omega, \cF, \{\cF_t\}_{t\geq 0}, \mbP)$ be a
stochastic basis with the usual assumptions and $w_k=w_k(t), \
k\geq 1,\ t\geq 0$, a collection of standard Wiener processes on
$\mbF$. As in Section \ref{sec:.:clas}, let $(V, H, V')$ be a
normal triple of Hilbert spaces and
 $\cA(t): V\to V'$,
$\cM_k(t): V\to H$, linear bounded operators; $t\in [0,T]$.

In this section we study  the linear equation
\begin{equation}
\label{eq:linH} u(t)=u_0+\int_0^t (\cA u(s)+f(s))ds + \int_0^t
(\cM_k u (s)+g_k(s))
 dw_k,\ 0\leq t\leq T,
\end{equation}
under the following {\bf {assumptions}}:
\begin{enumerate}
\item[\textbf{{A1}}] There exist positive numbers $C_1$ and
$\delta$ so that
\begin{equation}
\label{eq:linHa1} \langle\cA(t) v, v\rangle + \delta \|v\|_V^2
\leq C_1\|v\|_H^2,\ v\in V,\ t\in [0,T].
\end{equation}
\item[\textbf{{A2}}] There exists a real number $C_2$ so that
\begin{equation}
\label{eq:linHa2} 2\langle\cA(t) v, v\rangle +\sum_{k\geq
1}\|\cM_k(t) v\|_H^2 \leq C_2\|v\|^2_H, \ v\in V,\ t\in [0,T].
\end{equation}
\item[\textbf{A3}] The initial condition $u_0$ is non-random and
belongs to  $ H$; the process  $f=f(t)$ is  deterministic and
$\int_0^T\|f(t)\|_{V'}^2dt < \infty$; each $g_k=g_k(t)$ is a
deterministic processes and $\sum_{k\geq 1}
\int_0^T\|g_k(t)\|_H^2dt < \infty$.
\end{enumerate}

Note that condition (\ref{eq:linHa2}) is weaker than
(\ref{eq:ParabCll}). Traditional analysis of equation
(\ref{eq:linH}) under (\ref{eq:linHa2}) requires additional
regularity assumptions on the input data and additional Hilbert
space constructions beyond the normal triple \cite[Section
3.2]{Roz}. In particular, no existence of a traditional solution
is known under assumptions \textbf{A1}-\textbf{A3}, and
 the Wiener chaos approach provides new existence and regularity results for
equation (\ref{eq:linH}). A different version of the following
theorem is presented in \cite{LR_AP}.

\begin{theorem}
\label{th:linH1} Under assumptions {\textbf{A1}}--{\textbf{A3}},
 for every $T>0$,  equation (\ref{eq:linH})
has a unique $w(V,V')$ Wiener Chaos solution.
 This  solution $u=u(t)$ has the following properties:
\begin{enumerate}
\item There exists a weight sequence $Q$ so that
$$u
 \in L_{2,Q}(\mbW; L_2((0,T); V))\bigcap L_{2,Q}(\mbW; \bC((0,T); H)).
 $$
\item For every $0\leq t\leq T$, $u(t) \in L_2(\Omega; H)$
 and
\begin{equation}
\label{eq:linHn} \mbE \|u(t)\|_H^2 \leq
3e^{C_2t}\left(\|u_0\|_{H}^2+C_f\int_0^t\|f(s)\|_{V'}^2ds
+\sum_{k\geq 1} \int_0^t\|g_k(s)\|_H^2ds \right),
\end{equation}
where the number $C_2$ is from (\ref{eq:linHa2}) and the positive
number $C_f$ depends only on $\delta$ and $C_1$ from
(\ref{eq:linHa1}). \item For every $0\leq t\leq T$,
\begin{equation}
\label{eq:iter_int1}%
\begin{split}
u(t)=u_{(0)}\;+\; &  \sum_{n\geq 1}\sum_{k_{1},
\ldots, k_{n}\geq1}\int_{0}^{t}\int_{0}^{s_{n}}\ldots\int_{0}^{s_{2}}\\
& \!\!\!\!\!\!\!\!\!\!\!\! \Phi_{t,s_{n}}{{\mathcal{M}}}_{k_{n}}\cdots
\Phi_{s_{2},s_{1}}\left( {{\mathcal{M}}}_{k_{1}}u_{(0)}+
g_{k_{1}}(s_{1})\right) dw_{k_{1}}(s_{1})\cdots dw_{k_{n}}(s_{n}),
\end{split}
\end{equation}
where $\Phi_{t,s}$ is the semi-group of  the operator $\cA$.
\end{enumerate}
\end{theorem}
\begin{proof} Assumption {\textbf{A2}} and the properties of the normal
triple imply that there exists a positive  number $C^*$ so that
\begin{equation}
\label{eq:contMk} \sum_{k\geq 1}\|\cM_k (t) v\|_H^2 \leq C^*
\|v\|_V^2, \ v \in V, \ t\in [0,T].
\end{equation}
Define the sequence $Q$ so that
\begin{equation}
q_k=\left(\frac{\mu\delta}{C^*}\right)^{1/2}:=q,\ k\geq 1,
\end{equation}
where $\mu\in (0,2)$ and
 $\delta$ is from Assumption {\textbf{A1}}. Then, by Assumption {\textbf{A2}},
\begin{equation}\label{eq:Roz}
2 \langle Av, v \rangle + \sum_{k\geq 1}
q^2\|\cM_k v\|_H^2 \leq -(2-\mu)\delta\|v\|_V^2 + C_1\|v\|_H^2.
\end{equation}
It follows from Theorem \ref{th:ParabCl} that equation
\begin{equation}\label{eq:Roz1}
v(t)=u_0+\int_0^t (\cA v + f)(s)ds
+\sum_{k\geq 1} \int_0^t q(\cM_k v + g_k)(s)dw_k(s)
\end{equation}
has a unique solution
$$
v \in L_2(\mbW; L_2 ((0,T); V))\bigcap L_2(\mbW; \bC((0,T); H)).
$$
Comparison of the propagators for equations (\ref{eq:linH}) and
(\ref{eq:Roz1}) shows that $u=\cQ^{-1} v$ is the unique $w(V,V')$
solution of (\ref{eq:linH}) and \begin{equation}\label{eq:regLQ} u
\in L_{2,Q}(\mbW; L_2 ( (0,T); V))\bigcap L_{2,Q}(\mbW; \bC((0,T);
H)). \end{equation}

If $C^* < 2\delta$, then equation (\ref{eq:linH}) is strongly
parabolic and
 $q>1$ is an admissible choice of the weight. As a result, for
 strongly parabolic equations, the result (\ref{eq:regLQ}) is  stronger than
 the conclusion of Theorem \ref{th:ParabCl}.

The proof of (\ref{eq:linHn}) is based on the analysis of the
propagator
\begin{equation}
\label{eq:def2h}
\begin{split}
 u_{\alpha}(t)&=u_{0}I(|\alpha|=0) + \int_0^t \Big(\cA u_{\alpha}(s)+f(s)I(|\alpha|=0)
 \Big)ds\\
&+ \int_0^t\sum_{i,k}\sqrt{\alpha_i^k}(\cM_k
u_{\alpha^-(i,k)}(s)+g_k(s)I(|\alpha|=1)) m_i(s)ds.
\end{split}
\end{equation}
We  consider three particular cases: (1) $f=g_k=0$ (the
homogeneous equation); (2) $u_0=g_k=0$; (3) $u_0=f=0$. The general
case  will then follow by linearity and the triangle inequality.

Denote by  $(\Phi_{t,s}, \ t\geq s\geq 0)$ the semi-group
generated by the operator $\cA(t)$; $\Phi_t:=\Phi_{t,0}$. One of
the consequence of Theorem \ref{th:ParabCl}  is that, under
Assumption {\textbf{A1}}, this semi-group exists and is strongly
continuous in $H$.

 Consider the homogeneous equation: $f=g_k=0$.
 By Corollary \ref{cor:ind},
 \begin{equation}
\label{eq:L2norm}
\sum_{|\alpha|=n} \|u_\alpha(t)\|_H^2 =\!\!\!\!
 \sum_{k_1, \ldots, k_n\geq 1} \int_0^t \int_0^{s_n}\!\!\!\!\cdots\!\!\int_0^{s_2}
 \!\!\! \|\Phi_{t,s_n}\cM_{k_n}\!\cdots \Phi_{s_2,s_1}\cM_{k_1}\Phi_{s_1}u_0\|_H^2ds^n,
 \end{equation}
 where $ds^n=ds_1\ldots ds_n$.
 Define  $F_n(t)=\sum_{|\alpha|=n} \|u_\alpha(t)\|_{H}^2$, $n\geq 0$.
Direct application of (\ref{eq:linHa2}) shows that
\begin{equation}
\frac{d}{dt}F_0(t)\leq C_2 F_0(t) - \sum_{k\geq
1}\|\cM_k\Phi_{t}u_0\|_H^2.
\end{equation}
  For $n\geq 1$,  equality (\ref{eq:L2norm}) implies
\begin{equation}
\label{eq:LinHp1}
\begin{split}
\frac{d}{dt}F_n(t)=\sum_{k_1, \ldots, k_n\geq 1} \int_0^t
\int_0^{s_{n-1}}\cdots\int_0^{s_2}\!\!
 \|\cM_{k_n} \Phi_{t,s_{n-1}}\cdots
 \cM_{k_1} \Phi_{s_1}u_0\|_{H}^2ds^{n-1}\\
 +\sum_{k_1, \ldots, k_n\geq 1}\int_0^t \int_0^{s_n}\ldots\int_0^{s_2}
 \langle\cA \Phi_{t,s_n}\cM_{k_n} \ldots
  \Phi_{s_1}u_0, \Phi_{t,s_n}\cM_{k_n} \ldots
  \Phi_{s_1}u_0\rangle ds^n.
 \end{split}
 \end{equation}
By (\ref{eq:linHa2}),
\begin{equation}
\label{eq:LinHp2}
\begin{split}
\sum_{k_1, \ldots, k_n\geq 1}\int_0^t
\int_0^{s_n}\ldots\int_0^{s_2}
 \langle\cA \Phi_{t,s_n}\cM_{k_n} \ldots
  \Phi_{s_1}u_0, \Phi_{t,s_n}\cM_{k_n} \ldots
  \Phi_{s_1}u_0\rangle ds^n\\
 \leq
 -\sum_{k_1, \ldots, k_{n+1}\geq 1}
\int_0^t \int_0^{s_{n}}\ldots\int_0^{s_2}
 \|\cM_{k_{n+1}} \Phi_{t,s_{n}}\cM_{k_{n}} \ldots
 \cM_{k_1} \Phi_{s_1}u_0\|_{H}^2ds^{n}\\
 + C_2\sum_{k_1, \ldots, k_n\geq 1}
\int_0^t \int_0^{s_{n}}\ldots\int_0^{s_2}
 \| \Phi_{t,s_{n}}\cM_{k_{n}} \ldots
 \cM_{k_1} \Phi_{s_1}u_0\|_{H}^2ds^{n}.
 \end{split}
 \end{equation}

As a result, for $n\geq 1$,
\begin{equation}
\label{eq:LinHp3}
\begin{split}
&\frac{d}{dt}F_n(t) \leq C_2F_n(t)\\
+&\sum_{k_1, \ldots, k_n\geq 1} \int_0^t
\int_0^{s_{n-1}}\ldots\int_0^{s_2}
 \|\cM_{k_n} \Phi_{t,s_{n-1}}\cM_{k_{n-1}} \ldots
 \cM_{k_1} \Phi_{s_1}u_0\|_{H}^2ds^{n-1}\\
-&\sum_{k_1, \ldots, k_{n+1}\geq 1} \int_0^t
\int_0^{s_{n}}\ldots\int_0^{s_2}
 \|\cM_{k_{n+1}} \Phi_{t,s_{n}}\cM_{k_{n}} \ldots
 \cM_{k_1} \Phi_{s_1}u_0\|_{H}^2ds^{n}.
 \end{split}
\end{equation}
Consequently,
\begin{equation}
\label{eq:LinHp4} \frac{d}{dt} \sum_{n=0}^N \sum_{|\alpha|=n}
\|u_{\alpha}(t)\|_H^2 \leq C_2\sum_{n=0}^N \sum_{|\alpha|=n}
\|u_{\alpha}(t)\|_H^2,
\end{equation}
so that,  by the Gronwall inequality,
\begin{equation}
\label{eq:linHun} \sum_{n=0}^N \sum_{|\alpha|=n}
\|u_{\alpha}(t)\|_H^2 \leq e^{C_2t}\|u_0\|_{H}^2
\end{equation}
or
\begin{equation}
\label{eq:linHnu} \mbE \|u(t)\|_H^2 \leq e^{C_2t}\|u_0\|_{H}^2.
\end{equation}

Next, let us assume that $u_0=g_k=0$.
 Then
the propagator  (\ref{eq:def2h}) becomes
\begin{equation}
\label{eq:def2nhf}
 u_{\alpha}(t)= \int_0^t (\cA u_{\alpha}(s)+f(s)I(|\alpha|=0))ds
+ \int_0^t\sum_{i,k}\sqrt{\alpha_i^k}\cM_k
u_{\alpha^-(i,k)}(s)m_i(s)ds.
\end{equation}
Denote by  $u_{(0)}(t)$ the solution corresponding to
 $\alpha=0$. Note that
\begin{equation*}
\begin{split}
 \|u_{(0)}(t)\|_H^2 = 2  \int_0^t \langle \cA u_{(0)}(s),u_{(0)}(s)\rangle ds+
2\int_0^t\langle f(s),u_{(0)}(s)\rangle ds \\
\leq C_2\int_0^t\|u_{(0)}(s)\|_H^2ds- \int_0^t\sum_{k\geq 1}
\|\cM_k u_{(0)}(s)\|_H^2ds + C_f\int_0^t\|f(s)\|^2_{V'}ds.
\end{split}
\end{equation*}
By Corollary \ref{cor:ind},
\begin{equation}
\label{eq:L2normNH}
\sum_{|\alpha|=n} \|u_\alpha(t)\|_H^2 =\\
 \sum_{k_1, \ldots, k_{n}\geq 1} \int_0^t \int_0^{s_n}\ldots\int_0^{s_2}
 \|\Phi_{t,s_n}\cM_{k_n} \ldots \cM_{k_1}u_{(0)}(s_1)\|_H^2ds^n
 \end{equation}
for $n\geq 1$. Then, repeating the calculations
(\ref{eq:LinHp1})--(\ref{eq:LinHp3}), we conclude that
\begin{equation}
 \sum_{n=1}^N\sum_{|\alpha|=n}\|u_{\alpha}(t)\|_H^2
\leq C_f\int_0^t\|f(s)\|^2_{V'}ds +
C_2\int_0^t\sum_{n=1}^N\sum_{|\alpha|=n}\|u_{\alpha}(s)\|_H^2ds,
\end{equation}
and,  by the Gronwal inequality,
\begin{equation}
\label{eq:linHnf} \mbE \|u(t)\|_H^2 \leq
C_fe^{C_2t}\int_0^t\|f(s)\|_{V'}^2ds.
\end{equation}

Finally, let us assume that $u_0=f=0$. Then the propagator
(\ref{eq:def2h}) becomes
\begin{equation}
\label{eq:def2nh}
\begin{split}
 u_{\alpha}(t)&= \int_0^t \cA u_{\alpha}(s)ds
\\ &
+ \int_0^t\left(\sum_{i,k}\sqrt{\alpha_i^k}\cM_k
u_{\alpha^-(i,k)}(s) + g_k(s)I(|\alpha|=1)\right)m_i(s)ds.
\end{split}
\end{equation}
Even though  $u_{\alpha}(t)=0$ if $\alpha=0$, we have
\begin{equation}u_{(ik)}=\int_0^t\Phi_{t,s}g_k(s)m_i(s)ds, \end{equation} and
then the arguments from the proof of Corollary \ref{cor:ind}
apply, resulting  in
\begin{equation*}
\sum_{|\alpha|=n} \|u_\alpha(t)\|_H^2 =\\
 \sum_{k_1, \ldots, k_{n}\geq 1} \int_0^t \int_0^{s_n}\ldots\int_0^{s_2}
 \|\Phi_{t,s_n}\cM_{k_n} \ldots \Phi_{s_2,s_1}g_{k_1}(s_1)\|_H^2ds^n
 \end{equation*}
for $n\geq 1$. Note that
\begin{equation*}
\sum_{|\alpha|=1} \|u_\alpha(t)\|_H^2=  \sum_{k\geq 1} \int_0^t
\|g_k(s)\|^2_Hds+ 2\sum_{k\geq 1}\int_0^t\langle \cA
\Phi_{t,s}g_k(s), \Phi_{t,s}g_k(s)\rangle ds.
\end{equation*}
Then, repeating the calculations
(\ref{eq:LinHp1})--(\ref{eq:LinHp3}), we conclude that
\begin{equation}
 \sum_{n=1}^N\sum_{|\alpha|=n}\|u_{\alpha}(t)\|_H^2
\leq \sum_{k\geq 1} \int_0^t\|g_k(s)\|_H^2ds +
C_2\int_0^t\sum_{n=1}^N\sum_{|\alpha|=n}\|u_{\alpha}(s)\|_H^2ds,
\end{equation}
and, by  the Gronwal inequality,
\begin{equation}
\label{eq:linHng} \mbE \|u(t)\|_H^2 \leq e^{C_2t}\sum_{k\geq 1}
\int_0^t\|g_k(s)\|_H^2ds.
\end{equation}

To derive (\ref{eq:linHn}), it remains to
 combine (\ref{eq:linHnu}), (\ref{eq:linHnf}),
and (\ref{eq:linHng}) with the inequality $(a+b+c)^2 \leq
3(a^2+b^2+c^2)$.

Representation  (\ref{eq:iter_int1}) of the Wiener chaos  solution
as a sum of iterated It\^{o} integrals now follows from Corollary
\ref{cor:ind}.
 Theorem \ref{th:linH1} is proved.
\end{proof}

\begin{corollary}
\label{cor:LinH1} If $ \displaystyle \sum_{\alpha\in \cJ}
\int_0^T\|u_{\a}(s)\|_V^2ds < \infty$, then $\displaystyle\sum_{\a \in
\cJ}\sup_{0\leq t\leq T} \|u_{\a}(t)\|^2_H < \infty.$
\end{corollary}
\begin{proof} The proof of Theorem \ref{th:linH1} shows that it is
enough to consider the homogeneous equation. Then by inequalities
(\ref{eq:LinHp2})--(\ref{eq:LinHp3}),
\begin{equation}
\label{eq:111}
\begin{split}
&\sum_{\ell=n+1}^{n_1}\sum_{|\a|=\ell}\|u_{\alpha}(t)\|_H^2=\sum_{\ell=n+1}^{n_1}
F_{\ell}(t)\\
&\leq e^{C_2T} \sum_{k_1, \ldots, k_{n+1}\geq 1} \int_0^T\int_0^t
\int_0^{s_{n}}\ldots\int_0^{s_2}
 \|\cM_{k_{n+1}} \Phi_{t,s_{n}}\cM_{k_{n}} \ldots
  \Phi_{s_1}u_0\|_{H}^2ds^{n}dt.
 \end{split}
 \end{equation}
 By Corollary \ref{cor:ind},
 \begin{equation}
 \label{eq:333}
 \begin{split}
& \int_0^T\|u_{\a}(s)\|_V^2ds
\\ =
&\sum_{n\geq 1}\sum_{k_1, \ldots, k_n\geq 1} \int_0^T\int_0^t
\int_0^{s_{n}}\ldots\int_0^{s_2}
 \|\cM_{k_n} \Phi_{t,s_{n}}\cM_{k_{n}} \ldots
  \Phi_{s_1}u_0\|_{V}^2ds^{n}dt < \infty.
 \end{split}
 \end{equation}
 As  a result, (\ref{eq:contMk}) and (\ref{eq:333}) imply
 $$
 \lim_{n\to \infty}
 \int_0^T\int_0^t \int_0^{s_{n}}\ldots\int_0^{s_2}
 \|\cM_{k_{n+1}} \Phi_{t,s_{n}}\cM_{k_{n}} \ldots
 \cM_{k_1} \Phi_{s_1}u_0\|_{H}^2ds^{n}dt =0,
 $$
 which, by (\ref{eq:111}), implies uniform, with respect to $t$, convergence of
 the series $\sum_{\alpha\in \cJ}\|u_{\a}(t)\|_H^2$. Corollary \ref{cor:LinH1}
 is proved.
 \end{proof}

 \begin{corollary}
\label{cor:meas_coef} Let $a_{ij},b_{i},c,\sigma_{ik},\nu_{k}$ be
deterministic measurable functions of $(t,x)$ so that
\[
|a_{ij}(t,x)|+|b_{i}(t,x)|+|c(t,x)|+|\sigma_{ik}(t,x)|+|\nu_{k}(t,x)|\leq
K,
\]
$\ i,j=1,\ldots,d,\ k\geq1,\ x\in{\mathbb{R}}^{d},\ 0\leq t\leq
T;$
\[
\left(
a_{ij}(t,x)-\frac{1}{2}\sigma_{ik}(t,x)\sigma_{jk}(t,x)\right)
y_{i}y_{j}\geq0,
\]
$\ x,y\in{\mathbb{R}}^{d},\ 0\leq t\leq T;$ and
\[
\sum_{k\geq1}|\nu_{k}(t,x)|^{2}\leq C_{\nu}<\infty,
\]
$x\in{\mathbb{R}}^{d},\ 0\leq t\leq T.$ Consider the equation
\begin{equation}
du=(D_{i}(a_{ij}D_{j}u)+b_{i}D_{i}u+c\;u+f)dt+(\sigma_{ik}D_{i}u+\nu
_{k}u+g_{k})dw_{k}. \label{eq:meas_coef}%
\end{equation}
Assume that the input data satisfy $u_{0}\in
L_{2}({\mathbb{R}}^{d})$, $f\in
L_{2}((0,T);H_{2}^{-1}({\mathbb{R}}^{d})),$\newline$\sum_{k\geq1}\Vert
g_{k}\Vert_{L_{2}((0,T)\times{\mathbb{R}}^{d})}^{2}<\infty$, and
there exists an $\varepsilon>0$ so that
\[
a_{ij}(t,x)y_{i}y_{j}\geq\varepsilon|y|^{2},\
x,y\in{\mathbb{R}}^{d},\ 0\leq t\leq T.
\]
Then there exists a unique Wiener Chaos solution $u=u(t,x)$ of
(\ref{eq:meas_coef}). The solution has the following regularity:
\begin{equation}
u(t,\cdot)\in L_{2}({\mathbb{W}};L_{2}({\mathbb{R}}^{d})),\ 0\leq
t\leq T,
\label{eq:meas_coef-reg1}%
\end{equation}
and
\begin{equation}
\label{eq:meas_coef-reg2}
\begin{split}
\mathbb{E}\Vert u\Vert_{L_{2}({\mathbb{R}}^{d})}^{2}(t)&\leq
C^{\ast}\Big(
\Vert u_{0}\Vert_{L_{2}({\mathbb{R}}^{d})}^{2}+\Vert f\Vert_{L_{2}%
((0,T);H_{2}^{-1}({\mathbb{R}}^{d}))}^{2}\\
&+\sum_{k\geq1}\Vert
g_{k}\Vert
_{L_{2}((0,T)\times{\mathbb{R}}^{d})}^{2}\Big),
\end{split}
\end{equation}
where the positive number $C^{\ast}$ depends only on
$C_{\nu},K,T,$ and $\varepsilon$.
\end{corollary}

\begin{remark}
\label{rm:LinH1}
$ {               }  $\\
(1) If (\ref{eq:ParabCll}) holds instead of (\ref{eq:linHa2}),
then the proof of Theorem \ref{th:linH1}, in particular,
(\ref{eq:LinHp2})--(\ref{eq:LinHp3}), shows that the term
$\mbE\|u(t)\|^2_H$ in the left-hand-side of inequality
(\ref{eq:linHn}) can be replaced with
$$
\mbE\left(\|u(t)\|^2_H+\varepsilon\int_0^t\|u(s)\|_V^2ds\right).
$$
(2) If $f=g_k=0$ and the equation is fully degenerate, that is,
$2\langle \cA(t) v, v\rangle + \sum_{k\geq 1}\|\cM_k(t)
v\|^2_H=0$, $t\in [0,T]$, then it is natural to expect
conservation of energy. Once again, analysis of
(\ref{eq:LinHp2})--(\ref{eq:LinHp3}) shows that equality
$$
\mbE\|u(t)\|_H^2 = \|u_0\|_H^2
$$
holds if and only if
$$
 \lim_{n\to \infty}
 \int_0^T\int_0^t \int_0^{s_{n}}\ldots\int_0^{s_2}
 \|\cM_{k_{n+1}} \Phi_{t,s_{n}}\cM_{k_{n}} \ldots
 \cM_{k_1} \Phi_{s_1}u_0\|_{H}^2ds^{n}dt =0.
 $$
 The proof of Corollary \ref{cor:LinH1} shows that a  sufficient condition for the
 conservation of energy in a fully degenerate homogeneous equation is
 $\mbE\int_0^T\|u(t)\|_V^2dt < \infty$.
\end{remark}
One of applications of the Wiener Chaos solution is new numerical
methods for solving the evolution equations. Indeed, an
approximation of the solution is obtained by truncating the sum
$\sum_{\a\in \cJ}u_{\alpha}(t)\xi_{\a}$. For the Zakai filtering
equation, these numerical methods were studied in
\cite{LMR,LR1,LR2}; see also Section \ref{sec:.:FLT} below.
 The main question in the analysis is the rate of convergence,
in $n$, of the series $\sum_{n\geq 1}\sum_{|\a|=n}\|u(t)\|_H^2$.
In general, this convergence can be arbitrarily slow. For example,
consider the equation
$$
du=\frac{1}{2}u_{xx}dt+u_xdw(t),\ t>0, \; x\in \bR,
$$
in the normal triple $(H^1_2(\bR), L_2(\bR), H^{-1}_2(\bR))$, with
initial condition $u|_{t=0}=u_0\in L_2(\bR)$. It follows from
 (\ref{eq:L2norm}) that
$$
F_n(t)=\sum_{|\alpha|=n}\|u\|_{L_2(\bR)}^2(t) =
\frac{t^n}{n!}\int_{\bR}|y|^{2n}e^{-y^2t}|\hat{u}_0|^2dy,
$$
where $\hat{u}_0$ is the Fourier transform of $u_0$. If
$$
|\hat{u}_0(y)|^2=\frac{1}{(1+|y|^2)^{\gamma}}, \ \gamma>1/2,
$$
then the rate of decay of $F_n(t)$ is close to
$n^{-(1+2\gamma)/2}$. Note that, in this example,
$\mbE\|u\|^2_{L_2(\bR)}(t)=\|u_0\|_{L_2(\bR)}^2$.

An exponential  convergence rate that is uniform in $\|u_0\|_H^2$
is  achieved under strong parabolicity  condition
(\ref{eq:ParabCll}).
 An even faster  factorial rate is achieved when the operators $\cM_k$
 are bounded on $H$.

\begin{theorem}
\label{th:conv-rate} Assume that the there exist a positive number
$\varepsilon$ and a real number $C_0$ so that
$$
2\langle \cA (t)v, v\rangle + \sum_{k\geq 1}\|\cM_k(t)v\|_H^2 +
\varepsilon \|v\|_V^2 \leq C_0\|v\|_H^2, \ t\in [0,T],\ v\in V.
$$
Then there exists a positive number $b$  so that, for all $t\in
[0,T]$,
\begin{equation}
\label{eq:123} \sum_{|\a|=n} \|u_{\alpha}(t)\|_H^2 \leq
\frac{\|u_0\|_H^2}{(1+b)^n}.
\end{equation}

If, in addition, $\sum_{k\geq 1}\|\cM_k(t) \varphi\|_H^2 \leq
C_3\|\varphi\|_{H}^2$, then
\begin{equation}
\label{eq:1234} \sum_{|\a|=n} \|u_{\alpha}(t)\|_H^2 \leq
\frac{(C_3t)^n}{n!}e^{C_1t}\|u_0\|_{H}^2.
\end{equation}
\end{theorem}
\begin{proof} If $C^*$ is from (\ref{eq:contMk}) and
 $b=\varepsilon/C^*$, then the operators $\sqrt{1+b}\cM_k$ satisfy
 $$
2\langle \cA(t) v, v\rangle + (1+b)\sum_{k\geq 1}\|\cM_k(t)\|_H^2
\leq C_0\|v\|_H^2.
$$
By Theorem \ref{th:linH1},
$$
 (1+b)^n\sum_{k_1, \ldots, k_n\geq 1} \int_0^t \int_0^{s_n}\ldots\int_0^{s_2}
 \|\Phi_{t,s_n}\cM_{k_n} \ldots \cM_{k_1}\Phi_{s_1}u_0\|_H^2ds^n
 \leq \|u_0\|_H^2,
 $$
 and (\ref{eq:123}) follows.

To establish (\ref{eq:1234}), note that, by (\ref{eq:linHa1}),
$$
\|\Phi_t f\|_{H}^2 \leq e^{C_1 t}\|f\|^2_{H},
$$
and therefore the result follows from (\ref{eq:L2norm}). Theorem
\ref{th:conv-rate} is proved.
\end{proof}

The Wiener Chaos solution of (\ref{eq:linH}) is not, in general, a
solution of the equation in the  sense of Definition
\ref{def:ParabCl}. Indeed, if  $u\not\in L_2(\Omega\times(0,T);
V)$, then
 the expressions $\langle \cA u(s), \varphi\rangle$ and
$(\cM_k u(s), \varphi)_H$ are not defined. On the other hand, if
there is a possibility to move the operators $\cA$ and $\cM$ from
the solution process $u$ to the test function $\varphi$, then
equation (\ref{eq:linH}) admits a natural analog of the
traditional weak formulation (\ref{eq:trad-def}).

\begin{theorem}
\label{th:linH2} In addition to {\textbf{A1}}--{\textbf{A3}},
 assume that there exist operators $\cA^*(t)$, $\cM^*_k(t)$ and a  dense subset
$V_0$ of the space $V$ so that
\begin{enumerate}
\item $\cA^*(t)(V_0) \subseteq H$, $\cM^*_k(t)(V_0) \subseteq H$,
 $t\in [0,T]$;
\item for every $v \in V$,  $\varphi \in V_0$, and $t\in [0,T]$,
$\langle \cA(t) v, \varphi\rangle=(v, \cA^*(t) \varphi)_H,$
$(\cM_k(t) v, \varphi)_H = (v, \cM^*_k(t) \varphi)_H$.
\end{enumerate}
If $u=u(t)$ is the Wiener Chaos solution  of (\ref{eq:linH}),
then,
 for every $\varphi \in V_0$ and every $t\in [0,T]$,
 the equality
\begin{equation}
\label{eq:linHw1}
\begin{split}
(u(t), \varphi)_H&=(u_0, \varphi)_H+\int_0^t(  u(s),
\cA^*(s)\varphi)_H
ds +\int_0^t\langle f(s), \varphi\rangle ds\\
&+ \int_0^t( u(s), \cM^*_k(s)\varphi)_Hdw_k(s) +
\int_0^t(g_k(s),\varphi)_Hdw_k(s)
\end{split}
\end{equation}
holds in $L_2(\mbW)$.
\end{theorem}

\begin{proof}
The arguments are identical to the proof of Theorem
\ref{th:WC-trad}(2).
\end{proof}

As was mentioned earlier, the Wiener Chaos solution can be
constructed for anticipating equations, that is, equations with
$\cF^W_T$-measurable input data. With obvious modifications,
inequality (\ref{eq:linHn}) holds if each of the input functions
$u_0, f$, and $g_k$ in (\ref{eq:linH}) is a {\em finite } linear
combination of the basis elements $\xi_{\alpha}$. The following
example demonstrates that inequality (\ref{eq:linHn}) is
impossible for general anticipating equation.

\begin{example}
\label{ex:ant1} {\rm Let  $u=u(t,x)$ be a Wiener Chaos solution of
an ordinary differential equation
\begin{equation}
\label{eq:ant_ex1} du=udw(t), 0<t\leq 1,
\end{equation}
with $u_0=\sum_{\alpha\in \cJ} a_{\alpha}\xi_{\alpha}$. For $n\geq
0$, denote by $(n)$ the multi-index with $\alpha_1=n$ and
$\alpha_i=0$, $i\geq 2$, and assume that $a_{(n)}>0$, $n\geq 0$.
Then
\begin{equation}
\label{eq:ant_ex2} \mbE u^2(1) \geq C \sum_{n\geq
0}e^{\sqrt{n}}a_{(n)}^2.
\end{equation}
Indeed, the first column of  propagator for $\alpha=(n)$ is
$u_{(0)}(t)=a_{(0)}$ and
$$
u_{(n)}(t)=a_{(n)}+\sqrt{n}\int_0^tu_{(n-1)}(s)ds,
$$
so that
$$
u_{(n)}(t)=\sum_{k=0}^n \frac{\sqrt{n!}}{\sqrt{(n-k)!k!}}\;
\frac{a_{(n-k)}}{\sqrt{k!}} t^k.
$$
Then $u_{(n)}^2(1)\geq \sum_{k=0}^n
\binom{n}{k}\frac{a^2_{(n-k)}}{k!}$ and
$$
\sum_{n\geq 0}u^2_{(n)}(1) \geq \sum_{n\geq 0}\left( \sum_{k\geq
0} \frac{1}{k!} \binom{n+k}{n}\right) a_{(n)}^2.
$$
Since
$$
\sum_{k\geq 0} \frac{1}{k!}\binom{n+k}{n} \geq \sum_{k\geq 0}
\frac{n^k}{(k!)^2} \geq Ce^{\sqrt{n}},
$$
the result follows. }
\end{example}

The consequence of  Example \ref{ex:ant1} is that it is possible,
in (\ref{eq:linH}), to have $u_0\in L_{2}^n(\mbW; H)$ for every
$n$, and still get $\mbE\|u(t)\|_H^2=+\infty$ for all $t>0$. More
generally, the solution operator for (\ref{eq:linH}) is not
bounded on any $L_{2,Q}$ or $(\cS)_{-\rho,-\gamma}$. On the other
hand, the following result holds.
\begin{theorem}
\label{th:ant} In addition to  Assumptions {\textbf{A1}},
{\textbf{A2}}, let $u_0$ be an element of $\cD'(\mbW; H)$, $f$, an
element of $ \cD'(\mbW; L_2((0,T), V'))$, and each $g_k$, an
element of $ \cD'(\mbW; L_2((0,T), H))$.  Then the  Wiener Chaos
solution of  equation (\ref{eq:linH}) satisfies
\begin{equation}
\label{eq:antN}
\begin{split}
\sqrt{\sum_{\alpha \in \cJ} \frac{\|u_{\alpha}(t)\|_H^2}{\alpha!}}
& \leq C \sum_{\alpha\in \cJ}\frac{1}{\sqrt{\alpha !}} \Bigg(
\|u_{0\alpha}\|_H + \left(\int_0^t \|f_{\alpha}(s)\|_{V'}^2ds
\right)^{1/2} \\
& + \left( \sum_{k\geq 1}\int_0^t \|g_{k,
\alpha}(s)\|^2_Hds\right)^{1/2}\Bigg),
\end{split}
\end{equation}
where $C>0$ depends only on $T$ and the numbers $\delta, C_1$, and
$C_2$ from (\ref{eq:linHa1}) and (\ref{eq:linHa2}).
\end{theorem}

\begin{proof} To simplify the presentation, assume that $f=g_k=0$.
For fixed $\gamma \in \cJ$, denote by $u(t;\varphi; \gamma) $ the
Wiener Chaos solution of the
 equation (\ref{eq:linH})  with initial condition $u(0;\varphi; \gamma)=
\varphi \xi_{\gamma}$. Denote by $(0)$ the zero multi-index. The
structure of the propagator implies the following  relation:
\begin{equation}
\label{eq:antM} \frac{u_{\alpha+\gamma}(t;\varphi;
\gamma)}{\sqrt{(\alpha+\gamma)!}} = \frac{u_{\alpha}\left( t;
\frac{\varphi}{\sqrt{\gamma!}}; (0)\right)} {\sqrt{\alpha!}}.
\end{equation}
Clearly, $u_{\alpha}(t;\varphi; \gamma)=0$ if $|\alpha|<
|\gamma|$. If
$$
\|v(t)\|_{(\cS)_{-1,0}(H)}^2= \sum_{\alpha \in \cJ}
\frac{\|v_{\alpha}(t)\|_H^2}{\alpha!},
$$
then, by linearity and triangle inequality,
$$
\|u(t)\|_{(\cS)_{-1,0}(H)}\leq \sum_{\gamma \in \cJ}
\|u(t;u_{0\gamma};\gamma)\|_{(\cS)_{-1,0}(H)}.
$$
We also have by (\ref{eq:antM}) and Theorem \ref{th:linH1}
\begin{equation*}
\begin{split}
\|u(t;u_{0\gamma};\gamma)\|_{(\cS)_{-1,0}(H)}^2 &= \left\| u\left(
t; \frac{u_{0\gamma}}{\sqrt{\gamma!}}; (0)\right)
\right\|_{(\cS)_{-1,0}(H)}^2\\
&\leq \mbE\left\| u\left( t; \frac{u_{0\gamma}}{\sqrt{\gamma!}};
(0)\right)\right\|^2_{H} \leq
e^{C_2t}\frac{\|u_{0\gamma}\|^2_H}{\gamma!}.
\end{split}
\end{equation*}
Inequality (\ref{eq:antN}) then follows. Theorem \ref{th:ant} is
proved.
\end{proof}
\begin{remark}
Using Proposition  \ref{prop:ConvSum} and the Cauchy-Schwartz
inequality, (\ref{eq:antN}) can be re-written in a slightly weaker
form to reveal  continuity of the solution operator for equation
(\ref{eq:linH}) from $(\cS)_{-1,\gamma}$ to $(\cS)_{-1,0}$ for
every $\gamma>1\!:$
\begin{equation*}
\begin{split}
\|u(t)\|_{(\cS)_{-1,0}(H)}^2&\leq C \Bigg(
\|u_0\|^2_{(\cS)_{-1,\gamma}(H)}+
\int_0^t\|f(s)\|_{(\cS)_{-1,\gamma}(V')}^2ds \\
&+ \sum_{k\geq 1}\int_0^t
\|g_k(s)\|_{(\cS)_{-1,\gamma}(H)}^2ds\Bigg).
\end{split}
\end{equation*}
\end{remark}

\section{Probabilistic Representation of Wiener Chaos Solutions}
\label{secPR}
\setcounter{equation}{0}
\setcounter{theorem}{0}

The general discussion so far has been dealing with the abstract
evolution equation
$$
du=(\cA u + f)dt + \sum_{k\geq 1} (\cM_k u + g_k)dw_k.
$$
By further specifying the operators $\cA$ and $\cM_k$, as well as
the input data $u_0, f, $ and $g_k$, it is possible to get
additional information about the Wiener Chaos solution of the
equation.

\begin{definition}
\label{def:weight} For $r\in \bR$,
 the space $L_{2,(r)}=L_{2,(r)}(\bR^d)$ is the
collection of real-valued measurable functions so that $f \in
L_{2,(r)}$ if and only if $\int_{\bR^d}|f(x)|^2(1+|x|^2)^{r} dx <
\infty. $ The space $H^1_{2,(r)}=H^1_{2,(r)}(\bR^d)$ is the
collection of real-valued measurable functions so that $f\in
H^1_{2,(r)}$ if and only if $f$ and all the first-order
generalized derivatives $D_if$ of $f$ belong to $L_{2,(r)}$.
\end{definition}

It is known, for example, from Theorem 3.4.7 in \cite{Roz}, that
$L_{2,(r)}$ is a Hilbert space with norm
$$
\|f\|_{0,(r)}^2=\int_{\bR^d}|f(x)|^2(1+|x|^2)^{r} dx,
$$
and $H^1_{2,(r)}$ is a Hilbert space with norm
$$
\|f\|_{1,(r)}=\|f\|_{0,(r)}+\sum_{i=1}^d \|D_if\|_{0,(r)}.
$$
Denote by $H^{-1}_{2,(r)}$ the dual of $H^{1}_{2,(r)}$ with
respect to the inner product in $L_{2,(r)}$. Then
$(H^1_{2,(r)},L_{2,(r)},H^{-1}_{2,(r)})$ is a normal triple of
Hilbert spaces.

Let $\mbF=(\Omega, \cF, \{\cF_t\}_{t\geq 0}, \mbP)$ be a
stochastic basis with the usual assumptions and $w_k=w_k(t), \
k\geq 1,\ t\geq 0$, a collection of standard Wiener processes on
$\mbF$. Consider the linear  equation
\begin{equation}
\label{eq:linG}
du=(a_{ij}D_iD_ju+b_iD_iu+cu+f)dt+(\sigma_{ik}D_iu+\nu_ku+g_k)dw_k
\end{equation}
under the following {\bf assumptions:}
\begin{enumerate}
\item[\textbf{{B0}}] All coefficients, free terms, and the initial
condition are non-random. \item[\textbf{{B1}}] The functions
$a_{ij}=a_{ij}(t,x)$ and their first-order derivatives with
respect to $x$ are uniformly bounded in $(t,x)$, and the matrix
$(a_{ij})$ is uniformly positive definite, that is, there exists a
$\delta>0$ so that, for all vectors $y\in \bR^d$ and all  $(t,x)$,
$a_{ij}y_iy_j\geq \delta|y|^2$. \item[\textbf{{B2}}] The functions
$b_i=b_i(t,x)$, $c=c(t,x)$, and $\nu_k=\nu_k(t,x)$ are measurable
and  bounded in $(t,x)$. \item[\textbf{{B3}}] The functions
$\sigma_{ik}=\sigma_{ik}(t,x)$ are continuous and bounded in
$(t,x)$. \item[\textbf{{B4}}] The functions $f=f(t,x)$ and
$g_k=g_k(t,x)$ belong to $L_2((0,T); L_{2,(r)})$ for some $r\in
\bR$. \item[\textbf{{B5}}] The initial condition $u_0=u_0(x)$
belongs to $L_{2,(r)}$.
\end{enumerate}

Under Assumptions \textbf{{B2}}--\textbf{{B4}}, there exists a
sequence $Q=\{q_k, k\geq 1\}$ of positive numbers with the
following properties:
\begin{enumerate}
\item[\textbf{{P1}}] The matrix $A$ with $A_{ij}=a_{ij}-
(1/2)\sum_{k\geq 1}q_{k}\sigma_{ik}\sigma_{jk}$ satisfies
$$
A_{ij}(t,x)y_{i}y_{j}\geq0,
$$
 $x,y\in{\mathbb{R}}^{d},$ $0\leq
t\leq T$. \item[\textbf{{P2}}] There exists a number $C>0$ so that
$$
\sum_{k\geq
1}\left(\sup_{t,x}|q_k\nu_k(t,x)|^2+\int_0^T\|q_kg_k\|_{0,(r)}^p(t)dt
\right)\leq C.
$$
\end{enumerate}

For the matrix $A$ and each $t,x$, we have $A_{ij}(t,x)=\tilde{\sigma}%
_{ik}(t,x)\tilde{\sigma}_{jk}(t,x)$, where the functions
$\tilde{\sigma}_{ik} $ are bounded. This representation might not
be unique; see, for example, \cite[Theorem III.2.2]{F.fi} or
\cite[Lemma 5.2.1]{StrVar}.
 Given any such representation of $A$, consider the
following backward It\^{o} equation
\begin{equation}
\label{eq:Schar}
\begin{split}
X_{t,x,i}\left( s\right) &=x_i+\int_s^tB_i\left(\tau,X_{t,x}\left(
\tau\right) \right) d\tau + \sum_{k\geq 1}q_k\sigma_{ik}\left(
\tau,X_{t,x}\left( \tau\right) \right)
\overleftarrow{dw_{k}}\left( \tau\right)\\
&+ \int_s^t \tilde{\sigma}_{ik}\left(\tau, X_{t,x}\left(
\tau\right)\right)
\overleftarrow{d\tilde{w}}%
_{k}\left( \tau\right);\ s\in (0,t),\; t\in (0,T],\; t-{\rm
fixed},
\end{split}
\end{equation}
where $B_i=b_i-\sum_{k\geq 1}q_k^2\sigma_{ik}\nu_k$ and
$\tilde{w}_k, \ k\geq 1,$ are independent standard Wiener
processes on $\mbF$
 that are  independent of $w_k, \ k\geq 1$. This equation
might not have a strong solution, but  does have weak, or
martingale,
 solutions due to
Assumptions \textbf{{B1}}--\textbf{{B3}} and properties
\textbf{{P1}} and \textbf{{P2}} of the sequence $Q$; this weak
solution is unique in the sense of probability law \cite[Theorem
7.2.1]{StrVar}.

The following result is a variation of Theorem 4.1 in
\cite{LR_AP}.

\begin{theorem}
\label{th:ProbSol} Under assumptions \textbf{{B0}}--\textbf{{B5}}
equation (\ref{eq:linG}) has a unique \\
$w(H^{1}_{2,(r)}, H^{-1}_{2,(r)}) $ Wiener Chaos solution.  If $Q$
is a sequence with  properties \textbf{{P1}} and \textbf{{P2}},
then the solution of  (\ref{eq:linG}) belongs to
$$
L_{2,Q}\left(\mbW; L_2((0,T); H^1_{2,(r)})\right)\bigcap
L_{2,Q}\left(\mbW; \bC((0,T); L_{2,(r)})\right)
$$
and has the  following representation:
\begin{equation}
\label{eq:FK}
\begin{split}
u(t,x)=\cQ^{-1}\mbE \Bigg( \int_0^tf(s,X_{t,x}(s))\gamma(t,s,x)ds
\\+
\sum_{k\geq 1}\int_0^tq_kg_k(s,
X_{t,x}(s))\gamma(t,s,x)\overleftarrow{dw_{k}}(s)+
u_0(X_{t,x}(0))\gamma(t,0,x)\Big{|}\cF^W_t\Bigg),\ t\leq T,
\end{split}
\end{equation}
where $X_{t,x}(s)$ is a weak solution of (\ref{eq:Schar}), and
\begin{equation}
\begin{split}
\gamma(t,s,x)=\exp\Bigg( \int_s^t c(\tau,X_{t,x}(\tau))d\tau
&+\sum_{k\geq 1}\int_s^t
q_k\nu_k(\tau,X_{t,x}(\tau))\overleftarrow{dw_{k}}(\tau)\\
&-\frac{1}{2} \int_s^t\sum_{k\geq
1}q_k^2|\nu_k(\tau,X_{t,x}(\tau))|^2d\tau \Bigg).
\end{split}
\end{equation}
\end{theorem}

\begin{proof} It is enough to establish (\ref{eq:FK}) when  $t=T$.
Consider the equation
\begin{equation}
\label{eq:linGQ} dU=(a_{ij}D_iD_jU+b_iD_iU+cU+f)dt+ \sum_{k\geq
1}(\sigma_{ik}D_iU+\nu_kU+g_k)q_kdw_k
\end{equation}
with initial condition $U(0,x)=u_0(x)$. Applying Theorem
\ref{th:ParabCl}
 in the normal triple  $(H^1_{2,(r)},L_{2,(r)},H^{-1}_{2,(r)})$,
we conclude that there is a unique solution
$$
U\in L_2\left(\mbW; L_2((0,T);H^1_{2,(r)})\right)\bigcap
L_{2}\left(\mbW; \bC((0,T); L_{2,(r)})\right)
$$ of this equation.
By Proposition \ref{prop:QT}, the process $u=\cQ^{-1}U$ is the
corresponding Wiener Chaos solution of (\ref{eq:linG}). To
establish representation (\ref{eq:FK}), consider the S-transform
$U_h$ of $U$. According to Theorem \ref{th:SW1}, the function
$U_h$ is the unique $w(H^1_{2,(r)},H^{-1}_{2,(r)})$ solution of
the equation
\begin{equation}
\label{eq:strQ1} dU_h=(a_{ij}D_iD_jU_h+b_iD_iU_h+cU_h+f)dt+
\sum_{k\geq 1}(\sigma_{ik}D_iU_h+\nu_kU_h+g_k)q_kh_kdt
\end{equation}
with initial condition $U_h|_{t=0}=u_0$. We also define
\begin{equation}
\label{eq:FKchar}
\begin{split}
Y(T,x)&=\int_0^Tf(s,X_{T,x}(s))\gamma(T,s,x)ds
\\&+
\sum_{k\geq 1}\int_0^Tg_k(s,
X_{T,x}(s))\gamma(T,s)q_k\overleftarrow{dw_{k}}(s)+
u_0(X_{T,x}(0))\gamma(T,0,x).
\end{split}
\end{equation}
By direct computation,
\begin{equation*}
\begin{split}
\mathbb{E}\left( \mbE\left(\mathcal{E}(h)Y(T,x)
|\mathcal{F}_{T}^{W} \right) \right) =\mathbb{E}\left(
\mathcal{E}(h)Y(T,x) \right)
 =\mathbb{E}^{\prime }Y(T,x),
\end{split}
\end{equation*}
where $\mathbb{E}^{\prime }$ is the expectation with respect to the measure $d%
\mathbb{P}_{T}^{\prime }=\mathcal{E}(h)d\mathbb{P}_{T\text{ \ }}$and $%
\mathbb{P}_{T}$ is the  restriction of  $\mathbb{P}$ to
$\mathcal{F}^W_{T}.$

To proceed, let us first  assume that the input data $u_0$, $f$,
and $g_k$ are all smooth functions with compact support. Then,
applying the Feynmann-Kac formula to the solution of equation
(\ref{eq:strQ1}) and using  the Girsanov theorem (see, for
example, Theorems 3.5.1 and 5.7.6 in \cite{KarShr}), we conclude
that $U_h(T,x)=\mathbb{E}^{\prime }Y(T,x)$ or
$$
\displaystyle\mathbb{E}\left( \mathcal{E}(h)\mathbb{E}Y(t,x)
 |\mathcal{F}_{T}^{W}\right)
 =\mathbb{E}\left( \mathcal{E}\left( h\right) U(T,x) \right).
 $$
  By Remark \ref{rm:dense}, the last equality implies
   $\displaystyle U \left( T,\cdot \right) =%
\mathbb{E}\left(Y(T,\cdot) |\mathcal{F}_{T}^{W}\right) $ as
elements of $L_{2}\left( \Omega
;L_{2,(r)}(\mathbb{R}^{d})\right).$

To remove the additional smoothness assumption on the input data,
let $u^n_0$, $f^n$, and $g^{n}_k$ be sequences of smooth compactly
supported functions so that
\begin{equation}
\begin{split}
\lim_{n\to \infty}\Bigg( \|u_0-u_0^n\|_{L_{2,(r)}(\bR^d)}^2 +
\int_0^T \|f-f^n\|^2_{L_{2,(r)}(\bR^d)}(t)dt\\
+ \sum_{k\geq 1}\int_0^T
q_k^2\|g_k-g^{n}_k\|^2_{L_{2,(r)}(\bR^d)}(t)dt\Bigg)=0.
\end{split}
\end{equation}
Denote by $U^{n}$ and $Y^{n}$ the corresponding objects defined by
(\ref{eq:linGQ}) and (\ref{eq:FKchar}) respectively. By Theorem
\ref{th:linH1}, we have
\begin{equation}
\label{eq:FPlim1} \lim_{n\to \infty}
\mbE\|U-U^{n}\|^2_{L_{2,(r)}(\bR^d)}(T)=0.
\end{equation}
To complete the proof, it remains to show that
\begin{equation}
\label{eq:FKlim3} \lim_{n\to \infty}
\mbE\left\|\mbE\left(Y(T,\cdot)-
Y^{n}(T,\cdot)\Big{|}\cF^W_T\right)\right\|^2_{L_{2,(r)}(\bR^d)}=0.
\end{equation}
To this end,  introduce a new probability measure
$\mathbb{P}_{T}^{''}$ by
\begin{equation*}
\begin{split}
d\mathbb{P}_{T}^{''}&=\exp\Bigg( 2\sum_{k\geq 1}\int_0^T
\nu_k(s,X^Q_{T,x}(s))q_k\overleftarrow{dw_{k}}(s)\\
&- 2\int_0^T\sum_{k\geq 1}q_k^2|\nu_k(s,X^Q_{T,x}(s))|^2ds \Bigg)
d\mbP_T.
\end{split}
\end{equation*}
 By Girsanov's theorem, equation
(\ref{eq:Schar}) can be rewritten as
\begin{equation}
\label{eq:Schar1}
\begin{split}
X_{T,x,i}\left( s\right) &=x_i+\int_{s}^{T}\sum_{k\geq 1}
\sigma_{ik}\left(\tau, X_{T,x}\left(
\tau\right) \right) h_{k}\left( \tau\right)q_k d\tau\\
&+ \int_s^t(b_i+\sum_{k\geq
1}q_k^2\sigma_{ik}\nu_k)\left(\tau,X_{T,x}\left( \tau\right)
\right) d\tau \\
&+ \int_s^t \sum_{k\geq 1} q_k\sigma_{ik}\left( \tau,X_{T,x}\left(
\tau\right) \right) \overleftarrow{dw''_{k}}\left( \tau\right) +
\int_s^t \tilde{\sigma}_{ik}\left(\tau, X_{T,x}\left(
\tau\right)\right) \overleftarrow{d\tilde{w''}}_{k}\left(
\tau\right),
\end{split}
\end{equation}
where $w''_k$ and $\tilde{w''}_k$ are independent Winer processes
with respect to the measure $\mbP''_T$. Denote by $p(s,y|x)$ the
corresponding distribution density of $X_{T,x}(s)$ and write
$\ell(x)=(1+|x|^2)^r$.
 It then follows by the H\"{o}lder and Jensen inequalities
that
\begin{equation}
\label{eq:FKf}
\begin{split}
\mbE\left\|\mbE\left( \int_0^T
\gamma^2(T,s,\cdot)(f-f^n)(s,X_{T,\cdot}(s))ds\Big{|}
\cF^W_T\right) \right\|^2_{L_{2,(r)}(\bR^d)}
\\
\leq K_1\int_{\bR^d}\left(\int_0^T\mbE\left(
\gamma^2(T,s,x)(f-f^n)^2(s,X_{T,x}(s))\right) ds\right)\ell(x)dx
\\
\leq K_2\int_{\bR^d}\left(\int_0^T\mbE''
(f-f^n)^2(s,X_{T,x}(s)) ds\right)\ell(x)dx\\
=K_2\int_{\bR^d} \int_0^T
\int_{\bR^d}(f(s,y)-f^n(s,y))^2p(s,y|x)dy\; ds\  \ell(x)dx,
\end{split}
\end{equation}
where the number $K_1$ depends only on $T$, and the number $K_2$
depends only   on $T$ and $\sup_{(t,x)}|c(t,x)|+
 \sum_{k\geq 1}q_k^2\sup_{(t,x)} |\nu_k(t,x)|^2$.
 Assumptions {\textbf{B0}}--{\textbf{B2}}  imply that there exist positive
 numbers $K_3$ and $K_4$ so that
\begin{equation}
\label{eq:FPdens} p(s,y|x) \leq \frac{K_3}{(T-s)^{d/2}}\exp\left(
-K_4\frac{|x-y|^2}{T-s}\right);
\end{equation}
see, for example, \cite{Eid}. As a result,
$$
\int_{\bR^d}p(s,y|x)\ell(x)dx \leq K_5\ell(y),
$$ and
\begin{equation}
\label{eq:FKf1}
\begin{split}
\int_{\bR^d} \int_0^T \int_{\bR^d}(f(s,y)-f^n(s,y))^2p(s,y|x)dy\;
ds\  \ell(x)dx
\\
\leq K_5\int_0^T\|f-f^n\|^2_{L_{2,(r)}(\bR^d)}(s)ds \to 0, \ n\to
\infty,
\end{split}
\end{equation}
where the number $K_5$ depends only on $K_3, K_4$, $T$, and $r$.

Calculations similar to (\ref{eq:FKf})--(\ref{eq:FKf1}) show that
\begin{equation}
\begin{split}
&\mbE\left\|\mbE\left(
\gamma^2(T,0,\cdot)(u_0-u_0^n)(X_{T,\cdot}(0))\Big{|} \mbW\right)
\right\|^2_{L_{2,(r)}(\bR^d)}\\
&+ \mbE\left\|\mbE\left(\int_0^T\sum_{k\geq 1}(g_k-g^{n}_k) (s,
X_{T,\cdot}(s))\gamma(t,s,\cdot)q_k\overleftarrow{dw_{k}}(s)\Big{|}
\mbW\right) \right\|^2_{L_{2,(r)}(\bR^d)} \to 0
\end{split}
\end{equation}
as $\ n\to \infty$.
Then convergence   (\ref{eq:FKlim3}) follows, which, together with
(\ref{eq:FPlim1}),
implies that $\displaystyle U \left( T,\cdot \right) =%
\mathbb{E}\left(U^Q(T,\cdot) |\mathcal{F}_{T}^{W}\right) $ as
elements of $L_{2}\left( \Omega
;L_{2,(r)}(\mathbb{R}^{d})\right).$ It remains to note that
$u=\cQ^{-1}U$.
 Theorem \ref{th:ProbSol} is proved. \end{proof}

Given  $f \in L_{2,(r)}$, we say that $f\geq 0$ if and only if
$$
\int_{\bR^d}f(x)\varphi(x)dx \geq 0
$$
for every non-negative  $\varphi \in \bC^{\infty}_0(\bR^d)$.
 Then  Theorem \ref{th:ProbSol} implies the following result.

 \begin{corollary}
 \label{cor:pos}
 In addition to Assumptions {\textbf{B0}}--{\textbf{B5}}, let $u_0\geq 0$,
 $f\geq 0$, and $g_k=0$ for all $k\geq 1$. Then $u\geq 0$.
 \end{corollary}

 \begin{proof} This follows from (\ref{eq:FK}) and Proposition \ref{prop:pos}.
 \end{proof}

 \begin{example}

\label{ex:kvf}(Krylov-Veretennikov formula)

{\rm

 Consider the equation
\begin{equation}
du=\left(  a_{ij}D_{i}D_{j}u+b_{i}D_{i}u\right)
dt+\sum_{k=1}^{d}\sigma
_{ik}D_{i}udw_{k},\ u\left(  0,x\right)  =u_{0}\left(  x\right).\label{eq:bkr}%
\end{equation}
Assume \textbf{B0}--\textbf{B5 }and suppose that $a_{ij}(t,x)=\frac{1}{2}%
\sigma_{ik}(t,x)\sigma_{jk}(t,x)$.
 By Theorem \ref{th:linH1}, equation
(\ref{eq:bkr}) has a unique Wiener chaos solution  so  that
\[
\mathbb{E}\Vert u\Vert_{L_{2}({\mathbb{R}}^{d})}^{2}(t)\leq
C^{\ast}\Vert u_{0}\Vert_{L_{2}({\mathbb{R}}^{d})}^{2}
\]
and
\begin{equation}%
\begin{split}
u\left(  t,x\right)   &
=\sum_{n=1}^{\infty}\sum_{|\alpha|=n}u_{\alpha
}(t,x)\xi_{\alpha}=u_{0}\left(  x\right)
+\sum_{n=1}^{\infty}\sum_{k_{1},\ldots
,k_{n}=1}^d\int_{0}^{t}\int_{0}^{s_{n}}\ldots\int_{0}^{s_{2}}\\
&  \Phi_{t,s_{n}}\sigma_{jk_{n}}D_{j}\cdots \Phi_{s_{2},s_{1}}\sigma_{ik_{1}}%
D_{i}\Phi_{s_1,0}u_{0}(x)dw_{k_{1}}(s_{1})\cdots
dw_{k_{n}}(s_{n}),
\end{split}
\end{equation}
where $\Phi_{t,s}$\ \ is the semi-group generated by the operator $\mathcal{A}%
=a_{ij}D_{i}D_{j}u+b_{i}D_{i}u.$ On the other hand, in this case,
Theorem \ref{th:ProbSol} yields
\[
u(t,x)={\mathbb{E}}\Bigg(u_{0}(X_{t,x}(0))\left\vert {\mathcal{F}}_{t}%
^{W}\right.  \Bigg),
\]
where $W=\left(  w_{1},...,w_{d}\right)  $ and
\begin{equation}%
\begin{split}
X_{t,x,i}\left(  s\right)   &  =x_{i}+\int_{s}^{t}b_{i}\left(
\tau ,X_{t,x}\left(  \tau\right)  \right)
d\tau+\sum_{k=1}^{d}\sigma_{ik}\left( \tau,X_{t,x}\left(
\tau\right)  \right)  \overleftarrow{dw_{k}}\left(
\tau\right)  \\
&  \ s\in(0,t),\;t\in(0,T],\;t-\mathrm{fixed}.
\end{split}
\end{equation}
Thus, we have arrived at the Krylov-Veretennikov formula
\cite[Theorem 4]{KV}
\begin{equation}%
\begin{split}
\mathbb{E}\left(  u_{0}\left(  X_{t,x}\left(  0\right)  \right)
|\mathcal{F}_{t}^{W}\right)   &  =u_{0}\left(  x\right)  +\sum_{n=1}^{\infty}%
\sum_{k_{1},\ldots,k_{n}=1}^d\int_{0}^{t}\int_{0}^{s_{n}}\ldots\int
_{0}^{s_{2}}\\
& \!\!\!\!\!\!\!\!\!\!\!\!\!\!\!\!
 \Phi_{t,s_{n}}\sigma_{jk_{n}}D_{j}\cdots \Phi_{s_{2},s_{1}}\sigma_{ik_{1}}%
D_{i}\Phi_{s_1,0}u_{0}(x)dw_{k_{1}}(s_{1})\cdots
dw_{k_{n}}(s_{n}).
\end{split}
\end{equation}
}
\end{example}

\section{Wiener Chaos and Nonlinear Filtering}
\label{sec:.:FLT} \setcounter{equation}{0} \setcounter{theorem}{0}

In this section, we discuss some applications of the Wiener Chaos
expansion to numerical solution of the nonlinear filtering problem
for diffusion processes; the presentation is essentially based on
 \cite{LMR}.

Let $(\Omega,\cF,\mbP)$ be a complete probability space with
independent standard Wiener processes $W=W(t)$ and $V=V(t)$ of
dimensions $d_1$ and $r$ respectively. Let $X_0$ be a random
variable independent of $W$ and $V$. In the {\em diffusion
filtering model}, the unobserved $d$ - dimensional state (or
signal) process $X=X(t)$ and the $r$-dimensional observation
process $Y=Y(t)$ are defined by the  stochastic ordinary
differential equations \begin{equation}\label{fmod.sss2}
\begin{array}{l}
\ds dX(t)=b(X(t))dt+\sg(X(t))dW(t)+\rho(X(t))dV(t),\\
\ds dY(t)=h(X(t))dt+dV(t),\ 0 < t \leq T;\\
\ds X(0)=X_0,\quad Y(0)=0,
\end{array}
\end{equation} where $b(x)\in\bR^d$, $\sigma(x)\in\bR^{d\times d_1}$,
 $\rho(x)\in\bR^{d\times r}$,
$h(x)\in \bR^r$.

Denote by $\bC^n(\bR^d)$ the Banach space of bounded, $n$ times
continuously differentiable functions on $\bR^d$ with finite norm
$$
\|f\|_{\bC^n(\bR^d)}=\sup_{x\in \bR^d}|f(x)|+\max_{1\leq k\leq n}
\sup_{x\in \bR^d} |D^kf(x)|.
$$

{\bf Assumption  R1}. The the components of the functions $\sigma$
and $\rho$  are in $\bC^2(\bR^d)$,  the components of the
functions $b$ are in $\bC^1(\bR)$, the components of the function
$h$ are bounded measurable, and the random variable $X_0$ has a
density $u_0$.

{\bf Assumption R2}. The matrix $\sigma\sigma^*$ is uniformly
positive definite: there exists an $\varepsilon>0$ so that
$$
\sum_{i,j=1}^d\sum_{k=1}^{d_1}\sigma_{ik}(x)\sigma_{jk}(x)y_iy_j
\geq \varepsilon |y|^2,\ x,y\in \bR^d.
$$

Under Assumption \textbf{R1} system (\ref{fmod.sss2}) has a unique
strong solution \cite[Theorems 5.2.5 and 5.2.9]{KarShr}. Extra
smoothness of the coefficients in assumption \textbf{R1}  insure
the existence of a convenient
  representation of the optimal filter.

If $f=f(x)$ is a scalar measurable function on $\bR^d$ so that \\
 $ \sup_{0 \leq t\leq T}\mbE|f(X(t))|^2 < \infty$,
then the {\em filtering problem} for (\ref{fmod.sss2}) is to find
the best mean square estimate $\hat{f}_t$ of $f(X(t)),\ t \leq T,$
given the observations $Y(s),\ 0< s \leq t$.

 Denote by $\cF^Y_t$ the $\sg$-algebra generated  by
$Y(s),\ 0 \leq s \leq t$. Then the properties of the conditional
expectation imply that the solution of the filtering problem is
$$
\hat{f}_t=\mbE\left(f(X(t))|\cF^Y_t \right).
$$
To derive an alternative representation of $\hat{f}_t$, some
 additional constructions will be necessary.

Define a new probability measure $\widetilde{\mbP}$ on
$(\Omega,\cF)$ as follows: for $A \in \cF$,
$$
\widetilde{\mbP}(A)=\int_AZ_T^{-1}d\mbP,
$$
where
$$
Z_t = \exp \left\{ \int_0^t h^*(X(s))dY(s)-\frac{1}{2} \int_0^t
|h(X(s))|^2ds \right\}
$$
(here and below, if $\zeta \in \bR^k$, then $\zeta$ is a {\it
column} vector, $\zeta^*=(\zeta_1, \ldots,\zeta_k),$ and
$|\zeta|^2=\zeta^*\zeta$). If  the function $h$ is bounded, then
the measures $\mbP$ and $\widetilde{\mbP}$ are equivalent. The
expectation with respect to the measure $\widetilde{\mbP}$ will be
denoted by $\widetilde{\mbE}$.

The following properties of the measure $\widetilde{\mbP}$ are
well known \cite{Kal,Roz}:
\begin{enumerate}
\item[{\bf P1.}] Under the measure $\widetilde{\mbP}$, the
distributions of the Wiener process $W$ and the random variable
$X_0$ are
 unchanged, the observation process $Y$ is
a standard Wiener process, and, for $ 0 < t \leq T$,
 the state process $X$ satisfies
$$
\begin{array}{l}
\ds \!\!\!\!\!\!\!\!\!\!
 dX(t)=b(X(t))dt+\sigma(X(t))dW(t)+\rho(X(t))
 \left( dY(t)-h(X(t))dt\right),
\\
\ds \!\!\!\!\!\!\!\!\! X(0)=X_0;
\end{array}
$$

\item[{\bf P2.}] Under the measure $\widetilde{\mbP},$  the Wiener
processes $W$ and $Y$ and the random variable $X_0$ are
independent of one another;

\item[{\bf P3.}] The optimal filter $\hat{f}_t$ satisfies
\begin{equation}\label{bf} \hat{f}_t=\frac{\widetilde{\mbE} \left[
f(X(t))Z_t|\cF_t^Y \right] }{\widetilde{\mbE}[Z_t|\cF^Y_t]}.
\end{equation}
\end{enumerate}

Because of property {\bf P2} of the measure $\widetilde{\mbP}$ the
filtering problem will be studied on the probability space
$(\Omega, \cF, \widetilde{\mbP})$. In particular,  we will
consider the stochastic basis $\widetilde{\mbF}=\{\Omega,\cF,
\{\cF^Y_t\}_{0\leq t\leq T}, \widetilde{\mbP}\}$ and the Wiener
Chaos space $\widetilde{L}_2(\mbY)$ of $\cF^Y_T$-measurable random
variables $\eta$ with $\widetilde{\mbE}|\eta|^2< \infty$.

If the function $h$ is bounded, then, by the Cauchy-Schwarz
inequality, \begin{equation}\label{emb} \mbE |\eta|\leq
C(h,T)\sqrt{\widetilde{\mbE}|\eta|^2},\ \eta \in
\widetilde{L}_2(\mbY). \end{equation}

Next, consider the partial differential operators
$$
\cL g(x)=\frac{1}{2} \sum_{i,j=1}^d \left( (\sg(x)\sg^*(x))_{ij} +
(\rho(x)\rho^*(x))_{ij} \right) \frac{\pd^2 g(x)}{\pd x_i \pd
x_j}+ \sum_{i=1}^d b_i(x) \frac{\pd g(x)}{\pd x_i};
$$
$$
\cM_l g(x) = h_l(x)g(x)+\sum_{i=1}^d \rho_{il}(x) \frac{\pd
g(x)}{\pd x_i},\ l=1, \ldots, r;
$$
and their adjoints
\begin{equation*}
\begin{split}
\cL^* g(x)&=\frac{1}{2} \sum_{i,j=1}^d \frac{\pd^2}{\pd x_i \pd
x_j} \left( (\sg(x)\sg^*(x))_{ij}g(x)
+ (\rho(x)\rho^*(x))_{ij} g(x)\right) \\
&- \sum_{i=1}^d \frac{\pd}{\pd
x_i}\left( b_i(x)g(x) \right);
\end{split}
\end{equation*}
$$
\cM_l^* g(x) = h_l(x)g(x)-\sum_{i=1}^d \frac{\pd}{\pd x_i}\left(
\rho_{il}(x)g(x) \right),\  l=1, \ldots, r.
$$

 Note that, under the assumptions {\textbf{R1}} and {\textbf{R2}},
 the operators $\cL, \cL^* $ are bounded from $H^1_2(\bR^d)$
 to $H^{-1}_2(\bR^d)$, operators $\cM, \cM^*$ are bounded from $H^1_2(\bR^d)$
 to $L_2(\bR^d)$, and
\begin{equation}\label{eq:dissip} 2\langle\cL^*v,v\rangle+\sum_{l=1}^r
\|\cM_l^*v\|_{L_2(\bR^d)}^2 +\varepsilon \|v\|_{H_1^2(\bR^d)}^2
\leq C\|v\|_{L_2(\bR^d)}^2,\  v\in H^1_2(\bR^d), \end{equation}
where $\langle \cdot, \cdot \rangle$ is the duality between
$H^1_2{\bR^d}$ and $H^{-1}_2(\bR^d)$. The following result is well
known \cite[Theorem 6.2.1]{Roz}.

\begin{proposition}
\label{prop:egun} In addition to Assumptions \textbf{R1} and
\textbf{R1}
 suppose that the initial density $u_0$
belongs to  $L_2(\bR^d)$. Then there exists  a random field
$u=u(t,x),\ t \in [0,T],\ x \in \bR^d,$ with the following
properties:

1. $u\in \widetilde{L}_2(\mbY;L_2( (0,T); H^1_2(\bR^d)) )\cap
         \widetilde{L}_2(\mbY; \bC([0,T], L_2(\bR^d))). $

2. The function   $u(t,x)$ is a traditional  solution of  the
stochastic partial differential equation
\begin{equation}\label{ze}
\begin{array}{ll}
\ds du(t,x)&\ds =\cL^* u(t,x)dt+\sum_{l=1}^r\cM_l^*u(t,x)
dY_l(t),\ 0<t\leq T,\
 x \in \bR^d;\\
\ds u(0,x)&\ds =u_0(x).
\end{array}
\end{equation}

3. The equality \begin{equation}\label{uof} \widetilde{\mbE}
\left[ f(X(t))Z_t|\cF_t^Y \right] = \int_{\bR^d} f(x) u(t,x) dx
\end{equation} holds for all bounded measurable functions $f$.
\end{proposition}

The random field $u=u(t,x)$ is called the {\em unnormalized
filtering density} (UFD) and the random variable
$\phi_t[f]=\widetilde{\mbE} \left[ f(X(t))Z_t|\cF_t^Y \right]$,
the {\em unnormalized optimal filter}.

A number of authors studied the nonlinear filtering problem using
the multiple It\^{o} integral version of the Wiener chaos
\cite[etc.]{BdKl,Kunita81,Ocn,Wong}. In what follows, we construct
approximations of $u$ and $\phi_t[f]$ using the Cameron-Martin
version.

By Theorem \ref{th:WC-trad},
\begin{equation}
\label{eq:WC-flt} u(t,x)=\sum_{\a \in \cJ} u_{\a}(t,x)\xi_{\a},
\end{equation}
where
\begin{equation}
\label{eq:CMB-flt}
\xi_{\a}=\frac{1}{\sqrt{\a!}}\prod_{i,k}H_{\a^k_i}(\xi_{ik}),\
\xi_{ik}=\int_0^Tm_i(t)dY_k(t),\ k=1,\ldots, r;
\end{equation}
as before, $H_{n}(\cdot)$ is the Hermite polynomial
(\ref{eq:HerPol}) and $m_i\in \mathfrak{m}$, an orthonormal basis
in $L_2((0,T))$. The functions  $u_{\a}$ satisfy the corresponding
propagator
\begin{equation}
\label{eq:prop-flt}
\begin{split}
\frac{\partial}{\partial t}u_{\a}(t,x)& =\cL^* u_{\a}(t,x)\\
&+\sum_{k,i}\sqrt{\a^k_i}\cM_k^*u_{\a^{-}(i,k)}(t,x)m_i(t),\
0<t\leq T,\
 x \in \bR^d;\\
 u(0,x)& =u_0(x)I(|\a|=0).
\end{split}
\end{equation}

Writing
$$
f_{\a}(t)=\int_{\bR^d}f(x)u_{\a}(t,x)dx,
$$
we also get a Wiener chaos expansion for the unnormalized optimal
filter:
\begin{equation}
\label{eq:WCE-f(x)} \phi_t[f]=\sum_{\a\in \cJ} f_{\a}(t)\xi_{\a},\
t\in [0,T].
\end{equation}

For a positive integer $N$, define
\begin{equation}
\label{eq:flt-apporx1} u_N(t,x)=\sum_{|\a|\leq
N}u_{\a}(t,x)\xi_{\a}.
\end{equation}

\begin{theorem}
\label{th:flt-appr1} Under Assumptions \textbf{R1} and
\textbf{R2}, there exists a positive number $\nu$, depending only
on the functions $h$ and $\rho$, so that
\begin{equation}
\label{eq:flt-err1}
\widetilde{\mbE}\|u-u_N\|_{L_2(\bR^d)}^2(t)\leq
 \frac{\|u_0\|_{L_2(\bR^d)}^2}{\nu(1+\nu)^N}, \ t\in [0,T].
 \end{equation}
 If, in addition, $\rho=0$, then there exists a real number $C$, depending
 only on the functions $b$ and $\sigma$, so that
\begin{equation}
\label{eq:flt-err10}
\widetilde{\mbE}\|u-u_N\|_{L_2(\bR^d)}^2(t)\leq
 \frac{(4h_{\infty}t)^{N+1}}{(N+1)!}e^{Ct}\|u_0\|_{L_2(\bR^d)}^2, \ t\in [0,T],
 \end{equation}
where $h_{\infty}=\max_{k=1,\ldots,r}\sup_x|h_k(x)|$.
\end{theorem}

For positive integers $N,n$, define a set of multi-indices
$$
\cJ^n_N=\{\a=(\a^k_i,\ k=1,\ldots, r,\ i=1,\ldots,n): \ |\a|\leq
N\}.
$$
and let
\begin{equation}
\label{eq:flt-approx2} u_N^n(t,x)=\sum_{\a \in
\cJ^n_N}u_{\a}(t,x)\xi_{\a}.
\end{equation}

Unlike Theorem \ref{th:flt-appr1}, to compute the approximation
error in this case
 we need to choose a  special basis $\mathfrak{m}$
--- to do the error analysis for the Fourier approximation in time. We also
need extra regularity of the coefficients in the state and
observation  equations
--- to have the semi-group generated by the operator $\cL^*$ continuous
not only in $L_2(\bR^d)$ but also in $H^2_2(\bR^d)$. The resulting
error bound is presented below; the proof can be found in
\cite{LMR}.

\begin{theorem}
\label{th:flt-appr2} Assume that
\begin{enumerate}
\item The  basis  $\mathfrak{m}$ is the Fourier cosine basis
\begin{equation}\label{basis} \!\!\!\!\!\!\!\!\!\!
 m_1(s)\!=\!\frac 1{\sqrt{T}};\  m_k(t)\!=
 \!\sqrt{\frac{2}{T} }
 \cos \left( \frac{\pi (k-1) t}{T} \right), \, k>1; \
0\leq t \leq T, \end{equation} \item The  components of the
functions $\sigma$ are in $\bC^4(\bR^d)$,  the components of the
functions $b$ are in $\bC^3(\bR)$, the components of the function
$h$ are in $\bC^2(\bR^d)$; $\rho=0$; $u_0\in H^2_2(\bR^d)$.
\end{enumerate}
Then there exist a positive number $B_1$ and a real number $B_2$,
both depending only on the functions $b$ and $\sigma$ so that
\begin{equation}
\label{eq:flt-err2}
\widetilde{\mbE}\|u-u_N^n\|_{L_2(\bR^d)}^2(T)\!\leq\! B_1e^{B_2T}
\left(
\frac{(4h_{\infty}T)^{N+1}}{(N+1)!}e^{Ct}\|u_0\|_{L_2(\bR^d)}^2+
\frac{T^3}{n}\|u_0\|_{H_2^2(\bR^d)}^2\right)\!,
 \end{equation}
 where $h_{\infty}=\max_{k=1,\ldots,r}\sup_x|h_k(x)|$.
\end{theorem}

\section{Passive Scalar in a Gaussian Field}
\label{secPS}
\setcounter{equation}{0}
\setcounter{theorem}{0}

This section presents the results from \cite{LR_AP} and
\cite{LR_ps} about the stochastic transport equation.

The following viscous transport equation is used to describe time
evolution of a scalar quantity $\theta$ in a given velocity field
${\mathbf{v}}$:
\begin{equation}
\label{eq:ps}\dot{\theta}(t,x)=\nu\Delta\theta(t,x) -
{\mathbf{v}}(t,x) \cdot\nabla\theta(t,x) +f(t,x);\
x\in{\mathbb{R}}^{d},\ d>1.
\end{equation}
The scalar $\theta$ is called passive because it does not affect
the velocity field $\mathbf{v}$.

We assume that ${\mathbf{v}}={\mathbf{v}}(t,x)\in{\mathbb{R}}^{d}$
is an isotropic Gaussian vector field with zero mean and
covariance
\[
{\mathbb{E}}(v^{i}(t,x)v^{j}(s,y))=\delta(t-s)C^{ij}(x-y),
\]
where $C=(C^{ij}(x), i,j=1, \ldots, d)$ is a matrix-valued
function so that $C(0)$ is a scalar matrix; with no loss of
generality we will assume that $C(0)=I,$ the identity matrix.

It is known from \cite[Section 10.1]{LeJan} that, for an isotropic
Gaussian vector field, the Fourier transform $\hat{C}=\hat{C}(z)$
of the function $C=C(x)$ is
\begin{equation}
\label{eq:FTC}\hat{C}(y)=\frac{A_{0}}{(1+|y|^{2})^{(d+\alpha)/2}}\left(
a\frac{yy^{*}}{|y|^{2}}+\frac{b}{d-1} \left(
I-\frac{yy^{T}}{|y|^{2}}\right) \right)  ,
\end{equation}
where $y^{*}$ is the row vector $(y_{1}, \ldots, y_{d})$, $y$ is
the corresponding column vector, $|y|^{2}=y^{*}y$;
$\gamma>0,\,a\geq 0,\,b\geq0,\,A_{0}>0$ are real numbers. Similar
to \cite{LeJan}, we assume that $0<\gamma<2$. This range of values
of $\gamma$ corresponds to a turbulent velocity field
${\mathbf{v}}$, also known as the generalized Kraichnan model
\cite{GV}; the original Kraichnan model \cite{Krch} corresponds to
$a=0$. For small $x$, the asymptotics of $C^{ij}(x)$ is
$(\delta_{ij}-c^{ij}|x|^{\gamma })$ \cite[Section 10.2]{LeJan}.

By direct computation (cf. \cite{BH}), the vector field
${\mathbf{v}}=(v^{1}, \ldots, v^{d})$ can be written as
\begin{equation}
\label{eq:v}v^{i}(t,x)= \sigma^{i}_{k}(x)\dot{w}_{k}(t),
\end{equation}
where $\{\sigma_{k}, \ k\geq1\}$ is an orthonormal basis in the
space $H_{C}$, the reproducing kernel Hilbert space corresponding
to the kernel function $C$. It is known from \cite{LeJan} that
$H_{C}$ is all or part of the Sobolev space
$H^{(d+\gamma)/2}({\mathbb{R}}^{d}; {\mathbb{R}}^{d})$.

If $a>0$ and $b>0$, then the matrix $\hat{C}$ is invertible and
\[
H_{C}=\left\{  f\in{\mathbb{R}}^{d}: \int_{{\mathbb{R}}^{d}} \hat{f}%
^{*}(y)\hat{C}^{-1}(y)\hat{f}(y)dy < \infty\right\}  =
H^{(d+\gamma )/2}({\mathbb{R}}^{d}; {\mathbb{R}}^{d}),
\]
because $\|\hat{C}(y)\| \sim(1+|y|^{2})^{-(d+\gamma)/2}$.

If $a>0$ and $b=0$, then
\[
H_{C}=\left\{  f\in{\mathbb{R}}^{d}: \int_{{\mathbb{R}}^{d}} |\hat{f}%
(y)|^{2}(1+|y|^{2})^{(d+\gamma)/2}dy < \infty; \ yy^{*}\hat{f}(y)=|y|^{2}%
\hat{f}(y) \right\}  ,
\]
the subset of gradient fields in
$H^{(d+\gamma)/2}({\mathbb{R}}^{d}; {\mathbb{R}}^{d})$, that is,
vector fields $f$ for which $\hat{f}(y)=y\hat {F}(y)$ for some
scalar $F \in H^{(d+\gamma+2)/2}({\mathbb{R}}^{d})$.

If $a=0$ and $b>0$, then
\[
H_{C}=\left\{  f\in{\mathbb{R}}^{d}: \int_{{\mathbb{R}}^{d}} |\hat{f}%
(y)|^{2}(1+|y|^{2})^{(d+\gamma)/2}dy < \infty; \ y^{*}\hat{f}(y)=0
\right\} ,
\]
the subset of divergence-free fields in
$H^{(d+\gamma)/2}({\mathbb{R}}^{d}; {\mathbb{R}}^{d})$.

By the embedding theorems, each $\sigma_{k}^{i}$ is a bounded
continuous function on ${\mathbb{R}}^{d}$; in fact, every
$\sigma_{k}^{i}$ is H\"{o}lder continuous of order $\gamma/2$. In
addition, being an element of the corresponding space $H_{C}$,
each $\sigma_{k}$ is a gradient field if $b=0$ and is divergence
free if $a=0$.

Equation (\ref{eq:ps}) becomes
\begin{equation}
\label{eq:ps2}d{\theta}(t,x)=(\nu\Delta\theta(t,x) +f(t,x))dt- \sum_{k}%
\sigma_{k}(x)\cdot\nabla\theta(t,x) dw_{k}(t).
\end{equation}

We summarize the above constructions in the following
\textbf{assumptions}:

\begin{enumerate}
\item[\textbf{{S1}}] There is a fixed stochastic basis
$\mathbb{F}=(\Omega, \mathcal{F}, \{\mathcal{F}_{t}\}_{t\geq0},
\mathbb{P})$ with the usual assumptions and
$(w_{k}(t),k\geq1,t\geq0)$ is a collection of independent standard
Wiener processes on $\mathbb{F}$.

\item[\textbf{{S2}}] For each $k$, the vector field $\sigma_{k}$
is an element of the Sobolev space
\newline$H^{(d+\gamma)/2}_{2}(\mathbb{R}^{d}; \mathbb{R}^{d})$,
$0<\gamma<2$, $d\geq2$.

\item[\textbf{{S3}}] For all $x,y$ in $\mathbb{R}^{d}$, $\sum_{k}
\sigma _{k}^{i}(x)\sigma_{k}^{j}(y)=C^{ij}(x-y)$ so that the
matrix-valued function $C=C(x)$ satisfies (\ref{eq:FTC}) and
$C(0)=I$.

\item[\textbf{{S4}}] The input data $\theta_{0}, f$ are
deterministic and satisfy
\[
\theta_{0}\in L_{2}({\mathbb{R}}^{d}),\ f \in L_{2}((0,T); H_{2}%
^{-1}({\mathbb{R}}^{d}));
\]
$\nu> 0$ is a real number.
\end{enumerate}

\begin{theorem}
\label{th:ps} Let $Q$ be a sequence with $q_{k}=q< \sqrt{2\nu}$,
$k\geq1$.

Under assumptions \textbf{S1}--\textbf{S4}, there exits a unique $w(H^{1}%
_{2}({\mathbb{R}}^{d}),H^{-1}_{2}({\mathbb{R}}^{d}))$ Wiener Chaos
solution of (\ref{eq:ps2}). This solution is an
${\mathcal{F}}_{t}^{W}$-adapted process and satisfies
\[%
\begin{split}
\|\theta\|_{L_{2,Q}({\mathbb{W}}; L_{2}((0,T);H^{1}_{2}({\mathbb{R}}^{d}%
)))}^{2}  &  + \|\theta\|_{L_{2,Q}({\mathbb{W}};
{\mathbf{C}}((0,T);
L_{2}({\mathbb{R}}^{d})))}^{2}\\
&  \leq C(\nu,q,T)\left(
\|\theta_{0}\|_{L_{2}({\mathbb{R}}^{d})}^{2} + \|f\|_{L_{2}((0,T);
H_{2}^{-1}({\mathbb{R}}^{d}))}^{2}\right)  .
\end{split}
\]

\end{theorem}

Theorem \ref{th:ps} provides new information about the solution of
equation (\ref{eq:ps}) for all values of $\nu>0$. Indeed, if
$\sqrt{2\nu}>1$, then $q>1$ is an admissible choice of the
weights, and, by Proposition \ref{prop:QT}(1), the solution
$\theta$ has Malliavin derivatives of every order. If
$\sqrt{2\nu}\leq1$, then equation (\ref{eq:ps2}) does not have a
square-integrable solution.

Note that if the weight is chosen so that $q=\sqrt{2\nu}$, then
equation (\ref{eq:ps}) can still be analyzed using Theorem
\ref{th:linH1} in the normal
triple $(H^{1}_{2}({\mathbb{R}}^{d}),L_{2}({\mathbb{R}}^{d}),H^{-1}%
_{2}({\mathbb{R}}^{d}))$.

If $\nu=0$, equation (\ref{eq:ps2}) must be interpreted in the
sense of Stratonovich:
\begin{equation}
\label{eq:psStr}du(t,x)=f(t,x)dt-\sigma_{k}(x)\cdot\nabla\theta(t,x)\circ
dw_{k}(t).
\end{equation}
To simplify the presentation, we assume that $f=0$. If
(\ref{eq:FTC}) holds with $a=0$, then each $\sigma_{k}$ is
divergence free and (\ref{eq:psStr}) has an equivalent It\^{o}
form
\begin{equation}
\label{eq:psDF}d{\theta}(t,x)= \frac{1}{2}\Delta\theta(t,x)dt-\sigma_{k}%
^{i}(x)D_{i}\theta(t,x)dw_{k}(t).
\end{equation}
Equation (\ref{eq:psDF}) is a model of non-viscous turbulent
transport \cite{EVE}. The propagator for (\ref{eq:psDF}) is
\begin{equation}
\frac{\partial}{\partial t}\theta_{\alpha}(t,x)= \frac{1}{2}
\Delta
\theta_{\alpha}(t,x)-\sum_{i,k}\sqrt{\alpha^{k}_{i}} \sigma_{k}^{j}D_{j}%
\theta_{\alpha^{-}(i,k)}(t,x)m_{i}(t),\ 0<t\leq T,
\end{equation}
with initial condition
$\theta_{\alpha}(0,x)=\theta_{0}(x)I(|\alpha|=0). $

The following result about solvability of (\ref{eq:psDF}) is
proved in \cite{LR_AP} and, in a slightly weaker form, in
\cite{LR_ps}.

\begin{theorem}
\label{th:psDF} In addition to \textbf{S1}--\textbf{S4}, assume
that each
$\sigma_{k}$ is divergence free. Then there exits a unique $w(H^{1}%
_{2}({\mathbb{R}}^{d}),H^{-1}_{2}({\mathbb{R}}^{d}))$ Wiener Chaos
solution $\theta=\theta(t,x)$ of (\ref{eq:psDF}). This solution
has the following properties:

(A) For every
$\varphi\in{\mathbf{C}}_{0}^{\infty}(\mathbb{R}^{d})$ and all
$t\in[0,T]$, the equality
\begin{equation}
\label{eq:psDF1}(\theta,\varphi)(t)=(\theta_{0},\varphi)+\frac{1}{2}\int
_{0}^{t}(\theta,\Delta\varphi)(s)ds+\int_{0}^{t}(\theta,\sigma_{k}^{i}%
D_{i}\varphi)dw_{k}(s)
\end{equation}
holds in $L_{2}({\mathcal{F}}^{W}_{t})$, where $(\cdot, \cdot)$ is
the inner product in $L_{2}({\mathbb{R}}^{d})$.

(B) If $X=X_{t,x}$ is a weak solution of
\begin{equation}
\label{eq:psDFchar}X_{t,x} =x+\int_{0}^{t}\sigma_{k}\left(
X_{s,x}\right) {dw_{k}}\left(  s\right)  ,
\end{equation}
then, for each $t\in[0,T]$,
\begin{equation}
\label{eq:psDF2}\theta\left(  t,x\right)  =\mathbb{E}\left(
\theta_{0}\left( X_{t,x} \right)  |\mathcal{F}_{t}^{W}\right)  .
\end{equation}

(C) For $1\leq p <\infty$ and $r \in{\mathbb{R}}$, define $L_{p,(r)}%
({\mathbb{R}}^{d})$ as the Banach space of measurable functions
with norm
\[
\|f\|_{L_{p,(r)}({\mathbb{R}}^{d})}^{p} = \int_{{\mathbb{R}}^{d}}
|f(x)|^{p}(1+|x|^{2})^{pr/2}dx
\]
is finite. Then there exits a number $K$ depending only on $p,r$
so that, for each $t>0$,
\begin{equation}
\label{eq:psDF3}{\mathbb{E}}
\|\theta\|^{p}_{L_{p,(r)}({\mathbb{R}}^{d})}(t) \leq
e^{Kt}\|\theta_{0}\|^{p}_{L_{p,(r)}({\mathbb{R}}^{d})}.
\end{equation}
In particular, if $r=0$, then $K=0$.
\end{theorem}

It follows that, for all $s,t$ and almost all $x,y,$
\begin{align*}
\mathbb{E}\theta\left(  t,x\right)   &  =\theta_{\alpha}\left(
t,x\right)
I_{\left\vert \alpha\right\vert =0}\\
&  \text{\textrm{and}}\\
\mathbb{E}\theta\left(  t,x\right)  \mathbb{\theta}\left(
s,y\right)   &
=\sum_{\alpha\in\mathcal{J}}\mathbb{\theta}_{\alpha}\left(
t,x\right) \mathbb{\theta}_{\alpha}\left(  s,y\right)  .
\end{align*}

If the initial condition $\theta_{0}$ belongs to
$L_{2}({\mathbb{R}}^{d})\cap L_{p}({\mathbb{R}}^{d})$ for
$p\geq3$, then, by (\ref{eq:psDF3}), higher order moments of
$\theta$ exist. To obtain the expressions of the higher-order
moments in terms of the coefficients $\theta_{\alpha}$, we need
some auxiliary constructions.

For $\alpha,\ \beta\in{\mathcal{J}}$, define $\alpha+\beta$ as the
multi-index with components $\alpha^{k}_{i}+\beta_{i}^{k}$.
Similarly, we define the multi-indices $|\alpha-\beta|$ and
$\alpha\wedge\beta=\min(\alpha,\beta)$. We write $\beta\leq\alpha$
if and only if $\beta^{k}_{i}\leq\alpha^{k}_{i}$ for all
$i,k\geq1$. If $\beta\leq\alpha$, we define
\[
\binom{\alpha}{\beta}:= \prod_{i,k}\frac{\alpha^{k}_{i}!}{\beta^{k}%
_{i}!(\alpha^{k}_{i}-\beta^{k}_{i})!}.
\]

\begin{definition}
\label{def:comp} We say that a triple of multi-indices $\left(
\alpha ,\beta,\gamma\right)  $ is complete and write $\left(
\alpha,\beta ,\gamma\right)  \in\triangle$ if all the entries of
the multi-index $\alpha+\beta+\gamma$ are even numbers and
$\left\vert \alpha-\beta\right\vert \leq\gamma\leq\alpha+\beta.$
For fixed $\alpha, \beta\in\mathcal{J},$ we write
\[
\triangle\left(  \alpha\right)  :=\left\{
\gamma,\mu\in\mathcal{J}:\left( \alpha,\gamma, \mu\right)
\in{\triangle}\right\}
\]
and
\[
\triangle(\alpha,\beta):=\{\gamma\in{\mathcal{J}}:(\alpha,
\beta,\gamma )\in\triangle\}.
\]

\end{definition}

For $\left(  \alpha,\beta,\gamma\right)  \in\triangle,$ we define
\begin{equation}
\label{eq:Psi}\Psi\left(  \alpha,\beta,\gamma\right)
:=\sqrt{\alpha !\beta!\gamma!}\left(  \left(
\frac{\alpha-\beta+\gamma}{2}\right) {\text{{\Large {!}}}}\left(
\frac{\beta-\alpha+\gamma}{2}\right) {\text{{\Large {!}}}} \left(
\frac{\alpha+\beta-\gamma}{2}\right) {\text{{\Large {!}}}}\right)
^{-1}.
\end{equation}
Note that the triple $(\alpha, \beta, \gamma)$ is complete if and
only if any permutation of the triple $(\alpha, \beta, \gamma)$ is
complete. Similarly, the value of $\Psi\left(
\alpha,\beta,\gamma\right)  $ is invariant under permutation of
the arguments.

We also define
\begin{equation}
\label{eq:C}C\left(  \gamma,\beta,\mu\right)  : =\left[
\binom{\gamma
+\beta-2\mu}{\gamma-\mu}\binom{\gamma}{\mu}\binom{\beta}{\mu}\right]
^{1/2},\ \mu\leq\gamma\wedge\beta.
\end{equation}
It is readily checked that if $f$ is a function on $\mathcal{J},$
then for
$\gamma,\beta\in\mathcal{J},$%
\begin{equation}
\label{kap}\sum_{\mu\leq\gamma\wedge\beta}C\left(
\gamma,\beta,p\right) f\left(  \gamma+\beta-2\mu\right)
=\sum_{\mu\in\left(  \gamma,\beta\right)  } f\left(  \mu\right)
\Phi\left(  \gamma,\beta,\mu\right)
\end{equation}

The next theorem presents the formulas for the third and fourth
moments of the solution of equation (\ref{eq:psDF}) in terms of
the coefficients $\theta_{\alpha}$.

\begin{theorem}
\label{th:h-moments} In addition to \textbf{S1}--\textbf{S4},
assume that each $\sigma_{k}$ is divergence free and the initial
condition $\theta_{0}$ belongs to $L_{2}({\mathbb{R}}^{d})\cap
L_{4}({\mathbb{R}}^{d})$. Then
\begin{equation}
\label{eq:third}{\mathbb{E}}\theta(t,x)\theta\left(
t^{\prime},x^{\prime }\right)  \theta(s,y) =\sum_{\left(
\alpha,\beta,\gamma\right)  \in\triangle }\Psi\left(
\alpha,\beta,\gamma\right)  \theta_{\alpha}(t,x)\theta_{\beta
}(t^{\prime},x^{\prime}) \theta_{\gamma}\left(  s,y\right)
\end{equation}
and
\begin{align}
\label{eq:fourth} &  {\mathbb{E}}
\theta(t,x)\theta(t^{\prime},x^{\prime
})\theta\left(  s,y\right)  \theta\left(  s^{\prime},y^{\prime}\right) \\
&  =\sum_{\rho\in\triangle\left(  \alpha,\beta\right)
\cap\triangle\left( \gamma,\kappa\right)  } \Psi\left(
\alpha,\beta,\rho\right)  \Psi\left( \rho,\gamma,\kappa\right)
\theta_{\alpha}\left(  t,x\right)  \theta_{\beta
}(t^{\prime},x^{\prime}) \theta_{\gamma}\left(  s,y\right)
\theta_{\kappa }\left(  s^{\prime},y^{\prime}\right)  .\nonumber
\end{align}

\end{theorem}

\begin{proof} It is known \cite{meyer} that
\begin{equation}
\label{zproduct}\xi_{\gamma}\xi_{\beta}=\sum_{\mu\leq\gamma\wedge\beta
}C\left(  \gamma, \beta,\mu\right)  \xi_{\gamma+\beta-2\mu}.
\end{equation}

Let us consider the triple product
$\xi_{\alpha}\xi_{\beta}\xi_{\gamma}.$ By
(\ref{zproduct}),%
\begin{equation}
{\mathbb{E}}\xi_{\alpha}\xi_{\beta}\xi_{\gamma}={\mathbb{E}}\sum_{\mu
\in\triangle(\alpha,\beta)} \xi_{\gamma}\xi_{\mu}\Psi\left(
\alpha,\beta ,\mu\right)  =
\begin{cases}
\Psi\left(  \alpha,\beta,\gamma\right)  , & (\alpha,\beta,\gamma)
\in
\triangle;\\
0, & \mathrm{otherwise}.
\end{cases}
\label{eq;3z}%
\end{equation}
Equality (\ref{eq:third}) now follows.

To compute the fourth moment, note that%
\begin{equation}
\label{zzz}%
\begin{split}
\xi_{\alpha}\xi_{\beta}\xi_{\gamma}  &
=\sum_{\mu\leq\alpha\wedge\beta
}C\left(  \alpha,\beta,\mu\right)  \xi_{\alpha+\beta-2\mu}\xi_{\gamma}\\
&  =\sum_{\mu\leq\alpha\wedge\beta}C\left(
\alpha,\beta,\mu\right) \sum_{\rho\leq\left(
\alpha+\beta-2\mu\right)  \wedge\gamma} C\left(
\alpha+\beta-2\mu,\gamma,\rho\right)
\xi_{\alpha+\beta+\gamma-2\mu-2\rho}.
\end{split}
\end{equation}
Repeated applications of $\left(  \ref{kap}\right)  $ yield
\[%
\begin{split}
\xi_{\alpha}\xi_{\beta}\xi_{\gamma}  &
=\sum_{\mu\leq\alpha\wedge\beta }C\left(  \alpha,\beta,\mu\right)
\sum_{\rho\in\triangle\left(  \alpha +\beta-2\mu,\gamma\right)
}\xi_{\rho}\Psi\left(  \alpha+\beta-2\mu
,\gamma,\rho\right) \\
&  =\sum_{\mu\in\triangle\left(  \alpha,\beta\right)
}\sum_{\rho\in \triangle\left(  \mu,\gamma\right)  } \Psi\left(
\alpha,\beta,\mu\right)
\Psi\left(  \mu,\gamma,\rho\right)  \xi_{\rho}%
\end{split}
\]
Thus,
\begin{align*}
{\mathbb{E}}\xi_{\alpha}\xi_{\beta}\xi_{\gamma}\xi_{\kappa}  &  =
\sum_{\mu \in\triangle\left(  \alpha,\beta\right)  }
\sum_{\rho\in\triangle\left( \mu,\gamma\right)  } \Psi\left(
\alpha,\beta,\mu\right)  \Psi\left(
\mu,\gamma,\rho\right)  I_{\left\{  \mu=\kappa\right\}  }\\
&  =\sum_{\rho\in\triangle\left(  \alpha,\beta\right)
\cap\triangle\left( \gamma,\kappa\right)  } \Psi\left(
\alpha,\beta,\rho\right)  \Psi\left( \rho,\gamma,\kappa\right)  .
\end{align*}
Equality (\ref{eq:fourth}) now follows. \end{proof}

In the same way, one can get formulas for fifth- and higher-order
moments.

\begin{remark}
\label{rm:moments} Expressions (\ref{eq:third}) and
(\ref{eq:fourth}) do not depend on the structure of equation
(\ref{eq:psDF}) and can be used to compute the third and fourth
moments of any random field with a known Cameron-Martin expansion.
The interested reader should keep in mind that the formulas for
the moments of orders higher then two should be interpreted with
care. In fact, they represent the pseudo-moments (for detail see
\cite{MR2}).
\end{remark}

We now return to the analysis of the passive scalar equation
(\ref{eq:ps2}). By reducing the smoothness assumptions on
$\sigma_{k}$, it is possible to consider velocity fields
${\mathbf{v}}$ that are more turbulent than in the Kraichnan
model, for example,
\begin{equation}
\label{eq:vWN}v^{i}(t,x)=\sum_{k\geq0}
\sigma^{i}_{k}(x)\dot{w}_{k}(t),
\end{equation}
where $\{\sigma_{k}, \ k\geq1\}$ is an orthonormal basis in $L_{2}%
({\mathbb{R}}^{d};{\mathbb{R}}^{d})$. With ${\mathbf{v}}$ as in (\ref{eq:vWN}%
), the passive scalar equation (\ref{eq:ps2}) becomes
\begin{equation}
\label{eq:ps2WN}\dot{\theta}(t,x)=\nu\Delta\theta(t,x) +f(t,x)-
\nabla \theta(t,x) \cdot\dot{W}(t,x),
\end{equation}
where $\dot{W}=\dot{W}(t,x)$ is a $d$-dimensional space-time white
noise and the It\^{o} stochastic differential is used. Previously,
such equations have been studied using white noise approach in the
space of Hida distributions \cite{Pot_Tr,Pot}. A summary of the
related results can be found in \cite[Section 4.3]{HOUZ}.

The $Q$-weighted Wiener chaos spaces allow us to state a result
that is fully analogous to Theorem \ref{th:ps}. The proof is
derived from Theorem \ref{th:linH1}; see \cite{LR_AP} for details.

\begin{theorem}
\label{th:psWN} Suppose that $\nu> 0$ is a real number, each $|\sigma_{k}%
^{i}(x)| $ is a bounded measurable function, and the input data
are deterministic and satisfy $u_{0}\in L_{2}(\mathbb{R}^{d})$,
$f\in L_{2}\left( (0,T); H^{-1}_{2}({\mathbb{R}}^{d})\right)  $.

Fix $\varepsilon>0$ and let $Q=\{q_{k}, \; k\geq1\}$ be a sequence
so that, for all $x,y\in{\mathbb{R}}^{d}$,
\[
2\nu|y|^{2}-\sum_{k\geq1}q_{k}^{2}\sigma_{k}^{i}(x)\sigma_{k}^{j}(x)y_{i}y_{j}
\geq\varepsilon|y|^{2}.
\]
Then, for every $T>0$, there exits a unique $w(H^{1}_{2}({\mathbb{R}}%
^{d}),H^{-1}_{2}({\mathbb{R}}^{d}))$ Wiener Chaos solution
$\theta$ of equation
\begin{equation}
\label{eq:ps2.7}d{\theta}(t,x)=(\nu\Delta\theta(t,x) +f(t,x))dt-
\sigma _{k}(x)\cdot\nabla\theta(t,x) dw_{k}(t),
\end{equation}
The solution is an ${\mathcal{F}}_{t}$-adapted process and
satisfies
\[%
\begin{split}
\|\theta\|_{L_{2,Q}({\mathbb{W}}; L_{2}((0,T);H^{1}_{2}({\mathbb{R}}^{d}%
)))}^{2}+ \|\theta\|_{L_{2,Q}({\mathbb{W}}; {\mathbf{C}}((0,T); L_{2}%
({\mathbb{R}}^{d})))}^{2}\\
\leq C(\nu,q,T)\left(
\|\theta_{0}\|_{L_{2}({\mathbb{R}}^{d})}^{2} + \|f\|_{L_{2}((0,T);
H_{2}^{-1}({\mathbb{R}}^{d}))}^{2}\right)  .
\end{split}
\]

\end{theorem}

If $\max_{i}\sup_{x} |\sigma_{k}^{i}(x)|\leq C_{k}$, $k\geq1$,
then a possible choice of $Q$ is
\[
q_{k}= (\delta\nu)^{1/2}/(d2^{k}C_{k}),\ 0<\delta<2.
\]

If $\sigma_{k}^{i}(x)\sigma_{k}^{j}(x) \leq C_{\sigma}<+\infty$,
$i,j=1,\ldots,d$, $x\in{\mathbb{R}}^{d}$, then a possible choice
of $Q$ is
\[
q_{k}= \varepsilon\left(  2\nu/(C_{\sigma} d)\right)  ^{1/2},\
0<\varepsilon< 1.
\]

\section{Stochastic Navier-Stokes Equation}
\label{sec:.:SNS} \setcounter{equation}{0} \setcounter{theorem}{0}

In this section, we review the main facts about the stochastic
Navier-Stokes equation and indicate how the Wiener Chaos approach
can be used in the study of non-linear equations. Most of the
results of this section come from the two papers \cite{MR2} and
\cite{MR_AP}.

A priori, it is not clear in what sense the motion described by
Kraichnan's velocity (see Section \ref{secPS}) might fit into the
paradigm of Newtonian mechanics. Accordingly, relating the
Kraichnan velocity field $\mathbf{v}$ to classic fluid mechanics
naturally leads to the question whether we can compensate
$\mathbf{v}\left(  t,x\right)  $ by a field $\mathbf{u}\left(
t,x\right)  $ that is more regular with respect to the time
variable, so that there is a balance of momentum for the resulting
field $\mathbf{U}\left( t,x\right)  =\mathbf{u}\left(  t,x\right)
+\mathbf{v}\left(  t,x\right)  $ or, equivalently, that the motion
of a fluid particle in the velocity field $\mathbf{U}\left(
t,x\right)  $\ \ satisfies the Second Law of Newton.

A positive answer to this question is given in \cite{MR2}, where
it is shown that the equation for the smooth component
$\mathbf{u}=(u^{1},\ldots,u^{d})$ of the velocity is given by
\begin{equation}
\left\{
\begin{array}
[c]{l}%
du^{i}=[\nu\Delta u^{i}-u^{j}D_{j}u^{i}-D_{i}P+f_{i}]dt\\
\\
+\left(
g_{k}^{i}-D_{i}\tilde{P}_{k}-D_{j}\sigma_{k}^{j}u^{i}\right)
dw_{k},\ i=1,\ldots,d,\ 0<t\leq T;\\
\\
\text{ div\thinspace}\mathbf{u}=0,\text{
}\mathbf{u}(0,x)=\mathbf{u}_{0}(x).
\end{array}
\right.  \label{eq:sns0}%
\end{equation}
where $w_{k},\ k\geq1$ are independent standard Wiener processes
on a stochastic basis ${\mathbb{F}}$, the functions
$\sigma_{k}^{j}$ are given by (\ref{eq:v}), the known functions
$\mathbf{f}=(f^{1},\ldots,f^{d})$, $\mathbf{g}_{k}=(g_{k}^{i})$,
$i=i,\ldots,d$, $k\geq1$ are, respectively, the drift and the
diffusion components of the free force, and the unknown functions
$P$, $\tilde{P}_{k}$ are the drift and diffusion components of the
pressure.

\begin{remark}
It is useful to study equation (\ref{eq:sns0}) for more general
coefficients $\sigma_{k}^{j}.$ So, in the future, $\sigma_{k}^{j}$
are not necessarily the same as in Section \ref{secPS}.
\end{remark}

We make the following assumptions:

\begin{enumerate}
\item[{\textbf{NS1}}\ ] The functions
$\sigma^{i}_{k}=\sigma^{i}_{k}(t,x)$ are deterministic and
measurable,
\[
\sum_{k\geq1}\left(  \sum_{i=1}^{d} |\sigma_{k}^{i}(t,x)|^{2} +|D_{i}%
\sigma^{i}_{k}(t,x)|^{2}\right)  \leq K,
\]
and there exists $\varepsilon>0$ so that, for all
$y\in{\mathbb{R}}^{d}$,
\[
\nu|y|^{2}-\frac{1}{2}\sigma^{i}_{k}(t,x)\sigma^{j}_{k}(t,x)y_{i}y_{j}
\geq\varepsilon|y|^{2},
\]
$t\in[0,T]$, $x\in{\mathbb{R}}^{d}$.

\item[{\textbf{NS2}}\ ] The functions $f^{i},g_{k}^{i}$ are
non-random and
\[
\sum_{i=1}^{d}\left(  \Vert f^{i}\Vert_{L_{2}((0,T);H_{2}^{-1}({\mathbb{R}%
}^{d}))}^{2}+\sum_{k\geq1}\Vert g_{k}^{i}\Vert_{L_{2}((0,T);L_{2}({\mathbb{R}%
}^{d}))}^{2}\right)  <\infty.
\]
\end{enumerate}

\begin{remark}
In \textbf{NS1, }the derivatives $D_{i}\sigma_{k}^{i}$ are
understood as Schwartz distributions, but it is assumed that
$div\mathbf{\sigma:=}\sum
_{i=1}^{d}\mathbf{\partial}_{i}\sigma^{i}$ is a bounded
$l_{2}-$valued function. Obviously, the latter assumption holds in
the important case when
$\sum_{i=1}^{d}\mathbf{\partial}_{i}\sigma^{i}=0.$
\end{remark}

Our next step is to use the divergence-free property of
$\mathbf{u}$ to eliminate the pressure $P$ and $\tilde{P}$ from
equation (\ref{eq:sns0}). For that, we need the decomposition of
$L_{2}({\mathbb{R}}^{d};{\mathbb{R}}^{d})$ into potential and
solenoidal components.

Write $\mathfrak{S}({L}_{2}({\mathbb{R}}^{d};{\mathbb{R}}^{d}))=\{\mathbf{V}%
\in{L}_{2}({\mathbb{R}}^{d};{\mathbb{R}}^{d}): \mathrm{div}\;{\mathbf{V}%
}=0\}.$ It is known (see e.g. \cite{Kato}) that%
\[
L_{2}({\mathbb{R}}^{d};{\mathbb{R}}^{d})=\mathfrak{G}(L_{2}({\mathbb{R}}%
^{d};{\mathbb{R}}^{d}))\oplus\mathfrak{S}({L}_{2}({\mathbb{R}}^{d}%
;{\mathbb{R}}^{d})),
\]
where $\mathfrak{G}({L}_{2}({\mathbb{R}}^{d};{\mathbb{R}}^{d}))$
is a Hilbert
subspace orthogonal to $\mathfrak{S}({L}_{2}({\mathbb{R}}^{d};{\mathbb{R}}%
^{d}))$.

The functions $\mathfrak{G}(\mathbf{V})$ and
$\mathfrak{S}(\mathbf{V})$ can be
defined for $\mathbf{V}$ from any Sobolev space $H_{2}^{\gamma}({\mathbb{R}%
}^{d};{\mathbb{R}}^{d})$ and are usually referred to as the
potential and the divergence free (or solenoidal), projections,
respectively, of the vector field $\mathbf{V}$.

Now let $\mathbf{u}$ be a solution of equation (\ref{eq:sns0}).
Since $\mathrm{div}\;\mathbf{u}=0,$ we have
\[
D_{i}(\nu\Delta u^{i}- u^{j}D_{j}u^{i}- D_{i}P+f^{i}) =0;\
D_{i}(\sigma ^{j}_{k}D_{j}u^{j}u^{i} +g^{i}_{k}
-D_{i}\tilde{P}_{k}) =0,\ k\geq1.
\]
As a result,
\[
D_{i}P=\mathfrak{G}(\nu\Delta u^{i}-u^{j}D_{j}u^{i}+f^{i});\
D_{i}\tilde
{P}_{k}=\mathfrak{G}(\sigma_{k}^{j}D_{j}u^{i}+g_{k}^{i}),\
i=1,\ldots ,d,\ k\geq1.
\]
So, instead of equation (\ref{eq:sns0})$,$ we can and will
consider its equivalent form for the unknown vector
$\mathbf{u}=(u^{1}, \ldots, u^{d})$:
\begin{equation}
\label{eq:navs0}d{\mathbf{u}} =\mathfrak{S}(\nu\Delta{\mathbf{u}}-u^{j}%
D_{j}\mathbf{u}+\mathbf{f})dt
+\mathfrak{S}(\sigma^{j}_{k}D_{j}\mathbf{u}
+\mathbf{g}_{k})dw_{k}, \ 0<t\leq T,
\end{equation}
with initial condition $\mathbf{u}|_{t=0}=\mathbf{u}_{0}$.

\begin{definition}
\label{def:navs} An ${\mathcal{F}}_{t}$-adapted random process
$\mathbf{u}$ from the space
$L_{2}(\Omega\times[0,T];H^{1}_{2}({\mathbb{R}}^{d};{\mathbb{R}}^{d}))$
is called a solution of equation (\ref{eq:navs0}) if

\begin{enumerate}
\item With probability one, the process $\mathbf{u}$ is weakly
continuous in $L_{2}({\mathbb{R}}^{d};{\mathbb{R}}^{d})$.

\item For every $\mathbf{\varphi} \in \bC^{\infty}_{0}({\mathbb{R}}%
^{d},{\mathbb{R}}^{d})$, with $\mathrm{{div}\;\mathbf{\varphi}=0}$
there exists a measurable set $\Omega^{\prime}\subset\Omega$ so
that, for all $t\in[0,T]$, the equality
\begin{equation}
\label{eq:def-navs}%
\begin{split}
(u^{i},\varphi^{i})(t)  &  = ({u}^{i}_{0},\varphi^{i})+
\int_{0}^{t}\big( (\nu
D_{j}u^{i}, D_{j}\varphi^{i})(s)+\langle f^{i},\varphi^{i}\rangle(s) \big)ds\\
&  \int_{0}^{t}\big( \sigma_{k}^{j}D_{j}u^{i}+g^{i},\varphi^{i})
dw_{k}(s)
\end{split}
\end{equation}
holds on $\Omega^{\prime}$. In (\ref{eq:def-navs}), $(\cdot,
\cdot)$ is the inner product in $L_{2}({\mathbb{R}}^{d})$ and
$\langle\cdot, \cdot, \rangle$
is the duality between $H^{1}_{2}({\mathbb{R}}^{d})$ and $H^{-1}%
_{2}({\mathbb{R}}^{d})$.
\end{enumerate}
\end{definition}

The following existence and uniqueness result is proved in
\cite{MR_AP}.

\begin{theorem}
\label{th:navs-ex} In addition to \textbf{NS1} and \textbf{NS2},
assume that the initial condition $\mathbf{u}_{0}$ is non-random
and belongs to $L_{2}({\mathbb{R}}^{d};{\mathbb{R}}^{d})$. Then
there exist a stochastic basis
${\mathbb{F}}=(\Omega,\mathcal{F},\{\mathcal{F}_{t}\}_{t\geq
0},\mathbb{P})$ with the usual assumptions, a collection
$\{w_{k},k\geq1\}$ of independent standard Wiener processes on
${\mathbb{F}}$, and a process $\mathbf{u}$ so that $\mathbf{u}$ is
a solution of (\ref{eq:navs0}) and
\[
\mathbb{E}\left(  \sup_{s\leq T}\|\mathbf{u}(s)\|_{L_{2}({\mathbb{R}}%
^{d};{\mathbb{R}}^{d})}^{2}+\int_{0}^{T}\|\nabla\mathbf{u}(s)\|_{L_{2}%
({\mathbb{R}}^{d};{\mathbb{R}}^{d})}^{2}\,ds\right)  <\infty.
\]

If, in addition, $d=2$, then the solution of (\ref{eq:navs0})
exists on any
prescribed stochastic basis, is strongly continuous in $t$, is ${\mathcal{F}%
}^{W}_{t}$-adapted, and is unique, both path-wise and in
distribution.
\end{theorem}

When $d\geq3,$ existence of a strong solution as well as
uniqueness (strong or weak) for equation (\ref{eq:navs0}) are
important open problems.

By the Cameron-Martin theorem,
\[
{\mathbf{u}}(t,x)=\sum_{\alpha\in{\mathcal{J}}}\mathbf{u}_{\alpha}%
(t,x)\xi_{\alpha}.
\]
If the solution of (\ref{eq:navs0}) is
${\mathcal{F}}^{W}_{t}$-adapted, then, using the It\^{o} formula
together with relation (\ref{eq:xi(t)}) for the time evolution of
${\mathbb{E}}(\xi_{\alpha}|{\mathcal{F}}^{W}_{t})$ and relation
(\ref{zproduct}) for the product of two elements of the
Cameron-Martin basis, we can derive the propagator system for
coefficients $\mathbf{u}_{\alpha}$ \cite[Theorem 3.2]{MR_AP}:

\begin{theorem}
\label{th:navs-prop} In addition to \textbf{NS1} and \textbf{NS2},
assume that $\mathbf{u}_{0}\in
L_{2}({\mathbb{R}}^{d};{\mathbb{R}}^{d})$ and equation
(\ref{eq:navs0}) has an ${\mathcal{F}}^{W}_{t}$-adapted solution
$\mathbf{u}$ so that
\begin{equation}
\sup_{t\leq T}\mathbb{E}\|\mathbf{u}\|_{L_{2}({\mathbb{R}}^{d};{\mathbb{R}%
}^{d})}^{2}(t)<\infty. \label{ineqm}%
\end{equation}
Then
\begin{equation}
\mathbf{u}\left(  t, x\right)
=\sum_{\alpha\in\mathcal{J}}\mathbf{u}_{\alpha
}\left(  t,x\right)  \xi_{\alpha}, \label{hermfur}%
\end{equation}
and the Hermite-Fourier coefficients $\mathbf{u}_{\alpha}(t,x)$
are $L_{2}({\mathbb{R}}^{d};{\mathbb{R}}^{d})$-valued weakly
continuous functions so that
\begin{equation}
\sup_{t\leq T}\sum_{\alpha\in\mathcal{J}} \|\mathbf{u}_{\alpha}\|_{L_{2}%
({\mathbb{R}}^{d};{\mathbb{R}}^{d})}^{2}(t)+\int_{0}^{T}\sum_{\alpha
\in\mathcal{J}} \|\nabla\mathbf{u}_{\alpha}\|_{L_{2}({\mathbb{R}}^{d}%
;{\mathbb{R}}^{d\times d})}^{2}(t)\,dt<\infty. \label{estal}%
\end{equation}
The functions $\mathbf{u}_{\alpha}\left(  t,x\right)
,\alpha\in\mathcal{J}, $ satisfy the (nonlinear) propagator
\begin{equation}%
\begin{split}
\label{eq:fura}\frac{\partial}{\partial t}\mathbf{u}_{\alpha}  &
=\mathfrak{S}\Big(\Delta\mathbf{u}_{\alpha} -\sum_{\gamma,\beta\in
\Delta\left(  \alpha\right)  }\Psi\left(
\alpha,\beta,\gamma\right)  \left(
\mathbf{u}_{\gamma},\nabla\mathbf{u}_{\beta}\right)  +I_{\left\{
\left\vert
\alpha\right\vert =0\right\}  }\mathbf{f}\\
\\
&  +\sum_{j,k}\sqrt{\alpha_{j}^{k}}\left(  \left(  \mathbf{\sigma}^{k}%
,\nabla\right)  \mathbf{u}_{\alpha^{-}\left(  j,k\right)  }
+I_{\left\{ \left\vert \alpha\right\vert =1\right\}
}\mathbf{g}^{k}\right)  m_{j}\left(
t\right)  \Big),\ 0<t\leq T;\\
\mathbf{u}_{\alpha}|_{t=0}  &  =\mathbf{u}_{0}I_{\left\{
\left\vert \alpha\right\vert =0\right\}  };
\end{split}
\end{equation}
recall that the numbers $\Psi(\alpha, \beta, \gamma)$ are defined
in (\ref{eq:Psi}).
\end{theorem}

One of the questions in the theory of the Navier-Stokes equation
is computation of the mean value
${\bar{\mathbf{u}}}={\mathbb{E}}\mathbf{u}$ of the solution. The
traditional approach relies on the Reynolds equation for the mean
\begin{equation}%
\begin{array}
[c]{l}%
\partial_{t}{\mathbf{\bar{u}}}-\nu\Delta{{\bar{\mathbf{u}}}}+\overline{\left(
\mathbf{\ u},\nabla\right)  \mathbf{\ u}}=0,
\end{array}
\label{reynolds}%
\end{equation}
which is not really an equation with respect to
$\mathbf{\bar{u}}$. Decoupling $\left(  \ref{reynolds}\right)  $
has been an area of active research: Reynolds approximations,
coupled equations for the moments, Gaussian closures, and so on
(see e.g. \cite{Monin}, \cite{VF} and the references therein)

Another way to compute $\bar{\mathbf{u}}\left(  t,x\right)  $ is
to find the distribution of $\mathbf{v}\left(  t,x\right)  $ using
the infinite-dimensional Kolmogorov equation associated with
(\ref{eq:navs0}). The complexity of this Kolmogorov equation is
prohibitive for any realistic application, at least for now.

The propagator provides a third way: expressing the mean and other
statistical moments of $\mathbf{u}$ in terms of
$\mathbf{u}_{\alpha}$. Indeed, by Cameron-Martin Theorem,
\begin{align*}
{\mathbb{E}}\mathbf{u}(t,x)  &  =\mathbf{u}_{0}(t,x),\\
{\mathbb{E}} u^{i}(t,x)u^{^{_{j}}}(s,y)  &  =\sum_{\alpha\in \cJ}u_{\alpha}%
^{i}(t,x)u_{\alpha}^{j}(s,y)
\end{align*}
If exist, the third- and fourth-order moments can be computed
using (\ref{eq:third}) and (\ref{eq:fourth}).

The next theorem, proved in \cite{MR_AP}, shows that the existence
of a solution of the propagator (\ref{eq:fura}) is not only
necessary but, to some extent, sufficient for the global existence
of a probabilistically strong solution of the stochastic
Navier-Stokes equation (\ref{eq:navs0}).

\begin{theorem}
\label{th:global} Let \textbf{NS1} and \textbf{NS2} hold and $\mathbf{u}%
_{0}\in L_{2}({\mathbb{R}}^{d};{\mathbb{R}}^{d})$. Assume that the
propagator $\left(  \ref{eq:fura}\right)  $ has a solution
$\left\{  \mathbf{u}_{\alpha }\left(  t,x\right)
,\;\alpha\in\mathcal{J}\right\}  $ on the interval $(0,T]$ so
that, for every $\alpha$, the process $\mathbf{u}_{\alpha} $ is
weakly continuous in $L_{2}({\mathbb{R}}^{d};{\mathbb{R}}^{d})$
and the inequality
\begin{equation}
\sup_{t\leq T}\sum_{\alpha\in\mathcal{J}} \|\mathbf{u}_{\alpha}\|_{L_{2}%
({\mathbb{R}}^{d};{\mathbb{R}}^{d})}^{2}(t)
+\int_{0}^{T}\sum_{\alpha
\in\mathcal{J}} \|\nabla\mathbf{u}_{\alpha}\|_{L_{2}({\mathbb{R}}%
^{d};{\mathbb{R}}^{d\times d})}^{2}(t)\,dt<\infty\label{estal1}%
\end{equation}
holds. If the process
\begin{equation}
\mathbf{\bar{U}}\left(  t,x\right)  :=\sum_{\alpha\in\mathcal{J}}%
\mathbf{u}_{\alpha}\left(  t,x\right)  \xi_{\alpha}\text{ }%
\end{equation}
is $\mathcal{F}_{t}^{W}$-adapted, then it is a solution of
(\ref{eq:navs0}).

The process $\mathbf{\bar{U}}$ satisfies
\[
\mathbb{E}\left(  \sup_{s\leq T}\|\mathbf{\bar{U}}(s)\|_{L_{2}({\mathbb{R}}%
^{d};{\mathbb{R}}^{d})}^{2}+\int_{0}%
^{T}\|\nabla\mathbf{\bar{U}}(s)\|_{L_{2}({\mathbb{R}}%
^{d};{\mathbb{R}}^{d\times d})}^{2}\,ds\right)  <\infty
\]
and, for every
$\mathbf{v\in{L}_{2}({\mathbb{R}}^{d};{\mathbb{R}}^{d})},$
${\mathbb{E}}\left(  \mathbf{\bar{U}},\mathbf{v}\right)  $ is a
continuous function of $t$.
\end{theorem}

Since $\mathbf{\bar{U}}$ is constructed on a prescribed stochastic
basis and over a prescribed time interval $\left[  0,T\right]  $,
this solution of (\ref{eq:navs0}) is strong in the probabilistic
sense and is global in time. Being true in any space dimension
$d$, Theorem \ref{th:global} suggests another possible way to
study equation (\ref{eq:navs0}) when $d\geq3$. Unlike the
propagator for the linear equation, the system (\ref{eq:fura}) is
not lower-triangular and not solvable by induction, so that
analysis of (\ref{eq:fura}) is an open problem.

 \section{First-Order It\^{o} Equations}
 \label{sec:.:frd}
\setcounter{equation}{0}
\setcounter{theorem}{0}

The objective of this section is to study  equation
\begin{equation}
\label{eq:nclassical1} du(t,x)= u_x(t,x)dw(t),\ t>0,\ x \in \bR,
\end{equation}
and its analog for $x\in \bR^d$.

 Equation (\ref{eq:nclassical1})  was first encountered in Example \ref{ex2};
 see also \cite{Gikhman}.
 With a non-random  initial condition  $u(0,x)=\varphi(x)$,
 direct computations show that, if exists, the Fourier transform $\hat{u}=\hat{u}(t,y)$
 of the solution must satisfy
 \begin{equation}
 \label{eq718}
 d\hat{u}(t,y)=\sqrt{-1}y\hat{u}(t,y)dw(t), \ {\rm or}\
 \hat{u}(t,y)=\hat{\varphi}(y)e^{\sqrt{-1}yw(t)+\frac{1}{2}y^2t}.
 \end{equation}
 The last equality
  shows that the properties of the solution essentially depend on the
 initial condition, and, in general, the solution is not in $L_2(\mbW)$.

The S-transformed equation, $v_t=h(t)v_x$, has a unique solution
$$
v(t,x)=\varphi\left(x+\int_0^th(s)ds\right), \
h(t)=\sum_{i=1}^Nh_im_i(t).
$$
The results of Section \ref{sec:.:WN} imply that a white noise
solution of the equation can exist only if $\varphi$ is a real
analytic function. On the other hand,
 if $\varphi$ is infinitely differentiable, then, by Theorem \ref{th:SW111},
the Wiener Chaos solution exists and can be recovered from $v$.

\begin{theorem}
\label{th:f-order} Assume  that the initial condition $\varphi$
belongs to the Schwarz space $\cS=\cS(\bR)$ of tempered
distributions. Then there exists a generalized random process
$u=u(t,x)$, $t\geq 0$, $x \in \bR$, so that, for every  $\gamma\in
\bR$ and  $T>0$,  the process $u$ is the unique
$w(H^{\gamma}_2(\bR), H^{\gamma-1}_2(\bR))$ Wiener Chaos solution
of equation (\ref{eq:nclassical1}).
\end{theorem}

\begin{proof} The  propagator for (\ref{eq:nclassical1})  is
\begin{equation}
\label{eq410} u_{\alpha}(t,x)=\varphi(x)I(|\alpha|=0)+\int_0^t
\sum_{i} \sqrt{\alpha_i}(u_{\alpha^-(i)}(s,x))_xm_i(s)ds.
\end{equation}
Even though  Theorem \ref{th:soft} is not applicable, the system
can be solved by induction if $\varphi$ is sufficiently smooth.
 Denote by $C_{\varphi}(k)$, $k\geq 0$,
the square of the $L_2(\bR)$ norm of the $k^{{\rm th}}$ derivative
of $\varphi$:
\begin{equation}
\label{eq:ex3}
C_{\varphi}(k)=\int_{-\infty}^{+\infty}|\varphi^{(k)}(x)|^2dx.
\end{equation}
By  Corollary \ref{cor:ind},  for every $k\geq 0$ and $n\geq 0$,
\begin{equation}
\label{eq411} \sum_{|\alpha|=k}
\|(u^{(n)}_{\alpha})_x\|^2_{L_2(\bR)}(t) = \frac{t^k
C_{\varphi}(n+k)}{k!}.
\end{equation}
The statement of the theorem now follows.\end{proof}

\begin{remark}
Once interpreted in a suitable sense, the Wiener Chaos solution of
(\ref{eq:nclassical1}) is $\cF^W_t$-adapted and does not depend on
the choice of the Cameron-Martin basis in $L_2(\mbW)$. Indeed,
choose the wight sequence so that
$$
r_{\a}^2=\frac{1}{1+C_{\varphi}(|\a|)}.
$$
By (\ref{eq411}), we have  $u\in \cR L_2(\mbW; L_2(\bR))$.

Next, define
$$
\psi_N(x)=\frac{1}{\pi}\frac{\sin(Nx)}{x}.
$$
Direct computations show that the Fourier transform of $\psi_N$ is
supported in $[-N,N]$ and $\int_{\bR}\psi_N(x)dx=1$. Consider
equation  (\ref{eq:nclassical1}) with initial condition
$$
\varphi_{N}(x)=\int_{\bR}\varphi(x-y)\psi_N(y)dy.
$$
By (\ref{eq718}), this equation has a unique solution $u_N$ so
that  $u_N(t,\cdot)\in L_2(\mbW; H^{\gamma}_2(\bR))$, $t\geq 0$,
$\gamma \in \bR$. Relation (\ref{eq411}) and the definition of
$u_N$ imply
$$
\lim_{N\to \infty} \sum_{|\a|=k}
\|u_{\a}-u_{N,\a}\|_{L_2(\bR)}^2(t)=0,\ t\geq 0,\ k \geq 0,
$$
so that, by the Lebesgue dominated convergence theorem,
$$
\lim_{N\to \infty} \|u - u_N\|_{\cR L_2(\mbW;L_2(\bR))}^2(t)=0, \
t\geq 0.
$$
In other words,  the solution of the propagator (\ref{eq410})
corresponding to any basis $\mathfrak{m}$ in $L_2((0,T))$ is a
limit in $\cR L_2(\mbW;L_2(\bR))$ of the sequence $\{u_N, \; N\geq
1\}$ of $\cF^W_t$-adapted processes.
\end{remark}

The properties  of the Wiener Chaos solution of
(\ref{eq:nclassical1}) depend on the
 growth rate of the numbers  $C_{\varphi}(n)$. In particular,
\begin{itemize}
\item  If $C_{\varphi}(n) \leq  C^n(n!)^{\gamma}, \ C>0,\ 0\leq
\gamma < 1,$
then \\
$ u \in L_2\left(\mbW; L_2((0,T); H^n_2(\bR))\right)$ for all
$T>0$ and every $n\geq 0$. \item If $C_{\varphi}(n) \leq C^nn!, \
C>0,$  then
\begin{itemize}
\item for every $n\geq 0$, there is a $T>0$ so that
 $ u \in L_2\left(\mbW; L_2((0,T);H^n_2(\bR))\right)$.
  In other words, the square-integrable
 solution exists only for sufficiently small $T$.
 \item for every $n\geq 0$ and every $T>0$, there exists a number
 $\delta\in (0,1)$ so that
 $u \in L_{2,Q}\left(\mbW; L_2((0,T);H^n_2(\bR))\right)$ with $Q=(\delta, \delta,
 \delta,\ldots)$.
 \end{itemize}
 \item If
 the numbers $C_{\varphi}(n)$ grow as $C^n(n!)^{1+\rho}$, $\rho\geq 0$,
 then, for every $T>0$, there exists a number  $\gamma>0$ so that \\
   $u \in (\cS)_{-\rho,-\gamma}\left(L_2(\mbW); L_2((0,T);H^n_2(\bR))\right)$.
   If $\rho>0$, then this solution does not belong to any
   $L_{2,Q}\left(\mbW; L_2((0,T);H^n_2(\bR))\right)$.
   If $\rho>1$, then this solution does not have an S-transform.
\item If
 the numbers $C_{\varphi}(n)$ grow faster than  $C^n(n!)^b$ for any $b,C>0$, then
 the Wiener Chaos solution of (\ref{eq:nclassical1})  does not belong to
 any \\
 $(\cS)_{-\rho,-\gamma}\left( L_2((0,T);H^n_2(\bR))\right)$,
 $\rho,\;\gamma>0$, or
 $L_{2,Q}\left(\mbW; L_2((0,T);H^n_2(\bR))\right)$.
 \end{itemize}

 To construct a function $\varphi$  with the required rate of growth of
 $C_{\varphi}(n)$,  consider
  $$
  \varphi(x)=\int_0^{\infty}\cos(xy)e^{-g(y)}dy,
  $$
  where $g$ is a suitable positive, unbounded, even function. Note that,  up to a
  multiplicative constant, the Fourier transform of $\varphi$ is $e^{-g(y)}$, and
  so $C_{\varphi}(n)$ grows with $n$ as $\int_0^{+\infty}|y|^{2n}e^{-2g(y)}dy$.

A more general first-order equation can be considered:
\begin{equation}
\label{eq:nclassSP} du(t,x)=\sigma_{ik}(t,x)D_iu(t,x)dw_k(t), \
t>0,\ x \in \bR^d.
\end{equation}

\begin{theorem}
Assume that in equation (\ref{eq:nclassSP}) the initial condition
 $u(0,x)$ belongs to $\cS(\bR^d)$  and each
$\sigma_{ik}$ is infinitely differentiable with respect to $x$ so
that $\sup_{(t,x)} |D^n\sigma_{ik}(t,x)| \leq C_{ik}(n)$, $n\geq
0$. Then there exists a generalized random process $u=u(t,x)$,
$t\geq 0$, $x \in \bR^d$, so that, for every  $\gamma\in \bR$ and
$T>0$,  the process $u$ is the unique $w(H^{\gamma}_2(\bR^d),
H^{\gamma-1}_2(\bR^d))$ Wiener Chaos solution of equation
(\ref{eq:nclassical1}).
\end{theorem}

\begin{proof} The arguments are identical to the proof of Theorem
\ref{th:f-order}.
\end{proof}

Note that the S-transformed equation (\ref{eq:nclassSP}) is
 $v_t=h_k\sigma_{ik}D_iv$ and  has a unique
solution if each $\sigma_{ik}$ is a Lipschitz continuous function
of $x$. Still, without additional smoothness, it is impossible to
relate this  solution to any generalized random process.


\providecommand{\bysame}{\leavevmode\hbox
to3em{\hrulefill}\thinspace}
\providecommand{\MR}{\relax\ifhmode\unskip\space\fi MR }
\providecommand{\MRhref}[2]{%
  \href{http://www.ams.org/mathscinet-getitem?mr=#1}{#2}
} \providecommand{\href}[2]{#2}

\label{LR-REF}

\end{document}